\documentclass[11pt]{article}

\usepackage{yk}
\usepackage{amsmath,amsfonts,amssymb, dsfont, mathrsfs, comment, booktabs, numprint, accents, float,  hyperref, xr, breakurl, xurl, enumerate}
\usepackage[numbers]{natbib}
\usepackage{lscape}
\usepackage{pgfplotstable,filecontents}
\usepackage[margin=1.3in]{geometry}

\newcommand{\dbtilde}[1]{\accentset{\approx}{#1}}
\newcommand{\CPP}{\text{CPP}}
\newcommand{\s}{\text{sign}}
\newcommand{\init}{\text{init}}
\newcommand{\pen}{\text{pen}}
\newcommand{\new}{\text{new}}
\def\independenT#1#2{\mathrel{\rlap{$#1#2$}\mkern2mu{#1#2}}}
\newcommand{\sargmin}{\text{mid argmin}}
\newcommand\indep{\protect\mathpalette{\protect\independenT}{\perp}}
\title{On robust learning in the canonical change point problem under heavy tailed errors in finite and growing dimensions}

\author{Debarghya Mukherjee, Moulinath Banerjee and Ya'acov Ritov \\\\ Department of Statistics, University of Michigan, Ann Arbor, MI}
\begin{document}
\maketitle

\begin{abstract}%   <- trailing '%' for backward compatibility of .sty file
This paper presents a number of new findings about the canonical change point estimation problem. The first part studies the estimation of a change point on the real line in a simple stump model using the robust Huber estimating function which interpolates between the $\ell_1$ (absolute deviation) and $\ell_2$ (least squares) based criteria. While the $\ell_2$ criterion has been studied extensively, its robust counterparts and in particular, the $\ell_1$ minimization problem have not. We derive the limit distribution of the estimated change point under the Huber estimating function and compare it to that under the $\ell_2$ criterion. Theoretical and empirical studies indicate that it is more profitable to use the Huber estimating function (and in particular, the $\ell_1$ criterion) under heavy tailed errors as it leads to smaller asymptotic confidence intervals at the usual levels compared to the $\ell_2$ criterion. We also compare the $\ell_1$ and $\ell_2$ approaches in a parallel setting, where one has $m$ independent single change point problems and the goal is to control the maximal deviation of the estimated change points from the true values, and establish rigorously that the $\ell_1$ estimation criterion provides a superior rate of convergence to the $\ell_2$, and that this relative advantage is driven by the heaviness of the tail of the error distribution. Finally, we derive minimax optimal rates for the change plane estimation problem in growing dimensions and demonstrate that Huber estimation attains the optimal rate while the $\ell_2$ scheme produces a rate sub-optimal estimator for heavy tailed errors. In the process of deriving our results, we establish a number of properties about the minimizers of compound Binomial and compound Poisson processes which are of independent interest. 
%A comparison of the two procedures in multiple dimensions will require determination of the limit distributions which is a highly complex problem in its own right and will require independent investigation.
\end{abstract}

%\begin{keywords}
%change point estimation, heavy tailed error, robust learning
%\end{keywords}

\section{Introduction}

\label{sec:intro}
In the canonical change-point or change-boundary estimation problem, one posits a regression (or a classification) model in which the conditional distribution of the response given the covariate(s) changes from a constant value on one side of an unknown boundary in covariate space to another on the opposite side. Within the genre of regime change problems, the canonical model is a particularly convenient formulation for investigating the fundamentals of estimation and inference, and the challenges involved therein.  In particular, in the one-dimensional case, this gives us the so-called `stump model': 
$$
Y = \alpha_0 1(X \leq d_0) + \beta_0 1(X > d_0) + \eps 
$$ 
with $\alpha_0 \ne \beta_0$, where $X$ assumes values in $\mathbb R$. In the multidimensional scenario with a $p$-dimensional covariate $X$, a natural extension is given by 
$$
Y = \alpha_0 1(\psi(X, d_0) \leq 0) + \beta_0 1(\psi(X, d_0) > 0) + \epsilon \,,
$$ 
where $\psi(X, d_0) = 0$ defines a low dimensional smooth surface in $\mathbb R^p$. 
\\\\
\noindent
This paper deals with the estimation of change parameters in such models under different estimating functions in both fixed and growing dimensions along with the calibration of minimax optimal rates. The use of a variety of robust estimating functions is necessitated by the fact that heavy-tailed errors frequently drive data generating mechanisms associated with change-point problems in applications pertaining to finance (\cite{bradley2003financial}), hydrology (\cite{basso2015emergence}), climate and environmental science (\cite{weitzman2011fat}), internet data (\cite{hernandez2004variable}) and  genetics (\cite{slater2015power}).  We show in this paper that such robust criteria are essential for attaining optimal convergence rates when the number of parameters diverges with sample size. We also show that in the fixed dimension scenario the choice of the criterion function does not affect the convergence rate but \emph{does affect} the tails of the limit distribution of the estimated change-point in a way that makes the use of robust criteria more profitable for thick-tailed errors. 
%The simple linear regression model assumes a uniform linear relationship between the covariate and the response. In practice, the situation is often more complicated: for instance, several sub-population may have different $\beta_j$ parameter. Some common techniques to account for such heterogeneity includes mixed linear models, fitting different models among each sub-population which correspond to a supervised classification setting where the true groups are known. 
\\\\
\noindent
%The problem of studying heavy tailed random variables arises frequently in various fields of statistical applications, e.g. finance \cite{bradley2003financial}, hydrology \cite{basso2015emergence}, climate and environmental science \cite{weitzman2011fat}, internet data \cite{hernandez2004variable}, genetics \cite{slater2015power}.  Robust estimation is, of course, a well-developed area, and was pioneered by Huber in his seminal paper \cite{huber1964} (see, also  \cite{huber1992robust}), where the idea is to replace the squared error loss (that is highly sensitive to outliers) by less sensitive losses, e.g. the Huber loss. The use of such robust loss functions, including the simple $\ell_1$ loss where one considers the absolute distance between target and parameters has, since then, been widely adopted in various arenas of statistics. The excellent text by \cite{huber2004robust} provides an overview of the area.
We next focus on the organization of the manuscript and articulate the contributions of each section. But before that, we take a moment to introduce the (scaled) Huber estimating function (HEF) which is referred to below and used throughout the manuscript. The scaled HEF is defined as $\tilde H_k(x) := ((k+1)/k) H_k(x)$ where: 
$$
H_k(x) = 
\begin{cases}
\frac{x^2}{2} \,. & \text{ if } |x| \le k \\
k\left(|x| - \frac{k^2}{2}\right) \,, & \text{ otherwise } \,. 
\end{cases}
$$
The cost function corresponding to $\tilde{H}_k$ in a generic statistical problem can be written as $C_k(Z, \theta) := \tilde{H}_k(g(Z) - h(Z,\theta))$ where $g(Z)$ is some functional of the data vector $Z$ (say, the real-valued response in a regression model) and $h(Z,\theta)$ is some known function of $Z$ and the parameter $\theta$ (say, the regression function). Note that as $k \to 0$, we have $\tilde H_k(x) \to |x|$ and for $k \to \infty$, $\tilde H_k(x) = x^2/2$, therefore $C_k$ interpolates between the $\ell_1$ and $\ell_2$ cost functions via the parameter $k$. The function $H_k$  was introduced in the pioneering work of Peter Huber \cite{} for the robust estimation of parameters in presence of outliers. The key idea here is the observation that $\ell_1$ discrepancy is more robust to outliers than the $\ell_1$ discrepancy, whereas $\ell_2$ discrepancy has other attractive features like differentiability with constant curvature. The Huber function seeks to combine these two discrepancies and utilize the best of both worlds.
\\\\
\noindent
Section 2 presents a treatment of the canonical stump model with a one-dimensional covariate under HEF optimization as well as its limiting incarnations (the $\ell_1$ and the $\ell_2$ criteria) and provides explicit statements of asymptotic distributions which are seen to be the minimizers of various compound Poisson processes. While the limiting  behavior under $\ell_2$ has been long known in the literature, the study of the asymptotic properties under robust criteria is new. More interestingly, we are able to characterize the tail behaviors of the limit distributions in terms of the tail-indices of the corresponding error distributions which, to the best of our knowledge, was previously unknown.  We demonstrate that under the $\ell_2$ criterion the tail index of the error adversely affects the tail of the minimizer of the corresponding compound Poisson process: errors with polynomial decay of tails lead to polynomially decaying tails for the limit; while under HEF (including the $\ell_1$ criterion) the tail of the limit distribution is \emph{unaffected} by the tail of the error and is necessarily sub-exponential. This has direct implications for the construction of asymptotic confidence intervals as we discuss later. 
\\\\
\noindent
Section 3 explores the canonical problem for a growing number of change-point parameters. 
%indeed, it is the growing dimensional version of the problem that led to the current paper, with the one-dimensional aspect of the problem arising later. In the growing parameters scenario there are two natural versions. 
The first part pertains to situations where multiple change-point parameters are estimated in parallel from \emph{separate} data-sources, and this number is allowed to grow with the total sample size. The second version is the change-boundary problem alluded to at the beginning of our narrative. We explore, specifically, the case of a linear boundary, i.e. a model of the form
\[ Y = \alpha_0 \,\mathds 1(X^T d_0 \leq 0) + \beta_0 \mathds 1(X^T d_0 > 0) + \epsilon\,, \]
with $\|d_0\|_2 = 1$ (to enforce identifiability) and a $p$-dimensional covariate $X$. This is the so-called change-plane model which captures the core features of the change-boundary problem. Our motivation for studying change plane problems stems from the recent use of change-plane models in personalized medicine and related problems \cite{wei2015cox}, \cite{fan2017change}, as well as the use of change-plane models in econometrics (e.g. see \cite{seo2007smoothed}, \cite{li2018multi} and references therein). We assume that $n$ i.i.d. observations are available from this model and that either $p = o(n)$ or $p \gg n$ with the number of non-zero co-ordinates of $d_0$ constrained to be appropriately small. We show that in both the parallel change point and high dimensional change plane problems, the $\ell_2$ criterion based estimator  suffers from the curse of dimensionality unlike its robust counterparts. 
\\\\
Section 4 presents a range of simulation studies in the 1-dimensional case that compare the quantiles of the limit distributions under $\ell_1$ and $\ell_2$ criteria and discusses the observed patterns. Section 5 concludes, providing among other things an exposition of future challenges in this area. 

%{\noindent \em Remainder omitted in this sample. See http://www.jmlr.org/papers/ for full paper.}

% Acknowledgements should go at the end, before appendices and references

%\acks{We would like to acknowledge support for this project
%from the National Science Foundation (NSF grant IIS-9988642)
%and the Multidisciplinary Research Program of the Department
%of Defense (MURI N00014-00-1-0637). }


\section{Robust change point estimation in one dimension}
\label{sec:one_dim_findings}
\begin{figure}
\centering
\includegraphics[scale=0.25]{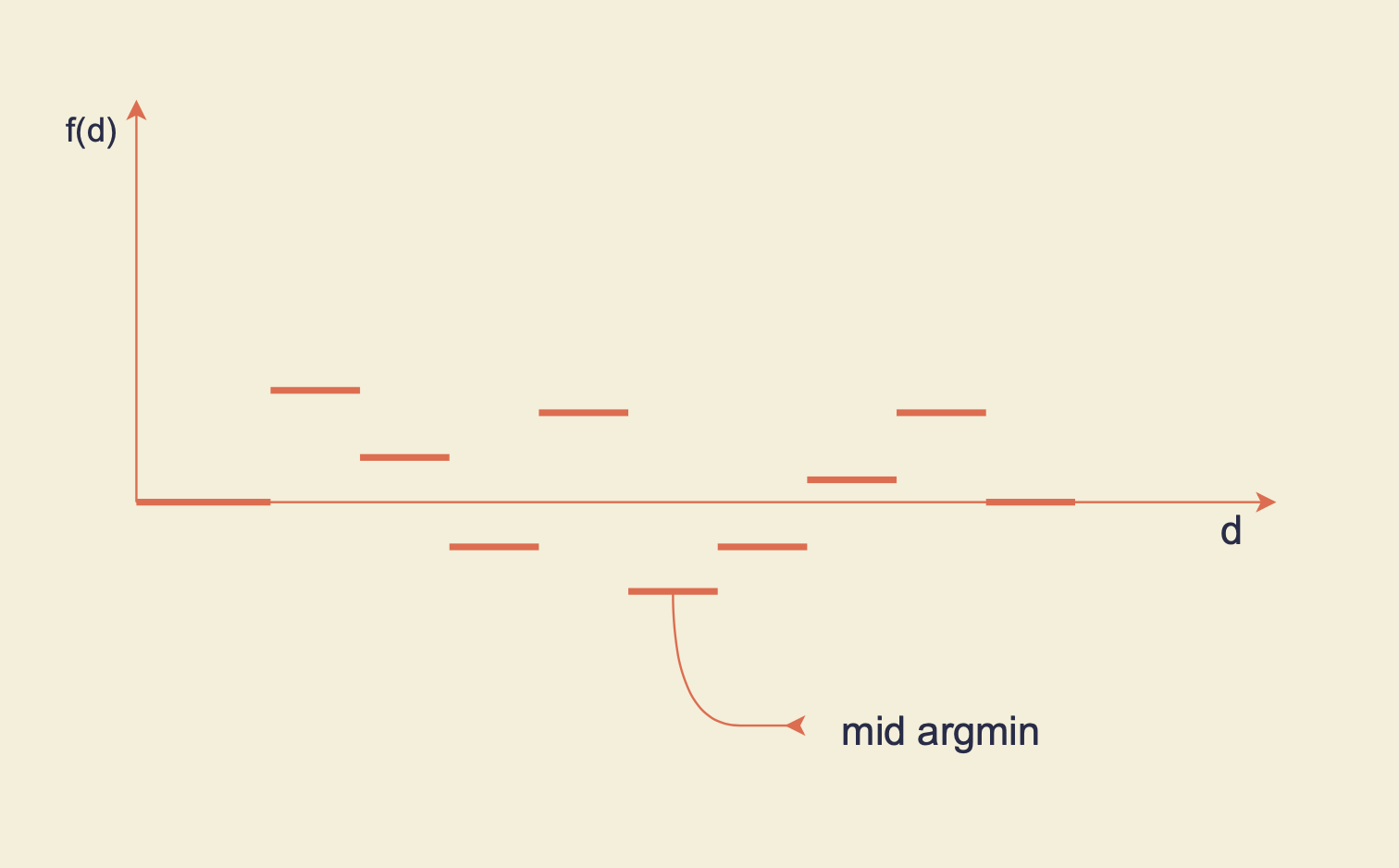}
\caption{The mid-argmin of a piecewise constant function}
\label{fig:mid_argmin}
\end{figure}

In this section we analyze the following canonical change point model in one dimension:
\begin{equation}
\label{eq:change_point_model_one_dim}
Y_i = \alpha_0 \mathds{1}_{X_i \le d_0} + \beta_0 \mathds{1}_{X_i > d_0} + \xi \,,
\end{equation}
for $1 \le i \le n$. The least squares estimators of the parameters of this model have been well-explored in the literature, but quite surprisingly, nothing is known about its  robust variant, and the trade-offs between the two approaches. To understand the difference, consider an even simpler model:  
$$
Y_i = \mathds{1}_{X_i > d_0} + \xi_i \,.
$$
where $X_i \in \reals$ is a real covariate and $\xi_i$ is a mean 0 error independent of $X_i$. Here $d_0$ is the parameter of interest, i.e. the \emph{change point} in the space of covariates. Traditionally, one minimizes the squared-error loss to obtain an estimator of $d_0$:
\begin{align*}
\hat d^{\ell_2} & = \sargmin_{d \in I} \frac{1}{n} \sum_{i=1}^n \left(Y_i - \mathds{1}_{X_i > d}\right)^2 \\
& =  \sargmin_{d \in I} \frac{1}{n} \sum_{i=1}^n\left(Y_i - \frac{1}{2}\right)\mathds{1}_{X_i \leq d} \\
& := \sargmin_{d \in I} \ f(d) \,.
\end{align*}
for some compact interval $I \subset \bbR$. Note that the function $f(d)$ is a right continuous step function with respect to $d$, therefore its minimizer is not unique, in fact it is an interval. By $\sargmin$, we denote the midpoint of the corresponding interval. (See Figure \ref{fig:mid_argmin} for an illustration.)
The statistical properties of this estimator are well-known; e.g. see Chapter 14 of \cite{kosorokintroduction} or Proposition 1 of \cite{lan2009change} and its preceding discussion. For example, if $X$ is compactly supported with density bounded away from $0$ and $\infty$ on its support, then:
$$
n\left(\hat d^{\ell_2} - d_0\right) \overset{\mathscr{L}}{\implies} \sargmin_{t \in \bbR} M(t)
$$
where $M(t)$ is a two-sided compound Poisson process with drift described thus: Let $N(t)$ be a homogeneous Poisson process with intensity parameter $f_X(d_0)$ on $[0, \infty)$ where $f_X(\cdot)$ is the density of $X$. Define two independent stochastic processes $V^+(t)$ on $[0, \infty)$ and $V^-(t)$ on $(-\infty, 0]$
 as follows: 
\begin{align*}
V^+(t) & = \sum_{i=1}^{N_1(t)}\left(\xi_i +\frac12\right) \\
V^-(t)  & = \sum_{i=1}^{N_2(-t)}\left(\xi_{-i} -\frac12\right)
\end{align*}
where $\left\{\xi_i\right\}_{i \in \bbZ \backslash \{0\}}$ are i.i.d. from the distribution of $\xi$ and $N_1(t), N_2(t)$ are i.i.d copies of $N(t)$, and are independent of the $\xi_i$'s. Then 
$$
M(t) = V^+(t)\mathds{1}_{t \ge 0} - V^-(t)\mathds{1}_{t < 0} \,,
$$ 
is a two sided compound Poisson process on the real line (we denote it by $CPP(\xi + 1/2, f_X(\theta_0))$ that drifts off to $\infty$ on either side, and is minimized almost surely on an interval of points. Taking the mid-argmin of this process ensures symmetry of the limiting distribution under the symmetry of the distribution of $\xi$. 
\\\\
\noindent
The asymptotics above require only a second moment for the errors and therefore are valid for many heavy-tailed errors. However, heavy tailed errors enlarge the spread of the limit distribution, resulting in wider confidence intervals for the change-point parameter. This is because the compound Poisson process is closely related to the two sided random walk on $\mathbb Z$ with step distribution given by $(\xi + 1/2)$ to the right of 0 and $(-\xi + 1/2)$ to its left. We quantify later in this section how the tail of the distribution of the minimizer of this compound Poisson process depends on the tail index of the error distribution with heavy tailed errors corresponding to a heavier tail for the minimizer which, in turn, implies a wider asymptotic confidence interval. 
\\\\
\noindent
The natural question, then, is what happens if one were to compute $d_0$ via  the robust HEF, in particular, say the $\ell_1$ criterion, and whether asymptotic efficiency relative to the $\ell_2$ criterion would accrue as a result in the case of heavy-tailed errors. So, consider:
$$
\hat d^{\ell_1} = \sargmin_{d \in I} \frac{1}{n} \sum_{i=1}^n \left|Y_i - \mathds{1}_{X_i \ge d}\right| \,.
$$
For consistency of $\hat d^{\ell_1}$ we need the assumption that $\med(\xi) = 0$. Since the $\ell_2$ criterion requires $\bbE(\xi) = 0$, \emph{the rest of the paper will be developed for symmetric errors} which simplifies the discussion without compromising conceptual issues. We show later (see Theorem \ref{cor:l1_loss}) that 
$$
n\left(\hat d^{\ell_1}- d_0\right) \overset{\mathscr{L}}{\implies} \sargmin_{t \in \bbR} M_R(t)
$$
where $M_R(t)$ is, again, a two sided compound Poisson process with intensity parameter $f_X(\theta_0)$ and the step-distribution given by that of $\left| \eps + 1\right| - \left| \eps \right|$. Observe that the random variable $\left| \eps + 1\right| - \left| \eps \right|$ is bounded in absolute value by $1$ irrespective of the tail index of the error and consequently sub-gaussian. This translates to a sub-exponential tail for the asymptotic distribution, resulting in a tighter asymptotic confidence interval than the one obtained via minimizing squared error loss. 
\\\\
\noindent
We next present our main results for more general stump model described in equation \eqref{eq:change_point_model_one_dim}.  %For consistency of the estimator, we henceforth assume that the distribution of $\xi$ is symmetric around 0. 
Minimizing the Huber estimating function yields the following estimator: 
\begin{equation}
\label{eq:huber_estimators}
\left(\hat \alpha^k, \hat \beta^k, \hat d^k\right) = \sargmin_{\alpha, \beta, d} \frac1n \sum_{i=1}^n \tilde H_k\left(Y_i - \alpha \mathds{1}_{X_i \le d} - \beta \mathds{1}_{X_i > d}\right)
\end{equation}
where, as mentioned earlier, we consider the midpoint of the minimizing interval of $d_0$. We next present the asymptotic distributions of $(\hat \alpha^k, \hat \beta^k, \hat d^k)$ upon proper centering and scaling. 
\begin{comment} 
We start with a lemma which is imperative to the proof of our main theorem, but may be of independent interest: 
\begin{lemma}
\label{lem:huber_lower_bound}
If $\xi$ follows a symmetric distribution around the origin with with continuous density $f_{\xi}$ satisfying $f_{\xi}(0) > 0$, then for any $k > 0, |\mu| < k$, we have: 
$$
\bbE\left[\tilde H_k(\xi + \mu) - \tilde H_k(\xi)\right] \ge \frac{\mu^2}{2}\bbP\left(-k \le \xi \le k- \mu \right) \ge \frac{\mu^2}{2}\bbP\left(-k \le \xi \le 0 \right) \,.
$$
\end{lemma}
\end{comment} 
%We now state our one of the main theorems, which quantifies the asymptotic distribution of $(\hat \alpha^k, \hat \beta^k, \hat d^k)$: 
\begin{theorem}
\label{thm:asymptotic_huber}
Suppose $\theta_0 \in I$ for some compact interval $I \subset \bbR^3$. Assume that the density of $X$ is continuous and strictly positive at $d_0$. Then the estimators $(\hat \alpha^k, \hat \beta^k, \hat d^k)$ obtained in equation \eqref{eq:huber_estimators} are asymptotically independent and satisy: 
\begin{align*}
\sqrt{n}(\hat \alpha^k - \alpha_0) & \overset{\mathscr{L}}{\implies}\cN\left(0, \frac{\sigma_k^2}{\mu_k^2F_X(d_0)}\right) \,, \\
\sqrt{n}(\hat \beta^k - \beta_0) & \overset{\mathscr{L}}{\implies} \cN\left(0, \frac{\sigma_k^2}{\mu_k^2\bar F_X(d_0)}\right) \,, \\
n(\hat d^k - d_0) &  \overset{\mathscr{L}}{\implies} \sargmin_{t \in \bbR}  \CPP\left(\tilde H_k\left(\xi + |\alpha_0 - \beta_0|\right) - \tilde H_k(\xi), f_X(\theta_0)\right)
\end{align*}
where the parameters $\mu_k$ and $\sigma_k$ are: 
\begin{align*}
\mu_k & = \frac{k+1}{k}\bbP\left(-k \le \xi \le k \right) \\
\sigma_k^2 & =\left( \frac{k+1}{k}\right)^2\left( \bbE\left[\xi^2\mathds{1}_{-k \le \xi \le k}\right] + 2k^2 \bbP\left(\xi > k\right)\right)  \,.
\end{align*}
\end{theorem}
\noindent
%\begin{remark}
Note that if $k \to 0$, then $\mu_k \to \mu^{\ell_1} = 2f_{\xi}(0)$ and $\sigma_k^2 \to (\sigma^{\ell_1})^2 = 1$. One the other hand, if $k \to \infty$, then $\mu_k \to \mu^{\ell_2} = 1$ and $\sigma_k^2 \to (\sigma^{\ell_2})^2 = \sigma_{\xi}^2$, as long as $E(\xi^2)$ is finite, which is a requirement for the $\ell_2$ based estimation strategy to work. The following two theorems present the asymptotic distribution of the estimated parameters (upon proper centering and scaling) for these special cases: $\ell_1$ and $\ell_2$ criteria, where we see that the limiting parameters are indeed $\mu^{\ell_1}, \sigma^{\ell_1}$ and $\mu^{\ell_2}, \sigma^{\ell_2}$ respectively.  We note that the proofs do not directly follow by taking the limit of $k$ in the proof ofTheorem \ref{thm:asymptotic_huber}, but rely on similar techniques.  %employed in the proof of Theorem \ref{thm:asymptotic_huber}. 
%\end{remark}

\begin{theorem}
\label{cor:l1_loss}
Consider minimizing the $\ell_1$ criterion function to obtain:   
$$
\hat \theta^{\ell_1} = \argmin_{\theta \in I} \frac1n \sum_{i=1}^n \left|Y_i - \alpha\mathds{1}_{X_i \le d} - \beta \mathds{1}_{X_i > d}\right|
$$
Then, under the assumptions of Theorem \ref{thm:asymptotic_huber}, we have:  
\begin{align*}
\sqrt{n}(\hat \alpha^{\ell_1} - \alpha_0) & \overset{\mathscr{L}}{\implies}\cN\left(0, \frac{1}{4f_{\xi}^2(0)F_X(d_0)}\right) \,, \\
\sqrt{n}(\hat \beta^{\ell_1} - \beta_0) & \overset{\mathscr{L}}{\implies} \cN\left(0, \frac{1}{4f_{\xi}^2(0)\bar F_X(d_0)}\right) \,, \\
n(\hat d^{\ell_1} - d_0) &  \overset{\mathscr{L}}{\implies} \sargmin_{t \in \bbR}  \CPP\left(|\xi + |\alpha_0 - \beta_0|| - |\xi|, f_X(\theta_0)\right) \,,
\end{align*}
and the estimates of the parameters are asymptotically independent. 
\end{theorem}

\begin{theorem}
\label{cor:l2_loss}
Consider minimizing the $\ell_2$ criterion function to obtain: 
$$
\hat \theta^{\ell_2} = \argmin_{\theta} \frac1n \sum_{i=1}^n \left(Y_i - \alpha\mathds{1}_{X_i \le d} - \beta \mathds{1}_{X_i > d}\right)^2
$$
Then, under the assumptions of Theorem \ref{thm:asymptotic_huber}, we obtain: 
\begin{align*}
\sqrt{n}(\hat \alpha^{\ell_2} - \alpha_0) & \overset{\mathscr{L}}{\implies}\cN\left(0, \frac{\sigma^2_{\xi}}{F_X(d_0)}\right) \,, \\
\sqrt{n}(\hat \beta^{\ell_2} - \beta_0) & \overset{\mathscr{L}}{\implies} \cN\left(0, \frac{\sigma^2_{\xi}}{\bar F_X(d_0)}\right) \,, \\
n(\hat d^{\ell_2} - \theta_0) &  \overset{\mathscr{L}}{\implies} \sargmin_{t \in \bbR}  \CPP\left(\xi + \frac{|\alpha_0 - \beta_0|}{2}, f_X(\theta_0)\right) \,,
\end{align*}
and the estimates of the parameters are asymptotically independent. 
\end{theorem}
\noindent
Observe from the above results that the asymptotic distributions of $\sqrt{n}(\hat \alpha - \alpha_0)$ and $\sqrt{n}(\hat \beta - \beta_0)$ are normal irrespective of the estimating function used for estimation, but the asymptotic variance depends upon the estimating function. 
%One may use the optimal $k$ which minimizes the variances if one is interested in efficient estimation of the mean levels $(\alpha_0, \beta_0)$. 
\\\\
\noindent
The more interesting part is how the asymptotic distributions of $n(\hat d - d_0)$ changes from the $\ell_1$ to the $\ell_2$ estimating function. In either case, the asymptotic distribution is characterized as the minimizer of a compound Poisson process, but the step-size is sensitive to the criterion. This has a bearing on the tail-behavior of the minimizer when $\xi$ is heavy-tailed as articulated below in Theorem \ref{thm:tail_bound_cpp}. 
%Observe that, irrespective of the distribution of the error $\xi$, the steps of the compound poisson process in Theorem \ref{cor:l1_loss} is bounded and hence sub-gaussian. This suggests that the limiting distribution of $n(\hat d - d_0)$ under $\ell_1$ loss is more concentration around $0$, especially under heavy-tailed error, as the tail index of $\xi$ does not affect the tail behavior of bounded random variable $|\xi + (\alpha_0 - \beta_0)| - |\xi|$. 
% prove in Theorem \ref{thm:tail_bound_cpp}.  
For notational simplicity, define $F_{\ell_i}$ as the limiting distribution of $n(\hat d - d_0)$ under the $\ell_i$ estimating function for $i \in \{1, 2\}$: 
\begin{align*}
F_{\ell_1}(x) & = \bbP\left(\sargmin_{t \in \bbR}  \CPP\left(|\xi + |\alpha_0 - \beta_0|| - |\xi|, f_X(\theta_0)\right) \le x\right) \,. \\
F_{\ell_2}(x) & = \bbP\left(\sargmin_{t \in \bbR}  \CPP\left(\xi + \frac{|\alpha_0 - \beta_0|}{2}, f_X(\theta_0) \right) \le x\right) \,.
\end{align*}
As we are working with the mid argmin, both $F_{\ell_1}$ and $F_{\ell_2}$ are symmetric around $0$ [e.g. see the discussion in Section 4.2 of \cite{lan2009change}].  
\newline
\newline
To compare the tail properties of $F_{\ell_1}$ and $F_{\ell_2}$ in presence of heavy tailed error, we assume the following distribution of $\xi$ in our subsequent analysis: 
\begin{equation}
\label{eq:xi_dist}
\bbP\left(|\xi| > x\right) = \frac{1}{1 + x^{\gamma}}
\end{equation}
and $\xi$ is symmetric around $0$. This ensures that $\bbE[|\xi|^{\gamma - \nu}] < \infty$ for all $0 < \nu \le \gamma$. %We start with a lower bound on the partial sums of independent random variables which may of independent interest:  
\begin{comment} 
\end{comment} 
We next present a theorem which quantifies the tails of the asymptotic distribution of $n(\hat d - d_0)$ under the $\ell_1$ and $\ell_2$ estimating functions for the above heavy-tailed errors.  
\begin{theorem}
\label{thm:tail_bound_cpp}
In our change point model equation \eqref{eq:change_point_model_one_dim}, under the error distribution specified in equation \eqref{eq:xi_dist}, we have for all $x \ge k_0$: 
$$
\bar F_{\ell_2}(x)  = 1- F_{\ell_2}(x) \ge \frac{c_0}{2f^{\gamma}_X(d_0)}x^{-\gamma} \,.
$$
for some constants $k_0, c_0, \mu_0$ explicitly mentioned in the proof. On the other hand, we have for all $x > 0$: 
$$
\bar F_{\ell_1}(x) = 1- F_{\ell_1}(x)  \le \frac{p^*}{e^{\frac{\mu_0^2}{8(\alpha_0 - \beta_0)^2}} - 1}\exp{\left(-xf_X\left(d_0\right)\left(1 - e^{-\frac{\mu_0^2}{8(\alpha_0 - \beta_0)^2}}\right)\right)} \,.
$$
where $p^* = \bbP\left(\min_{1 \le i < \infty} \sum_{j=1}^i \left(\left|\xi_j + |\alpha_0 - \beta_0|\right| - |\xi_j| \right) > 0\right) > 0$. 
\end{theorem}
\noindent
From Theorem \ref{thm:tail_bound_cpp}, it is immediate that the asymptotic distribution of $n(\hat d^{\ell_2} - d_0)$ is affected by the tail index of the error distribution of $\xi$: it can not decay faster than $x^{-\gamma}$, whereas the asymptotic distribution of of $n(\hat d^{\ell_1} - d_0)$ has a sub-exponential tail\footnote{Some of the constants involved in the sub-exponential tail bound of course depend on the distribution of $\xi$.}. Therefore for all large $x$, we have: 
$$
\bbP\left(-x \le \cD_{\ell_2}  \le x\right) \le \bbP\left(-x \le \cD_{\ell_1} \le x\right) \,.
$$
where $\cD_{\ell_2}$ (resp. $\cD_{\ell_1}$) is the limit of $n(\hat d^{\ell_2} - d_0)$ (resp. $n(\hat d^{\ell_1} - d_0)$). Therefore, it is preferable to use the change point estimator $\hat d^{\ell_1}$ to $\hat d^{\ell_2}$ for constructing an asymptotic confidence interval for all large enough levels of confidence. \emph{To the best of our knowledge, this is the first result characterizing the tail behavior of limiting compound Poisson processes, and can be expected to be of independent interest}. More detailed empirical comparisons are presented in Section \ref{sec:simulation}.

\noindent
\paragraph{Proof idea: } We now present a brief sketch of the proof of Theorem \ref{thm:tail_bound_cpp}. In Theorems \ref{cor:l1_loss} and \ref{cor:l2_loss}, we established the limiting distribution of $n(\hat d^{\ell_1} - d_0)$ and $n(\hat d^{\ell_2} - d_0)$ respectively. Both distributions are compound Poisson processes but with different step distributions: for the limit of the $\ell_2$ estimator, the step distribution is $\xi + 1/2$ and for the $\ell_1$ estimator, the step distribution is $|\xi + |\alpha_0 - \beta_0|| - |\xi|$. Hence, if $\xi$ is heavy-tailed (resp. light tailed), so is the step distribution of the limit of the $\ell_2$ estimator, whereas the steps of the limit of the $\ell_1$ estimator are bounded (and therefore sub-gaussian) irrespective of the tail of $\xi$. As a compound Poisson process is closely related to the random walk corresponding to its step-size, we first establish that the tail of the minimizer of the random walk depends on that of the error distribution. In particular, in Lemmas \ref{lem:lower_random_onesided} and \ref{lem:two_sided_rw_minimizer}, we show that if $\xi$ has a power tail structure, i.e. $\bbP(|\xi| > t) \sim t^{-\gamma}$ for some $\gamma > 0$, then the tail of the minimizer of the random walk is also lower bounded by $x^{-\gamma}$. This lower bound can be translated to a lower bound on the minimizer of the compound Poisson process. On the other hand, for the limit distribution of the $\ell_1$ estimator of the change point, the step distribution is sub-gaussian. Therefore, we first establish an exponential upper bound on the tail of the minimizer of a random walk with bounded steps and use it to obtain an exponential tail bound for a compound Poisson process with bounded steps. Details of the proof of Theorem \ref{thm:tail_bound_cpp} can be found in Appendix.

\begin{remark}
\label{rem:power_tail}
Although we have assumed a specific distribution for $\xi$ to establish our results, an inspection of the proofs shows that the only fact essential to the calculations is the \emph{power tail} structure of $\xi$, i.e. $\bbP(|\xi| > x) \sim x^{-\gamma}$ for some $\gamma > 0$. Our assumed functional form simply facilitates some routine computations and can be easily extended to the more general case. Therefore, the first conclusion of Theorem \ref{thm:tail_bound_cpp} is valid as long as $\xi$ has power tail with index $\gamma$. The second conclusion of Theorem \ref{thm:tail_bound_cpp} is agnostic to the tail index of $\xi$ and continues to hold for any $\xi$,  as long as it has finite variance. In fact, the broad conclusions of the above theorem are true for any Huber estimating function $\tilde H_k$ for $0 \le k < \infty$: any such Huber function based estimate yields a sub-exponential tail for the limiting minimizer. 
%Hence, the advantages of the limiting distribution of the change point estimator obtained by minimizing the Huber estimating function $\tilde H_k$ for $0 \le k < \infty$ (and in particular the $\ell_1$ criterion) over the estimator obtained by minimizing the $\ell_2$ criterion are also valid for any error distribution with this power tail structure. 
\end{remark}

\begin{remark}
\label{rem:subgaussian_errors}
%In Theorem \ref{thm:tail_bound_cpp}, we have established a lower bound on the tail of the limiting distribution of the change point estimator obtained by minimizing squared error loss in presence of heavy tailed error, whereas the $\ell_1$-estimator is ignorant of the tail index of the error. 
By using similar arguments to the proof of the above theorem, we can show that for a sub-gaussian $\xi$, both $\ell_1$ and $\ell_2$ criteria yield the sub-exponential concentration bound. Therefore, there is no significant gain in using robust criteria in comparison to the $\ell_2$ criterion in the presence of sub-gaussian errors. 
\end{remark}

\section{Estimation in multidimensional change-problems}
In the previous section, we have seen that with one-dimensional change point estimation, the advantage of using the more robust $\ell_1$ estimating function is expected to confer efficiency in terms of the spread of the limiting distribution (i.e. the length of the asymptotic confidence interval), but the rate of convergence is invariant to the estimating function used. In fact, this rate can be shown to be minimax optimal, i.e. one cannot get a better rate without any further assumptions. However, the effect of using a robust estimating function is more striking when the number of change points to be estimated grows with increasing sample size. \\\\
\noindent
In this section, we present two scenarios: one with many one-dimensional change points and the other with a high dimensional change-boundary, in both of which we estimate a diverging number of parameters and establish that it is possible to achieve faster rates of convergence in these situations in the presence of heavy-tailed errors using robust criteria, and in particular, the $\ell_1$ criterion.  
%In Subsection \ref{sec:parallel_change_point}, we analyze the \emph{free problem} which is a parallel change point problem, which is connected to distributed computing setup. In Subsection \ref{sec:findings_multi_dim}, we study a change plane problem, a multi-dimensional incarnation of the change point problem analyzed in Section \ref{sec:one_dim_findings}, where the covariate $X$ is assumed to be high dimensional (dimension may increase slowly with the sample size, or may even be larger than the sample size) and the parameter of interest is a hyperplane boundary. 

\subsection{Parallel change point estimation}
\label{sec:parallel_change_point}
Suppose we have $m$ parallel processes of one-dimensional change point models, with each process having $n$ independent observations. Specifically, the $i^{th}$ process carries $n$ pairs of covariate-response pairs from the following model:
$$
Y_{i, j} =  \mathds{1}_{X_{i, j} > d_{0, i}} + \xi_{i, j} \,,
$$
for $1 \le j \le n$ and $1 \le i \le m$.  Here, as before, we assume that $\{(X_{i,j}, \xi_{i, j})\}_{i,j}$ are i.i.d., $\xi_{i,j} \indep X_{i, j}$ and $\xi_{i, j}$ is symmetric around 0. Furthermore, we assume that all $nm$ pairs of observations are independent. The $\left\{d_{0, i}\right\}_{i=1}^m$'s are \emph{free parameters} to be estimated from the data. Due to the independence among the samples, $d_{0, i}$ is estimated only from the $n$ observations for the $i^{th}$ problem. Define $\hat d_{i, 1}$ and $\hat d_{i, 2}$ to be the smallest argmin estimators obtained for the $i$'th problem by minimizing the $\ell_1$ and $\ell_2$ criteria respectively. We would like to control the estimation errors across the different problems simultaneously, hence the natural metric to consider is the maximal loss over the $m$ problems. Specifically, we want to quantify the order of 
$$
\max_{1 \le i \le m} \left|\hat d_{i, k} - d_{0, i}\right|, \ \ k = 1, 2 \,. 
$$
We prove below that, for an appropriate growth rate of $n$ relative to $m$, the maximal error only inherits the slow factor $\log{m}$ for the robust estimators (i.e. $\{\hat d_{i, 1}\}_{i=1, \dots, p}$) irrespective of the tail of the error, whereas a factor of $m^{1/\gamma}$ in unavoidable with the $\ell_2$ estimates when $P(|\xi| \ge t) \sim t^{-\gamma}$. 
%More precisely we obtain: 
%\begin{align*}
%\max_{1 \le i \le m} \left|\hat d_{i, 1} - d_{0, i}\right| & = O_p\left(\frac{\log{m}}{n}\right) \\
%\max_{1 \le i \le m} \left|\hat d_{i, 2} - d_{0, i}\right| & = O_p\left(\frac{m^{1/\gamma}}{n}\right)\,.
%\end{align*}
%Furthermore, we establish that the factor $m^{1/\gamma}$ arising in the maximum $\ell_2$ estimating function is a \emph{tight} factor in the sense that the maximal $\ell_2$ deviation divided by $m^{1/\gamma}/n$ is not asymptotically degenerate at 0 as $n \rightarrow \infty$. Hence, the use of  the $\ell_1$ estimating function improves the rate of convergence, the implications being that one can obtain tighter uniform confidence bands for a growing number of change-points with $\ell_1$ based estimating function. %In contrast, \emph{had the errors been subgaussian}, we would have acquired the same $\log{m}$ factor for the $\ell_2$ loss as well.
\\\\
\noindent
We now present our theorem. As before, the distribution of $\xi$ is assumed to be symmetric and $|\xi|$ is distributed as: 
$$
\bbP\left(|\xi| \ge t\right) = \frac{1}{1 + t^{\gamma}} \,,
$$
for all $t \ge 0$. Echoing Remark \ref{rem:power_tail}, the core arguments of our proof only require the power tail structure of $\xi$, i.e. $\bbP(|\xi| > t) \sim t^{-\gamma}$. The following theorem highlights the disparity between the rates of convergence of the maximal deviations of the $\ell_2$ and $\ell_1$ based estimators. 
\begin{theorem}
\label{thm:parallel_change_point}
Suppose the change point estimator $\hat d^{\ell_2}_{i}$ for the $i^{th}$ problem is obtained by minimizing the squared error loss. If $n/m^{1/\gamma} \to \infty$, then for any $t > 0$: 
$$
\liminf_{n \to \infty} \bbP\left(\max_{1 \le i \le m}\frac{n}{m^{1/\gamma}}\left|\hat d^{\ell_2}_{i} - d_{0, i}\right| > t \right) \ge c(t) > 0 \,,
$$
where $c(t)$ is some positive constant depending on $t$ and other model parameters. On the other hand, if we obtain $\hat d^{\ell_1}_{i}$ by minimizing the $\ell_1$ estimating function, then we have: 
$$
\bbP\left(\frac{n}{\log{m}}\max_{1 \le i \le m}\left|\hat d^{\ell_1}_{i} - d_{0, i}\right| > t\right)  \le \frac{2e^{-c}}{1 - e^{-c}} e^{-\log{m}\left(t\frac{f_X(d_0)}{2}(1 - e^{-c}) - 1\right)}\,,
$$
as long as $n/\log{m} \to 0$ for some constant $c$ explicitly mentioned in the proof. 
\end{theorem}
\noindent
{\bf Proof idea: }We document the main ideas behind the proof of Theorem \ref{thm:parallel_change_point}.  The first part of the theorem establishes a lower bound on the rate of convergence of the $\ell_\infty$ error of all estimated change points across all the problems. The main idea of the proof is that, for a finite sample of size $n$, the distribution of $n(\hat d^{\ell_2}_i - d_{i, 0})$ (for any $1 \le i \le m$), is given by that of the minimizer of a \emph{compound binomial process} defined as below: 
$$
n\left(\hat d^{\ell_2}_i - d_{0,i}\right) \overset{d}{=} \sargmin_t \sum_{k=1}^{N_{n,+}^{i}(t)} \left(\varepsilon_k + \frac12\right)\mathds{1}_{t \ge 0} + \sum_{k=1}^{N_{n,-}^{i}(t)} \left(\tilde{\varepsilon}_k + \frac12\right)\mathds{1}_{t < 0} \,,
$$
where the binomial processes $N_{n, +}$ and $N_{n, -}$ are defined as: 
\begin{align*}
N_{n, +}^{i}(t) = \sum_{i=1}^n \mathds{1}_{d_{0,i} \le X_{i,j} \le d_{0,i} + \frac{t}{n}} & \sim \text{Bin} \left(n, F_X\left(d_{0, i} + \frac{t}{n}\right) - F_X(t)\right) \ \ \forall \ \ t > 0 \,, \\
N_{n, -}^{i}(t) = \sum_{i=1}^n \mathds{1}_{d_{0,i} + \frac{t}{n}\le X_{i,j} \le d_{0,i}} & \sim \text{Bin}  \left(n, F_X(t) - F_X\left(d_{0, i} + \frac{t}{n}\right)\right) \ \ \forall \ \ t < 0 \,,
\end{align*}
and then $\{\varepsilon_k\}$ are the $\xi_{i,j}$'s corresponding to the $X_{i,j}$'s satisfying $d_0 \le X_{i,j} \le d_0 + \frac{t}{n}$ and the $\{\tilde{\varepsilon}_k\}$ are the $-\xi_{i,j}$'s corresponding to the $X_{i,j}$'s satisfying $d_0 + \frac{t}{n}\le X_{i,j} \le d_0$. 
\\\\
\noindent
%Observe that, these processes are finite sample approximation of compound Poisson process, supported on $\{0, 1, \dots, n\}$. 
The distribution of $n(\hat d^{\ell_2}_i - d_{i, 0})$ is closely related to a random walk with step distribution $\xi + 1/2$, where the number of steps is derived from the binomial processes. Therefore, we first establish a lower bound on the tail of the minimizer of the random walk and then translate that lower bound to the tail of the distribution of $n|\hat d^{\ell_2}_i - d_{0,i}|$ (see Lemma \ref{lem:finite_sample_tail_bound}). Finally, we use the fact for any set of independent random variables $Z_1, \dots, Z_m$: 
$$
\bbP\left(\max_{1 \le i \le m} Z_i > t\right) = 1 - \Pi_{i=1}^m F_{Z_i}(t) = 1 - \Pi_{i=1}^m \left(1 - \bbP(Z_i > t)\right) \,.
$$
Hence, any lower bound on the tail of $Z_i$ yields a lower bound on the tail of $\max_{1 \le i \le m} Z_i$. Taking $Z_i = n|\hat d^{\ell_2}_i - d_{i, 0}|$ and converting the lower bound on the tail of $n|\hat d^{\ell_2}_i - d_{i, 0}|$ to the tail of $\max_{1 \le i \le m} n|\hat d^{\ell_2}_i - d_{i, 0}|$ concludes the first part of the proof. 
\\\\
\noindent
The proof of the second part is similar to the first, where instead of the lower bound we establish an upper bound on the tail of the $n|\hat d^{\ell_1}_i - d_{0,i}|$. Note that, in case of $\ell_1$ criterion: 
$$
n\left(\hat d^{\ell_1}_i - d_{0,i}\right) \overset{d}{=} \sargmin_t \sum_{i=1}^{N_{n,+}(t)} \left(\left|\xi_i + 1\right| - |\xi_i|\right)\mathds{1}_{t \ge 0} + \sum_{i=1}^{N_{n,-}(t)} \left(\left|\xi_i + 1\right| - |\xi_i|\right)\mathds{1}_{t < 0}
$$ 
The steps now are uniformly bounded and therefore sub-gaussian. Following the same line of arguments as in the first part of the proof, we first establish an upper bound on the tail of the minimizer of the random walk with bounded steps which is then translated to an upper bound on the tail of the $n|\hat d^{\ell_1}_i - d_{i, 0}|$ (see Lemma \ref{lem:finite_sample_tail_bound_upper}) and finally to the tail of $\max_{1 \le i \le m} n|\hat d^{\ell_1}_i - d_{0,i}|$ using a union bound. The detailed proof can be found in Appendix. 
\\

\begin{remark}
The above theorem shows the detrimental effect of the $\ell_2$ estimating funciton under heavy-tailed errors owing to the growing number of estimated parameters. The $\ell_1$ based estimator is only marginally affected (by the $\log m$ factor). While we don't establish this in the paper, the HEF based estimator used in the previous section will also yield the same rate of convergence as the $\ell_1$ based estimator. Further, the results are easily generalizable to the generic stump model with unknown levels on either side of the change-point with some standard technical modifications to our current proof.

%the rate of convergence of the maximal deviation of the least squares estimates becomes slower, whilst a much faster rate can be obtained via $\ell_1$ loss function. The basic idea is same as before: the distribution of the change point estimator is governed by a minimizer of a compound Poisson process, a close cousin of random walk. The distribution of the steps of this process are bounded (and hence subgaussian) when the estimators are obtained via minizing $\ell_1$ loss function, whilst, it is unbounded for $\ell_2$ loss. 
%Theorem \ref{thm:parallel_change_point} into two parts: the first one establishes the lower bound under $\ell_2$ loss which ensures one cannot obtain faster rate using the squared error loss without further assumptions, and the other for the tightness under $\ell_1$ loss. %
\end{remark}

\noindent
Finally we show that the rate obtained above (i.e. $n/\log{m}$) cannot be improved in general, even in the case of the zero error situation, i.e. this rate is minimax optimal, provided that we don't have any background information about the spread of the change points $\left\{d_{0, i}\right\}_{1 \le i \le m}$. 

\begin{theorem}
\label{thm:mlb_parallel}
Consider the above scenario of $m$ independent change point problems where for the $i$'th problem the observations are generated from the following stump model: 
$$
Y_{i, j} = \mathds{1}_{X_{i, j} > d_{0, i}} + \xi_{i, j} \,.
$$
Denote by $P_{d_{0 ,i}}$, the joint distribution of $(X, Y)$ or equivalently $(X, \xi)$ of the observations in $i^{th}$ problem which satisfies the conditions $\xi \indep X$ and $\xi$ has symmetric distribution around origin. Then we have: 
$$
\liminf_{n \to \infty} \frac{n}{\log{m}}\inf_{\left\{\hat d_i\right\}_{1 \le i \le m}} \sup_{\otimes_{i=1}^m P_{d_{i_0}}} \bbE\left[\max_{1 \le i \le m} \left|\hat d_i - d_{i, 0}\right|\right] \ge C > 0 \,,
$$
for some universal constant $C$. 
\end{theorem}

\subsection{Estimation of a change plane in growing dimensions}
\label{sec:findings_multi_dim}
As described in Section \ref{sec:intro}, a multi-dimensional version of the canonical stump model is the so-called `change-plane' problem:
\begin{equation}
\label{eq:change_point_model_multi_dim}
Y_i = \alpha_0\mathds{1}_{X_i^{\top}d_0 \le 0} + \beta_0 \mathds{1}_{X_i^{\top}d_0 > 0} + \xi_i \,.
\end{equation}
where $X_i, d_0 \in \reals^p$ and $p$ is assumed growing with $n$.  As $d_0$ is only identifiable up to its scale, we assume $d_0 \in S^{p-1}$. As before, we assume that $\left\{(X_i, \xi_i)\right\}_{i=1}^n$ are i.i.d and that $\xi_i$ is independent of $X_i$ with a symmetric distribution around the origin. We analyze the above canonical change plane model in two regimes: (i) when $p/n \to 0$ (Subsection \ref{sec:slowly_growing}) and (ii) when $p \gg n$ (Subsection \ref{sec:fast_growing}). In both regimes, the dimension of $d_0$ is increasing with sample size, but with one fundamental difference: when $p/n \to 0$, we have many more samples than parameters and should therefore be able to estimate $d_0$ consistently, whereas when $p \gg n$, the problem is ill-posed and as is customary in the high dimensional literature, we need to impose a sparsity assumption on $d_0$: i.e. an upper bound on the number of its non-null entries. 
%in a sense that all the parameters can not be estimated consistency without further assumptions. As we know from the existing literature of the standard high dimensional analysis of statistical model, the main philosophy here is to assume that the true parameter $d_0$ is sparse, i.e. most of its co-ordinates are 0 and number of non-zero elements are though growing with $n$, but much less than $n$. 
Mathematically speaking, we assume that $\|d_0\|_0 \le s_0$ for some unknown $s_0$ which satisfies $\left(s_0\log{p}\right)/n \to 0$. Our aim is to recover the non-zero signals in $d_0$ consistently. 
%Note that the problem would boil down to the regime of $p/n \to 0$ if we knew the support, i.e. which $s_0$ elements are of $d_0$ are non-zero. But here we know neither $s_0$, not its support, only the existence of some $s_0$ which satisfying assumption \ref{assm:sparsity}. Therefore, we resort to the Structural Risk Minimization (henceforth SRM) principle (see e.g. \cite{}, \cite{} and references therein) for model selection purpose. The idea of SRM is \textbf{should we reword the discussion of SRM from Manski paper?} 
We show that the rate of convergence of the change plane estimator obtained by minimizing the HEF (apart from $k = \infty$, i.e. the squared error loss) is minimax optimal in both the scenarios and is independent of the tail of the error distribution, whereas the $\ell_2$ criterion based analysis yields a slower convergence rate for  heavy tailed errors, which depends on the tail index of the error. 
%In this paper, our interest in the growing dimensional change-plane problem is primarily theoretical, though in future work (not within the scope of this grant) we would like to investigate the computational angles.
%\\\\
%\noindent
%
%
%
%
%
%
%
%In this section, we present our findings for estimating change plane in growing dimensions regime as described in equation \eqref{eq:change_point_model_multi_dim} of Section \ref{sec:intro}: 
%$$
%Y_i = \alpha_0 \mathds{1}_{X_i^{\top}d_0 \le 0} + \beta_0 \mathds{1}_{X_i^{\top}d_0 > 0} + \xi_i \,.
%$$
%Let $p$ be the underlying dimension of the covariate $X$. We analyze two different regimes in this section: (1) when $p$ increases with $n$ but slowly, i.e. $p/n \to 0$ and (2) when $p \gg n$. 

\subsubsection{When $d/n\to0$}
\label{sec:slowly_growing}
%In this section we analyze the following change model introduced in equation \eqref{eq:change_point_model_multi_dim} in Section \ref{}: 
%$$
%Y_i = \alpha_0 \mathds{1}_{X_i^{\top}d_0 \le 0} + \beta_0 \mathds{1}_{X_i^{\top}d_0 > 0} + \xi \,,
%$$
%for $1 \le i \le n$.
%{\bf ELABORATE ON THE ESTIMATION PROCEDURE, ONE VS TWO STAGES:} 
In the change plane estimation problem, we consider the semi-metric:  
\begin{align*}
& \dist\left((\alpha_1, \beta_1, d_1), (\alpha_2, \beta_2, d_2)\right) \\
& = \sqrt{(\alpha_1 - \alpha_2)^2 + (\beta_1 - \beta_2)^2 + \bbP\left(\s(X^{\top}d_1) \neq \s(X^{\top}d_2)\right)}\,,
\end{align*}
which is motivated by the one used by (\cite{kosorok2007introduction}, see Chapter 14) with the only difference being that instead of considering the Euclidean distance between two candidate change-plane vectors $d_1$ and $d_2$ we use the mass of the wedge bounded by the two corresponding corresponding hyperplanes to define a metric. This particular metric is geometrically convenient to analyze in the change-plane problem as will be seen in our subsequent computations and can be easily related to the $\ell_2$ distance under an additional condition which is satisfied under various distributional assumptions on the covariate $X$. 
%Using this semi-metric we established that: 
%$$
%\sqrt{n} \ \dist\left((\hat \alpha, \hat \beta, \hat d), (\alpha_0, \beta_0, d_0)\right) = O_p(1) \,.
%$$
%which immediately implies $\alpha_0, \beta_0$ can be estimated at $\sqrt{n}$ rate and $d_0$ can be estimated as $n$-rate (provided $X$ has mass near origin). 
\\\\
\noindent
Define $\theta = (\alpha, \beta, d)$. In the growing dimension regime, the rates of convergence of the estimates are affected by the underlying dimension.  We show later in this section (see Theorem \ref{thm:l2}) that for HEF with $0 \le k < \infty$ (i.e. excluding squared error loss), the corresponding Huber estimator satisfies: 
\begin{equation}
\label{eq:highdim_first_rate}
\frac{n}{p}\left(\log{\frac{n}{p}}\right)^{-1} \ \dist^2\left((\hat \alpha, \hat \beta, \hat d), (\alpha_0, \beta_0, d_0)\right) = O_p(1) \,.
\end{equation}
The above rate can be converted to a rate of convergence of the $\ell_2$ estimation error $\| \hat d - d_0\|$ of the change plane parameter via Assumption \ref{assm:wedge} stated below. In contrast, the rates for the least squares estimators are found to be slower and are non-trivially affected by the tail of the error distribution. Although the rate in equation \eqref{eq:highdim_first_rate} is shown to be minimax optimal for the change plane estimator  (Theorem \ref{thm:mlb_growing}), the rate of convergence of the one dimensional parameters $(\alpha_0, \beta_0)$ can be further boosted to $\sqrt{n}$ provided that the estimation error of the change plane estimator is smaller than $n^{-1/2}$ (i.e. $p \le \sqrt{n}/\log{(n/p)}$)using the following two step procedure:

\begin{enumerate}
%\item Split the data into two (almost) equal halves: $\cD = \cD_1 \cup \cD_2$. 
\item Get initial estimates of $(\alpha_0, \beta_0)$ and an estimate of $d_0$ as follows: 
$$(\hat \alpha^k_{\init}, \hat \beta^k_{\init}, \hat d^k) = \argmin_{\alpha, \beta, d} \sum_i \tilde H_k\left(Y_i - \alpha \mathds{1}_{X_i^{\top}d \le 0} - \beta \mathds{1}_{X_i^{\top}d > 0}\right) \,.$$
\item Update the estimates of $\alpha_0, \beta_0$ obtained in the previous step as follows: 
$$(\hat \alpha^k, \hat \beta^k) = \argmin_{\alpha, \beta} \sum_i \tilde H_k\left(Y_i - \alpha \mathds{1}_{X_i^{\top}\hat d^k \le 0} - \beta \mathds{1}_{X_i^{\top}\hat d^k > 0}\right) \,.$$
%\item Prescribe $(\hat \alpha^k, \hat \beta^k, \hat d^k)$ as the final estimates of the parameters. 
\end{enumerate}
%As mentioned previously we assume that $d \in S^{p-1}$ and $\alpha > \beta$ for identifiability purpose. 
The intuition for this rate acceleration is the following: if $\hat d^k$ converges to $d_0$ at a faster rate than $\sqrt{n}$, then we can re-estimate $(\alpha_0, \beta_0)$ at $\sqrt{n}$ - rate from the following surrogate model: 
$$
Y_i = \alpha_0 \mathds{1}_{X_i^{\top}\hat d^k \le 0} + \beta_0 \mathds{1}_{X_i^{\top}\hat d^k > 0} + \xi_i \,,
$$
where we simply replace $d_0$ by its estimate $\hat d^k$. If the estimation error of $\hat d^k$ is larger than $n^{-1/2}$, it is not possible to recover the parametric convergence rate for estimates of $(\alpha_0, \beta_0)$. %Note that in the fixed dimensional scenario, the change plane parameter $d_0$ can be estimated as $n$-rate, hence it is always possible to estimate $(\alpha_0, \beta_0)$ at $\sqrt{n}$ rate, whereas in growing dimension it depends on how fast the underlying dimension $p$ is growing with $n$. 
%Ignoring the log factor in the rate, it is easy to see that the rate acceleration upto $\sqrt{n}$ is possible when $n/p \gg \sqrt{n}$ i.e $p$ is growing slower than $\sqrt{n}$, otherwise the rate of estimation of $(\alpha_0, \beta_0)$ will be slower than $(n/p)$ up to a log factor. 
%It may be possible to achieve this acceleration in an one shot estimation method, but as this data-splitting does not influence the rate, we do not pursue this angle further. 
\\\\
\noindent
We now state our assumptions and the theorems. For technical simplicity, as before, we assume that $(\alpha_0, \beta_0) \in \Omega$ for a compact subset $\Omega \subseteq \bbR^2$ and for purposes of identifiability that $\alpha_0 > \beta_0$: 
\begin{assumption}
\label{assm:param_space}
Our parameter space $\Omega$ for $(\alpha, \beta)$ is a compact subset of $\bbR^2$ such that for any $(\alpha, \beta) \in \Omega$, $\alpha > \beta$. The hyperplane parameter $d_0 \in S^{p-1}$.  
\end{assumption}
Our next assumption (henceforth referred as \emph{wedge assumption}) relates the probability of $X$ lying in between two hyperplanes to the angle between those two hyperplanes. 
\begin{assumption}
\label{assm:wedge}
We assume there exists some $\delta > 0$ such that: 
\begin{align*}
\bbP\left(\s(X^{\top}d) \neq \s(X^{\top}d_0)\right) & \ge c \|d - d_0\|_2 \\
\bbP\left(X^{\top}d \wedge X^{\top}d_0 \ge 0\right) & \ge C_1 \\
\bbP\left(X^{\top}d \vee X^{\top}d_0 \le 0\right) & \ge C_2
\end{align*}
for all $\|d - d_0\| \le \delta$, where the constants $c, C_1, C_2, \delta$ do not depend on $n$. 
\end{assumption}
The first condition can be interpreted as saying that if we choose two hyperplanes $X^{\top}d = 0$ and $X^{\top}d_0 = 0$ the probability of $X$ falling in between these hyperplanes is bounded below, up to a constant, by the angle between the hyperplanes. This assumption can be thought as an analogue of the restricted eigenvalue assumption frequently used in the analysis of the high dimensional linear model (especially LASSO, see e.g. \cite{bickel2009simultaneous}) to obtain the estimation error from the prediction error and was also used in earlier work by the authors \cite{mukherjee2019non}, where it was shown that the condition is satisfied by several classes of distributions (e.g. under elliptical symmetry, log-concavity of densities).  The second and third inequalities are weak assumptions, which ensure that the support of $X$ is not restricted to the one side of the hyperplane. A sufficient condition for the above assumption is that the density of $X$ is uniformly bounded away from 0 on a ball of fixed radius around the origin, where the uniformity is over the dimension. 
\newline
\newline
We next state our theorems for this regime: 
%\begin{assumption}
%\label{assm:wedge_global}
%We assume that: 
%\begin{align*}
%\bbP\left(\s(X^{\top}\theta) \neq \s(X^{\top}\theta_0)\right) \ge c \|\theta - \theta_0\|_2 
%\end{align*}
%holds for all $\theta \in S^{p-1}$. 
%\end{assumption}
%\begin{remark}
%{\bf COMMENT ON THE ESTIMATION PROCEDURE} 
%\end{remark}
%\noindent

\begin{theorem}[Rate of convergence]
\label{thm:l2}
Suppose we estimate $\theta_0 = (\alpha_0, \beta_0, d_0)$ using the two-shot approach described above, i.e.  by minimizing the scaled HEF and then re-estimating $(\alpha_0, \beta_0)$. Then, under Assumptions \ref{assm:param_space}-\ref{assm:wedge}, we have for $0 \le k < \infty$: 
\begin{align*}
\left(\sqrt{n} \wedge \frac{n}{p}\left(\log{\frac{n}{p}}\right)^{-1}\right)\left(\hat \alpha^k - \alpha_0\right) = O_p(1), \  \ \left(\sqrt{n} \wedge \frac{n}{p}\left(\log{\frac{n}{p}}\right)^{-1}\right)\left(\hat \beta^k - \beta_0\right) & = O_p(1), \\
& \hspace{-25em} \frac{n}{p}\left(\log{\frac{n}{p}}\right)^{-1}\bbP\left(\s(X^{\top}\hat d^k) \neq \s(X^{\top}d_0)\right)  = O_p(1) \,,
\end{align*}
which along with Assumption \ref{assm:wedge} yields: 
$$
 \frac{n}{p}\left(\log{\frac{n}{p}}\right)^{-1}\left\|\hat d^k - d_0\right\|_2 = O_p(1) \,.
$$
For $k = \infty$, i.e. under squared error loss we have: 
\begin{align*}
& \left(\sqrt{n} \wedge \frac{n}{p\left\|\xi\right\|_{n, L_1}}\left(\log{\frac{n}{p\left\|\xi\right\|_{n, L_1}}}\right)^{-1}\right)\left(\hat \alpha^{\ell_2} - \alpha_0\right) = O_p(1)  \\
& \left(\sqrt{n} \wedge \frac{n}{p\left\|\xi\right\|_{n, L_1}}\left(\log{\frac{n}{p\left\|\xi\right\|_{n, L_1}}}\right)^{-1}\right)\left(\hat \beta^{\ell_2} - \beta_0\right) = O_p(1) \\
&\frac{n}{p\left\|\xi\right\|_{n, L_1}}\left(\log{\frac{n}{p\left\|\xi\right\|_{n, L_1}}}\right)^{-1}\left\|\hat d^{\ell_2} - d_0\right\|_2  = O_p(1) \,.
\end{align*}
where $\left\|\xi\right\|_{n, L_1} = \bbE\left[\max_{1 \le i \le n} |\xi_i|\right]$. 
%Therefore the rate is compromised when we use squared error loss instead of any other Huber estimating equation $\tilde H_k$. 
\end{theorem}

\noindent
Like the results of the previous subsection, Theorem \ref{thm:l2} shows that in a growing dimension setting the rate of convergence of the Huber estimator for any $0 \le k < \infty$ is faster than using the standard squared error loss: the rate of the least squares estimator of $d_0$ suffers from an additional factor $\left\|\xi\right\|_{n, L_1}$, which depends on the tail of the distribution of $\xi$. We note that this is not that an isolated phenomenon, e.g. in non-parametric regression, the rate of convergence of the least square estimators is similarly affected by the tail of the error, e.g. see \cite{han2019convergence}. Note that in the fixed $p$ regime this factor can be ignored via a different maximal inequality which, when used in growing dimensional regime, yields the rate $n/p^2$. More specifically, consider Lemma 2.14.1 of \cite{vdvw96}, which we state here for the ease of our readers: 
$$
\bbE\left[\sqrt{n}\sup_{f \in \cF} \left|(\bbP_n - P)f\right|\right] \lesssim \bbE\left[\cJ(\theta_n, \cF)\|F\|_{2, n}\right] \lesssim \cJ(1, \cF)\sqrt{\bbE[F^2]} \,,
$$
where $F$ is the envelope of $\cF$ and $\cJ$ quantifies the complexity of $\cF$ as follows: 
$$
\cJ(\delta, \cF) = \sup_Q \int_0^{\delta} \sqrt{1 + \log{N(\eps\|F\|_{Q, 2}, \cF, L_2(Q))}} \ d\eps \,.
$$
and $\theta_n = \sup_{f \in \cF} \left\|f/F\right\|_{2, n}$ (with the convention $0/0 = 0$). 
The rate of convergence of the least squares estimator in Theorem \ref{thm:l2} is obtained via a modified version of the first inequality (details can be found in the proof), whereas one may also use the weaker second inequality which, in this case, yields the rate $(n/p^2)$. Combining this with the rate obtained in Theorem \ref{thm:l2} leads to the following modified rate of convergence:  
$$
\left\|\hat d^{\ell_2} - d_0\right\|_2 = O_p\left(\frac{p\left\|\xi\right\|_{n, L_1}}{n}\log{\frac{n}{p\left\|\xi\right\|_{n, L_1}}} \wedge \frac{p^2}{n}\right) \,.
$$ 
When $p$ is fixed, the second weaker inequality yields a better rate of convergence as the factor $p^2$ is a constant. In the growing dimension regime, the VC dimension of the underlying function class is growing with the sample size, hence the interplay between the ambient dimension $p$ and the tail of the error distribution starts affecting the rate of convergence of the least squares estimator. However, the tail factor $\|\xi\|_{n, L_1}$ does not appear in the rate of the other Huber estimators, as the criterion function becomes bounded irrespective of the thickness of tail of the distribution. 
\\\\
\noindent
To summarize,  we have established in Section \ref{sec:one_dim_findings} that the robust Huber estimators (for any $0 \le k < \infty$) yield a more concentrated limiting distribution than the squared error loss, whereas in the growing dimension regime, the effect is more prominent: robust Huber estimators yield a faster rate of convergence, which is also minimax optimal as shown in our next theorem. This underscores the necessity of using robust estimators in high dimensional change plane problems, especially in presence of heavy tailed errors.

\begin{theorem}[Minimax lower bound]
\label{thm:mlb_growing}
Suppose $\cP = \left\{P_d: d \in S^{d-1} \right\}$ is the collection of all change plane models such that the distribution $P_{d}$ of $(X, Y)$ or equivalently the distribution of $(X, \xi)$ satisfies the following: 
$$
Y = \mathds{1}_{X^{\top}d > 0} + \xi 
$$
where $X$ is independent of $\xi$ and $\xi$ is symmetric around the origin. Then we have: 
$$
\inf_{\hat d} \sup_{P_d} \bbE_{\theta}\left(\dist^2(\hat d, d)\right) \ge K\frac{p}{n}\left(1 + \log{\frac{p}{n}}\right)  \,,
$$
where the semi-metric $\dist$ is defined as: 
$$
\dist(d_1, d_2) = \sqrt{\bbP\left(\s(X^{\top}d_1) \neq \s(X^{\top}d_2)\right)} \,.
$$
Hence, the change plane estimator obtained in Theorem \ref{thm:l2} via the Huber estimating equation $\tilde H_k$ for $0 \le k < \infty$ (i.e. excluding squared error loss) is minimax optimal. 

\end{theorem}

\begin{remark}
\label{rem:mlb_growing_dim}
Notice that we restrict our minimax calculation only to the change plane parameter $d_0$ assuming we know $(\alpha_0, \beta_0)$ (in fact, without loss of generality we assume $\alpha_0 = 0, \beta_0 = 1$), as the minimaxity of the rate of convergence of $(\alpha_0, \beta_0)$ is immediate and not interesting. Theorem \ref{thm:mlb_growing} indicates that any Huber estimator for $0 \le k < \infty$ is minimax optimal. The proof of this theorem relies on a clever construction of the local alternatives and an application of Fano's inequality (e.g. see \cite{yu1997assouad}).  
\end{remark}

\subsubsection{When $p \gg n$}
\label{sec:fast_growing}
In this section, we present our analysis of the change plane estimator in the regime $p \gg n$, i.e. the HDLSS (high dimension low sample size) setting. As is true for any high dimensional model, consistent estimate of $d_0$ is not information theoretically possible without further restrictions on the parameter space. A typical condition frequently imposed on the parameter space is that of sparsity: there exists some (unknown) $s$ such that only $s$ many elements of $d_0$ are non-zero, where $s$ may also increase with the sample size. We summarize this in the following assumption: 
\begin{assumption}
\label{assm:sparsity}
The true change plane direction $d_0$ is sparse, i.e. there exists $s$ such that $\|d_0\|_0 \le s$ where $s$ may slowly grow with $n$, satisfying $(s\log{p})/n \to 0$. 
\end{assumption}
To estimate $d_0$ under this sparsity constraint, we follow the \emph{structural risk minimization} method, an idea originally from \cite{vapnik1974theory} and later implemented in a series of work (e.g. \cite{massart2007concentration}, \cite{barron1999risk}, \cite{bartlett2002model} and references therein). The key idea is to use a penalty function to balance the bias-variance trade-off. To understand this, consider our stump model: 
$$
Y_i = \alpha_0 \mathds{1}_{X_i^{\top}d_0 \le  0} + \beta_0 \mathds{1}_{X_i^{\top}d_0 > 0} + \xi_i \,.
$$
%$$
%Y_i = \mathds{1}_{X_i^{\top}d_0 > 0} + \xi_i \,.
%$$
Define $\cF_m$ to be set of all hyperplanes with sparsity at-most $m$, i.e.: 
$$
\cF_m = \left\{f_{\theta}(X) = \alpha\mathds{1}_{X^{\top}d \le 0} + \beta\mathds{1}_{X^{\top}d > 0}: (\alpha, \beta) \in \Omega, \left\|d\right\|_0 \le m\right\} \,,
$$ 
for $1 \le m \le p$. Now for each $m$, we define the empirical minimizer $\hat \theta^k_m := \hat \theta^k_{m, n}$ as: 
\begin{align*}
\hat \theta^k_m & = \argmin_{\theta: f_{\theta} \in \cF_m} \frac{1}{n}\sum_{i=1}^n \tilde H_k\left(Y_i - f_\theta(X_i)\right) \\
& =  \argmin_{\theta: f_{\theta} \in \cF_m} \frac{1}{n}\sum_{i=1}^n \left[\tilde H_k\left(Y_i - f_\theta(X_i)\right) - \tilde H_k\left(Y_i - f_{\theta_0}(X_i)\right)\right] \\
& = \argmin_{\theta: f_{\theta} \in \cF_m} \frac{1}{n}\sum_{i=1}^n \left[\tilde H_k\left(Y_i - f_\theta(X_i)\right) - \tilde H_k\left(\xi_i\right)\right] 
\end{align*}
for $0 \le k \le \infty$. The corresponding population minimizer is defined as: 
$$
\theta^k_m = \argmin_{d: f_d \in \cF_m} \bbE\left[\tilde H_k\left(Y - f_\theta(X)\right) - \tilde H_k\left(\xi\right)\right]  \,.
$$
Note that, the larger the $m$, the more complex is the function class $\cF_m$ (as $\cF_{m_1} \subseteq \cF_{m_2}$ for any $m_1 \le m_2$), and consequently,  the variance starts dominating the bias for large values of $m$. In other words, $\hat \theta^k_m$ has smaller training error, but larger generalization error for large values of $m$. Therefore, to choose an optimal model $m$, we add a penalty $\pen(m)$ (which quantifies the complexity of the model $\cF_m$ and is increasing in $m$) to the training error of $\hat \theta^k_m$ and choose the one which minimizes the penalized training error: 
$$
\hat m^k = \argmin_{1 \le m \le p}  \frac{1}{n}\sum_{i=1}^n \left[\tilde H_k\left(Y_i - f_{\hat \theta^k_m}(X_i)\right) - \tilde H_k(\xi_i)\right] +  \pen(m) \,.
$$
and set the final estimator as $\hat \theta^k_{\hat m}$. The penalty function should be chosen carefully depending on the complexity of the underlying function class to balance the bias-variance tradeoff. We quantify the complexity of $\cF_m$ using its VC dimension. It follows from Lemma 1 of \cite{abramovich2017high} that the VC dimension of $\cF_m$ is, 
\begin{equation}
\label{eq:vcdim_highd}
V_m = VC\left(\cF_m\right) \asymp m\log{\frac{ep}{m}} \,.
\end{equation}
Based on the above notion of complexity, we use the following penalty function: 
\begin{equation}
\label{eq:penalty_huber}
\pen(m) = \kappa \left(\frac{V_m \log{(n/V_m)}}{n}\right) \,.
\end{equation}
for some constant $\kappa$ independent of $n$, while using HEF $\tilde H_k$ for $0 \le k < \infty$. For $k = \infty$, i.e. while using the least squares estimator, we use a slightly different penalty (see Theorem \ref{thm:l2_high} for more details). Therefore, our $\pen$ function is the VC dimension of the model under consideration (up to a constant and a log factor). Finally, we can accelerate the rate of convergence of $(\alpha_0, \beta_0)$ by following the same procedure as prescribed in Subsection \ref{sec:slowly_growing}): i.e., first estimate $d_0$ by minimizing the penalized criterion function, then re-estimate $(\alpha_0, \beta_0)$ using $\hat d^k_{\hat m^k}$ as a proxy for $d_0$. Henceforth, we denote by $(\hat \alpha^k, \hat \beta^k, \hat d_{\hat m^k})$ as these final estimators obtained via the two-shot procedure.  
\\\\
\noindent
Our next theorem presents the rate of convergence of the above estimates %$\hat d_{\hat m}$ obtained via method elaborated above using HEF $\tilde H_k$, 
for $0 \le k < \infty$: 

\begin{theorem}
\label{thm:huber_high}
%Let $(\hat \alpha^k, \hat \beta^k, \hat d_{\hat m^k})$ be the obtained via two-shot procedure described above, i.e. first by minimizing the penalized HEF $\tilde H_k$ (for any $0 \le k < \infty$) and then re-estimating $(\alpha_0, \beta_0)$.
Under Assumptions \ref{assm:param_space}, \ref{assm:wedge} and \ref{assm:sparsity} and using the penalty introduced in \eqref{eq:penalty_huber} we have: 
\begin{align*}
& \left(\sqrt{n} \wedge \frac{n}{V_{s} \log{\frac{n}{V_{s}}}} \right)\left(\hat \alpha - \alpha_0\right) = O_p(1), \\
& \left(\sqrt{n} \wedge \frac{n}{V_{s} \log{\frac{n}{V_{s}}}} \right) \left(\hat \beta - \beta_0\right) = O_p(1), \\
& \frac{n}{V_{s} \log{\frac{n}{V_{s}}}}\left\|\hat d - d_0\right\|_2  = O_p(1) \,.
\end{align*}
\end{theorem}

\begin{remark}
From equation \eqref{eq:vcdim_highd},  it is readily seen that the rate of convergence of the change plane estimator $\hat d$ is: 
$$
\left\| \hat d - d_0 \right\|_2 = O_p\left(\frac{s\log{(ep/s)}}{n}\log{\left(\frac{n}{s\log{(ep/s)}}\right)}\right)
$$ 
i.e. upto a log factor, the rate if $s\log{p}/n$, which can be thought as the high dimensional analogue of $1/n$ (the rate obtained for the change point estimator in finite dimension) in presence of sparsity. 
\end{remark}

We now present our results regarding the minimax lower bound for this change plane problem in this HDLSS scenario under the sparsity constraint. As before, we restrict our attention to the parameter of interest $d_0$ and assume we know $\alpha_0, \beta_0$, in particular setting $\alpha_0 = 0, \beta_0 = 1$. 
\begin{theorem}
\label{thm:mlb_high}
%Consider similar setup as in Theorem \ref{thm:mlb_growing}.
Assume $\cP = \left\{P_d: d \in S^{p-1}_s\right\}$ is the collection of all change plane models  with $S^{p-1}_{s}$ being the set of all unit vectors in dimension $p$ with sparsity at-most $s$, where the distribution $P_{d}$ of $(X, Y)$ or equivalently the distribution of $(X, \xi)$ satisfies the following: 
$$
Y =  \mathds{1}_{X^{\top}d > 0} + \xi \,,
$$
where $X$ is independent of $\xi$ and $\xi$ is symmetric around the origin. 
Then we have: 
$$
\inf_{\hat \theta} \sup_{P_\theta} \bbE_{\theta}\left(\|\hat d - d_0\|^2\right) \ge K\left(\frac{s\log{(ep/s)}}{n}\right)^2  \,.
$$
%where the semi-metric $\dist$ is defined as (similar to that of Theorem \ref{thm:mlb_growing}: 
%$$
%\dist(\theta_1, \theta_2) = \sqrt{(\alpha_1 - \alpha_2)^2 + (\beta_1 - \beta_2)^2 + \bbP\left(\s(X^{\top}d_1) \neq \s(X^{\top}d_2)\right)} \,.
%$$
\end{theorem}

\noindent
\paragraph{Squared error loss: }The rate of convergence of the least squares estimator for this regime is also compromised by the tail of the error distribution, which is in agreement with our findings in Subsection \ref{sec:slowly_growing}. To establish the theoretical properties of the LSE, we slightly strengthen our sparsity assumption below:  
\begin{assumption}
\label{assm:sparsity_l2}
The true change plane direction $d_0$ is sparse, i.e. there exists $s$ such that $\|d_0\|_0 \le s$ where $s$ may slowly grow with $n$, satisfying
$$
\frac{s(\log{p})^{(1+\delta)}\|\xi\|_{n, 2}}{n} \to 0 
$$ 
where $\|\xi\|_{n, 2} =  \sqrt{\bbE[\max_{1 \le i \le n} \xi_i^2]}$. 
\end{assumption}
Two comments on this modified sparsity assumption are in order: first note that, we need a slightly higher power of $\log{p}$ in comparison to its counterpart in Assumption \ref{assm:sparsity}. This is likely a technical artifact and possibly avoidable with more tedious analysis. Next, we also have an additional term $\|\xi\|_{n, 2}$ which captures the effect of the tail of the error distribution in the rate of the LSE, similar to what we see in Theorem \ref{thm:l2}.  This modified assumption necessitates changing our penalty to: 
\begin{equation}
\label{eq:penalty_l2}
\pen(m) = \frac{V_m(\log{p})^{\delta}\|\xi\|_{n, 2}}{n}\log{\frac{n}{V_m}}
\end{equation} 
where, as before, $V_m$ is the VC dimension of $\cF_m$. The following theorem establishes the rate of convergence of the LSE.  

\begin{theorem}
\label{thm:l2_high}
Suppose we estimate $\theta_0 = (\alpha_0, \beta_0, d_0)$ using the two-shot procedure under squared error loss. Then under Assumptions \ref{assm:param_space}, \ref{assm:wedge} and \ref{assm:sparsity_l2} we obtain: 
\begin{align*}
& \left(\sqrt{n} \wedge  \frac{n}{s (\log{p})^{(1+\delta)}\|\xi\|_{n, 2}}\left(\log{\frac{n}{s\log{p}}}\right)^{-1}\right)\left(\hat \alpha - \alpha_0\right) = O_p(1), \\
& \left(\sqrt{n} \wedge  \frac{n}{s (\log{p})^{(1+\delta)}\|\xi\|_{n, 2}}\left(\log{\frac{n}{s\log{p}}}\right)^{-1}\right)\left(\hat \beta - \beta_0\right) = O_p(1), \\
&  \frac{n}{s (\log{p})^{(1+\delta)}\|\xi\|_{n, 2}}\left(\log{\frac{n}{s\log{p}}}\right)^{-1}\left\|\hat d - d_0\right\|_2  = O_p(1) \,.
\end{align*}
\end{theorem}

\begin{remark}
\label{rem:mlb_high}
A remark similar to Remark \ref{rem:mlb_growing_dim} is in order: Theorem \ref{thm:mlb_high} conveys a similar message as Theorem \ref{thm:mlb_growing}, i.e. any Huber-estimator for $0 \le k < \infty$ is minimax optimal up to a log factor, whereas the least squares  estimator is not (as seen above), especially when the distribution of $\xi$ has a heavy tail. Therefore, as in the previous subsection, robust Huber-estimators are preferable to the least squares estimator in this high dimensional regime.  
\end{remark}

\section{An empirical study of the quantiles of the limiting distributions}
\label{sec:simulation}
In this section we present tables of quantiles of the limit distributions of change point estimator under both the $\ell_1$ and $\ell_2$ criteria. In Section \ref{sec:one_dim_findings}, we established theoretically (Theorem \ref{cor:l1_loss} and \ref{cor:l2_loss}) that in the presence of heavy tailed errors, the limiting distribution of the change point estimator under $\ell_1$ criterion has a thinner tail (i.e. more concentrated asymptotic confidence interval) than the change point estimator under $\ell_2$ criterion. We provide some illustrations of this phenomenon in our simulations below. 
\\\\
\noindent
We generate data from the following stump model: 
$$
Y_i = \mu \mathds{1}_{X_i \ge d_0} + \xi_i \,.
$$ 
where we have assumed $d_0 = 0$, $X_i \sim \cN(0, 1)$. Recall that the limiting distribution of $\hat d^{\ell_1}$ is (see Theorem \ref{cor:l1_loss}): 
$$
n(\hat d^{\ell_1} - d_0)  \overset{\mathscr{L}}{\implies} \sargmin_{t \in \bbR}  \CPP\left(|\xi + \mu| - |\xi|, f_X(d_0)\right)
$$
and the limiting distribution of $\hat d^{\ell_2}$ is (see Theorem \ref{cor:l2_loss}): 
$$
n(\hat d^{\ell_2} - d_0) \overset{\mathscr{L}}{\implies} \sargmin_{t \in \bbR}  \CPP\left(\xi + \frac{\mu}{2}, f_X(d_0)\right)
$$
For $\xi$, we consider seven different distributions: standardized $T_3, T_4, T_5, T_6, T_{10}, T_{15}$ (i.e. $\var = 1$) and $\cN(0, 1)$, while for the signal $\mu$ we consider four different values $\mu = 0.1, 0.5, 1, 2$. We present here 8 different tables: two tables ($\ell_1$ and $\ell_2$ quantiles) for each value of $\mu$.  
%where one for the limiting distribution under $\ell_1$ loss and the other for the limiting distribution under $\ell_2$ loss function. 
Each table consists of five different one sided quantiles  ($90\%, 95\%, 97.5\%, 99\%, 99.5\%$) for each of the five different distributions of $\xi$ (calculated using $10^6$ monte-carlo iterations). Recall that as we compute the \sargmin change-point estimator, the limit distributions are symmetric and it suffices to report the upper quantiles. The percentages presented inside the brackets following the quantiles in the even-numbered tables show the relative change in the $\ell_2$ based quantile as compared to the $\ell_1$ based counterpart. 

\begin{table}[H]
\centering
\caption{Quantiles of asymptotic distribution under $\ell_1$ criterion using $\mu = 0.1$}
\label{tab:l1_one_tenth}
%\resizebox{\textwidth}{!}{%
\begin{tabular}{@{}llllll@{}}
\toprule
Distributions       & 90\%         & 95\%         & 97.50\%      & 99\%         & 99.50\%      \\ \midrule
$T_3$    & 717.0 & 1152.3 & 1600.1 & 2100.3 & 2324.8 \\
$T_4$    & 943.7  & 1470.5 & 1924.0 & 2308.9 & 2432.3  \\
$T_5$    & 1062.0  & 1611.0 & 2045.6 & 2366.8 & 2460.7  \\
$T_6$    & 1133.6 & 1690.1  & 2110.6 & 2389.5 & 2475.5 \\
$T_{10}$    & 1247.8 & 1808.5  & 2196.3 & 2419.7 & 2489.5 \\
$T_{15}$    & 1294.4 & 1859.2  & 2229.8 & 2433.0 & 2499.5 \\
Normal & 1381.6 & 1944.1 & 2278.8 & 2449.6 & 2509.1  \\ \bottomrule
\end{tabular}%
%}
\end{table}

%\begin{table}[H]
%\centering
%\caption{Quantiles of asymptotic distribution under $\ell_2$ loss using $\mu = 0.1$}
%\label{ttab:l2_one_tenth}
%\resizebox{\textwidth}{!}{%
%\begin{tabular}{@{}llllll@{}}
%\toprule
%  Distributions     & 90\%        & 95\%        & 97.50\%     & 99\%        & 99.50\%     \\ \midrule
%$T_3$   & 669.9369204 & 1022.472695 & 1313.3702   & 1546.591579 & 1644.453057 \\
%$T_4$   & 802.4372688 & 1204.799208 & 1527.727669 & 1784.931242 & 1877.7586   \\
%$T_5$   & 857.6812917 & 1302.11395  & 1668.747998 & 1929.153078 & 2025.087224 \\
%$T_6$  & 897.5291741 & 1362.775772 & 1740.866282 & 2016.230318 & 2116.212543 \\
%Normal & 1097.997359 & 1667.913151 & 2116.975246 & 2457.623591 & 2562.040342 \\ \bottomrule
%\end{tabular}%
%}
%\end{table}

\begin{table}[H]
\centering
\caption{Quantiles of asymptotic distribution under $\ell_2$ criterion using $\mu = 0.1$}
\label{ttab:l2_one_tenth}
%\resizebox{\textwidth}{!}{%
\begin{tabular}{@{}llllll@{}}
\toprule
  Distributions     & 90\%        & 95\%        & 97.50\%     & 99\%        & 99.50\%     \\ \midrule
$T_3$   & 1045.7(+45.8\%) & 1584.7(+37.5\%) & 2029.7(+26.8\%)   & 2358.3(+12.3\%) & 2457.6(+5.7\%) \\
$T_4$   & 1050.1(+11.3\%) & 1598.2(+8.7\%) & 2034.3(+5.7\%) & 2355.6(+2\%) & 2456.3(+1\%)   \\
$T_5$   & 1055.6(-0.6\%) & 1600.9(-0.6\%)  & 2040.0(-0.3\%) & 2364.6(-0.1\%) & 2461.0(+0.01\%) \\
$T_6$  & 1056.2(-5.1\%) & 1601.2(-5.3\%) & 2043.0(-3.2\%) & 2366.0(-1\%) & 2461.1(-0.6\%) \\
$T_{10}$  & 1052.5(-15.6\%) & 1593.9(-11.9\%) & 2038.0(-7.2\%) & 2363.0(-2.3\%) & 2460.05(-1.2\%) \\
$T_{15}$  & 1054.8(-18.5\%) & 1601.0(-13.9\%) & 2046.2(-8.2\%) & 2366.1(-2.75\%) & 2461.0(-1.5\%) \\
Normal & 1051.4(-24\%) & 1600.9(-17.6\%) & 2044.5(-10.3\%) & 2363.9(-3.49\%) & 2460.0(-2\%) \\ \bottomrule
\end{tabular}%
%}
\end{table}

%\begin{table}[H]
%\centering
%\caption{Quantiles of asymptotic distribution under $\ell_1$ loss using $\mu = 0.5$}
%\label{tab:l1_one_half}
%\resizebox{\textwidth}{!}{%
%\begin{tabular}{@{}llllll@{}}
%\toprule
%Distributions       & 90\%        & 95\%        & 97.50\%     & 99\%        & 99.50\%     \\ \midrule
%$T_3$   & 17.49855479 & 29.93751147 & 43.70658659 & 64.85226543 & 82.50854664 \\
%$T_4$   & 27.92057467 & 46.71019611 & 67.93084899 & 98.01357333 & 121.0362425 \\
%$T_5$   & 34.7836682  & 58.0582129  & 84.17905743 & 123.0240624 & 154.5959854 \\
%$T_6$   & 40.92875357 & 67.732538   & 97.60320301 & 140.150837  & 175.9696741 \\
%Normal & 66.45295483 & 110.6197381 & 158.3821674 & 228.5698156 & 285.8810559 \\ \bottomrule
%\end{tabular}%
%}
%\end{table}

\begin{table}[H]
\centering
\caption{Quantiles of asymptotic distribution under $\ell_1$ criterion using $\mu = 0.5$}
\label{tab:l1_one_half}
%\resizebox{\textwidth}{!}{%
\begin{tabular}{@{}llllll@{}}
\toprule
Distributions       & 90\%        & 95\%        & 97.50\%     & 99\%        & 99.50\%     \\ \midrule
$T_3$   & 28.5 & 46.8 & 67.7 & 98.1 & 122.5 \\
$T_4$   & 38.2 & 62.7 & 89.9 & 129.7 & 162.2 \\
$T_5$   & 44.5  & 73.1  & 104.8 & 150.7 & 188.1 \\
$T_6$   & 48.3 & 78.9   & 113.5 & 163.3  & 203.6 \\
$T_{10}$   & 55.7 & 90.7   & 130.6 & 187.4  & 233.4 \\
$T_{15}$   & 59.0 & 97.0   & 139.7 & 200.5  & 250.0 \\
Normal & 66.1 & 108.3 & 155.7 & 224.8 & 279.9 \\ \bottomrule
\end{tabular}%
%}
\end{table}

%\begin{table}[H]
%\centering
%\caption{Quantiles of asymptotic distribution under $\ell_1$ loss using $\mu = 0.5$}
%\label{tab:l1_one_half}
%\resizebox{\textwidth}{!}{%
%\begin{tabular}{@{}llllll@{}}
%\toprule
% Distributions      & 90\%        & 95\%        & 97.50\%     & 99\%        & 99.50\%     \\ \midrule
%$T_3$   & 29.47839556 & 49.92846345 & 73.50554006 & 108.09319   & 135.1728329 \\
%$T_4$   & 34.93083131 & 58.95249007 & 84.81185902 & 121.406829  & 150.6424914 \\
%$T_5$   & 37.95710086 & 63.62867054 & 92.65865144 & 133.9274492 & 168.728778  \\
%$T_6$   & 40.43559765 & 67.37334858 & 96.2284684  & 137.3622588 & 169.2800555 \\
%Normal & 49.34790968 & 80.89201888 & 116.4527308 & 167.5209031 & 208.133535  \\ \bottomrule
%\end{tabular}%
%}
%\end{table}

\begin{table}[H]
\centering
\caption{Quantiles of asymptotic distribution under $\ell_2$ criterion using $\mu = 0.5$}
\label{tab:l1_one_half}
%\resizebox{\textwidth}{!}{%
\begin{tabular}{@{}llllll@{}}
\toprule
 Distributions      & 90\%        & 95\%        & 97.50\%     & 99\%        & 99.50\%     \\ \midrule
$T_3$   & 45.8(+60.7\%) & 77.3(+65.2\%) & 113.5(+67.6\%) & 166.9(+70.1\%)   & 207.7(69.5\%) \\
$T_4$   & 47.3(+23.8\%) & 78.3(+24.9\%) & 113.3(+26\%) & 163.5(+26.1\%)  & 204.5 (+26.1\%) \\
$T_5$   & 48.1(+8.1\%) & 78.9(+7.9\%) & 113.0(+7.8\%) & 163.3(+8.4\%) & 202.7(+7.8\%)  \\
$T_6$   & 48.3(+0\%) & 79.0(+0.1\%) & 113.5(+0\%)  & 162.8(-0.3\%) & 203.1(-0.25\%) \\
$T_{10}$   & 48.6(-12.7\%) & 79.1(-12.8\%)   & 113.3(-13.2\%) & 162.6(-13.2\%)  & 202.4(-13.3\%) \\
$T_{15}$   & 48.8(-17.3\%) & 79.3(-18.2\%)   & 113.7(-18.6\%) & 162.8(-18.8\%)  & 203.2(-18.72\%) \\
Normal & 48.8(-26.2\%) & 79.4(-26.7\%) & 113.8(-27\%) & 163.0(-27.5\%) & 202.7(-27.6\%)  \\ \bottomrule
\end{tabular}%
%}
\end{table}

%\begin{table}[H]
%\centering
%\caption{Quantiles of asymptotic distribution under $\ell_1$ loss using $\mu = 1$}
%\label{tab:l1_one}
%\resizebox{\textwidth}{!}{%
%\begin{tabular}{@{}llllll@{}}
%\toprule
% Distributions      & 90\%        & 95\%        & 97.50\%     & 99\%        & 99.50\%     \\ \midrule
%$T_3$   & 4.473446571 & 7.827882836 & 11.68386868 & 17.39970614 & 22.2382867  \\
%$T_4$   & 7.011291399 & 11.93218875 & 17.63746442 & 25.33109133 & 32.2561399  \\
%$T_5$   & 8.360983101 & 14.45156927 & 20.90820518 & 30.28423323 & 37.88892545 \\
%$T_6$   & 9.526928933 & 16.31017097 & 23.62834981 & 34.11996575 & 42.63179696 \\
%Normal & 15.15860988 & 25.54339525 & 37.41115434 & 54.69608119 & 68.92845491 \\ \bottomrule
%\end{tabular}%
%}
%\end{table}

\begin{table}[H]
\centering
\caption{Quantiles of asymptotic distribution under $\ell_1$ criterion using $\mu = 1$}
\label{tab:l1_one}
%\resizebox{\textwidth}{!}{%
\begin{tabular}{@{}llllll@{}}
\toprule
 Distributions      & 90\%        & 95\%        & 97.50\%     & 99\%        & 99.50\%     \\ \midrule
$T_3$   & 8.4 & 13.5 & 19.2 & 27.6 & 34.4  \\
$T_4$   & 10.5 & 16.9 & 24.0 & 34.3 & 42.7  \\
$T_5$   & 11.7 & 18.7 & 26.8 & 38.2 & 47.5 \\
$T_6$   & 12.5 & 20.2 & 28.9 & 41.2 & 51.3 \\
$T_{10}$   & 14 & 22.7   & 32.4 & 46.4  & 57.8 \\
$T_{15}$   & 14.7 & 23.9   & 34.1 & 48.9  & 60.9 \\
Normal & 16.1 & 26.2 & 37.2 & 53.7 & 66.7 \\ \bottomrule
\end{tabular}%
%}
\end{table}

%\begin{table}[H]
%\caption{Quantiles of asymptotic distribution under $\ell_2$ loss using $\mu = 1$}
%\label{tab:l2_one}
%\label{tab:my-table}
%\resizebox{\textwidth}{!}{%
%\begin{tabular}{@{}llllll@{}}
%\toprule
%Distributions      & 90\%        & 95\%        & 97.50\%     & 99\%        & 99.50\%               \\ \midrule
%$T_3$ & 6.751451354  & 12.469750627 & 19.098894855 & 29.014808708 & 36.915381378 \\
%$T_4$ & 8.4567381599 & 14.894298642 & 22.160527075 & 33.05404447  & 41.838441641 \\
%$T_5$ & 9.4225064777 & 16.239343504 & 23.828732086 & 34.428362797 & 43.193434901 \\
%$T_6$ & 10.091225819 & 17.22582063  & 24.795123564 & 36.158021188 & 45.619367878 \\
%Normal & 12.611931133 & 21.363317748 & 31.180733939 & 44.674729865 & 56.090374968 \\ \bottomrule
%\end{tabular}%
%}
%\end{table}

\begin{table}[H]
\centering
\caption{Quantiles of asymptotic distribution under $\ell_2$ criterion using $\mu = 1$}
\label{tab:l2_one}
\label{tab:my-table}
%\resizebox{\textwidth}{!}{%
\begin{tabular}{@{}llllll@{}}
\toprule
Distributions      & 90\%        & 95\%        & 97.50\%     & 99\%        & 99.50\%               \\ \midrule
$T_3$ & 11.8(+40.5\%)  & 20.3(+50.4\%) & 30.7(+59.9\%) & 46.1(+67.0\%) & 59.3(+72.4\%) \\
$T_4$ & 12.7(+21\%) & 21.0(+24.3\%) & 30.5(+27.1\%) & 44.6(+30.1\%)  & 56.0(+31.1\%) \\
$T_5$ & 13.0(+11.1\%) & 21.3(+14\%) & 30.5(+13.8\%) & 44.0(+15.2\%) & 55.0(+15.8\%) \\
$T_6$ & 13.1(+4.8\%) & 21.4(+5.9\%) & 30.6(+5.9\%) & 44.0(+6.8\%) & 55.0(+7.2\%) \\
$T_{10}$   & 13.4(-4.28\%) & 21.5(-5.29\%)   & 30.7(-5.25\%) & 43.8(-5.6\%)  & 54.1(-6.4\%) \\
$T_{15}$   & 13.4(-8.8\%) & 21.6(-9.6\%)   & 30.7(-10\%) & 43.9(-10.2\%)  & 54.0(-11.3\%) \\
Normal & 13.5(-16.1\%) & 21.7(-17.2\%) & 30.6(-17.1\%) & 43.5(-19\%) & 54.0(-19\%) \\ \bottomrule
\end{tabular}%
%}
\end{table}

%\begin{table}[H]
%\centering
%\caption{Quantiles of asymptotic distribution under $\ell_1$ loss using $\mu = 2$}
%\label{tab:l1_two}
%\resizebox{\textwidth}{!}{%
%\begin{tabular}{@{}llllll@{}}
%\toprule
%Distributions       & 90\%        & 95\%        & 97.50\%     & 99\%        & 99.50\%     \\ \midrule
%$T_3$   & 0.451262522 & 2.162571288 & 3.876515588 & 6.190487877 & 7.984534787 \\
%$T_4$   & 1.105959844 & 3.165713134 & 5.301901576 & 8.159910453 & 10.45177312 \\
%$T_5$   & 1.499639242 & 3.832926953 & 6.209819218 & 9.340622601 & 11.73344191 \\
%$T_6$   & 1.656518297 & 4.070744418 & 6.52471652  & 10.11533211 & 12.57110737 \\
%Normal & 2.849532378 & 5.928574203 & 9.162358656 & 13.60372323 & 17.31770314 \\ \bottomrule
%\end{tabular}%
%}
%\end{table}

\begin{table}[H]
\centering
\caption{Quantiles of asymptotic distribution under $\ell_1$ criterion using $\mu = 2$}
\label{tab:l1_two}
%\resizebox{\textwidth}{!}{%
\begin{tabular}{@{}llllll@{}}
\toprule
Distributions       & 90\%        & 95\%        & 97.50\%     & 99\%        & 99.50\%     \\ \midrule
$T_3$   & 1.7 & 4.8 & 7.5 & 11.1 & 13.8 \\
$T_4$   & 2.7 & 5.7 & 8.5 & 12.3 & 15.3 \\
$T_5$   & 3.1 & 6.1 & 9.0 & 13.0 & 16.0 \\
$T_6$   & 3.3 & 6.3 & 9.3  & 13.2 & 16.4 \\
$T_{10}$   & 3.6 & 6.7   & 9.8 & 14.0  & 17.2 \\
$T_{15}$   & 3.8 & 6.9   & 10.0 & 14.3  & 17.6 \\
Normal & 4.0 & 7.2 & 10.3 & 14.7 & 18.1 \\ \bottomrule
\end{tabular}%
%}
\end{table}

%\begin{table}[H]
%\centering
%\caption{Quantiles of asymptotic distribution under $\ell_2$ loss using $\mu = 2$}
%\label{tab:l2_two}
%\resizebox{\textwidth}{!}{%
%\begin{tabular}{@{}llllll@{}}
%\toprule
%Distributions       & 90\%        & 95\%        & 97.50\%     & 99\%        & 99.50\%     \\ \midrule
%$T_3$   & 0.719635921 & 2.792342472 & 5.137316634 & 8.805166323 & 11.7874725  \\
%$T_4$   & 1.257873137 & 3.630015237 & 6.140880704 & 9.842717575 & 13.01604768 \\
%$T_5$   & 1.537051033 & 4.004189541 & 6.736467324 & 10.33396216 & 13.44256397 \\
%$T_6$   & 1.666643461 & 4.187791951 & 6.948644662 & 10.63252291 & 13.46787624 \\
%Normal & 2.468869987 & 5.514223902 & 8.580450009 & 12.78174306 & 16.09263141 \\ \bottomrule
%\end{tabular}%
%}
%\end{table}

\begin{table}[H]
\centering
\caption{Quantiles of asymptotic distribution under $\ell_2$ criterion using $\mu = 2$}
\label{tab:l2_two}
%\resizebox{\textwidth}{!}{%
\begin{tabular}{@{}llllll@{}}
\toprule
Distributions       & 90\%        & 95\%        & 97.50\%     & 99\%        & 99.50\%     \\ \midrule
$T_3$   & 2.2(+29.4\%) & 5.9(+23\%) & 9.5(+26.7\%) & 15.0(+35.1\%) & 19.8(+43.5\%)  \\
$T_4$   & 2.9(+7.4\%) & 6.2(+8.8\%) & 10.0(+17.6\%) & 14.3(+16.3\%) & 18.2(+19\%) \\
$T_5$   & 3.2(+3.2\%) & 6.4(+3.3\%) & 9.6(+6.7\%) & 14.2(+9.2\%) & 17.8(+11.25\%) \\
$T_6$   & 3.3(+0\%) & 6.4(+1.6\%) & 9.6(+3.2\%) & 14.0(+6.1\%) & 17.5(+6.7\%) \\
$T_{10}$   & 3.5(-2.8\%) & 6.6(-1.5\%)   & 9.7(-1\%) & 14.0(+0\%)  & 17.4(+1.2\%) \\
$T_{15}$   & 3.6(-5.3\%) & 6.6(-4.3\%)   & 9.7(-3\%) & 13.9(-2.8\%)  & 17.3(-1.7\%) \\
Normal & 3.7(-7.5\%) & 6.7(-7\%) & 9.7(-5.8\%) & 13.9(-5.4\%) & 17.1(-5.5\%) \\ \bottomrule
\end{tabular}%
%}
\end{table}

\noindent
From the above tables, it is immediate that if the error distribution has heavy tails, say,  $T_3, T_4$, it is preferable to use $\hat d^{\ell_1}$ to $\hat d^{\ell_2}$ as the former has tighter limiting confidence interval for any of the levels presented in our tables.  On the other hand, if the error distribution is normal, then the $\ell_2$ estimator is more efficient in terms of the width of the asymptotic confidence interval, as it is maximum likelihood estimator of $d_0$. In fact, the $\ell_2$ estimator starts becoming efficient for $T$ distributions with higher degrees of freedom as is already evident from the above tables where we see a reduction is some of the $\ell_2$ quantiles for certain values of $\mu$ with $T_6$ error, and a systematic reduction with $T_{10}$ and $T_{15}$ errors.

\section{Conclusion}
\label{sec:conclusion}
In this paper, we have analyzed various estimators in the standard change point model and its multi-dimensional analogue by minimizing HEFs, especially in the presence of heavy tailed errors. We note that the robust Huber-estimators show varying degrees of advantage over the least squares estimator, depending on the dimensionality of the problem. 
\begin{enumerate} 
\item In one dimension, all estimators achieve the same rate of convergence, whereas the limiting distributions for the robust criteria based estimators are more concentrated around $0$ than that of the least squares estimator. This effect diminishes as the tail of the error distribution becomes lighter: in particular, for normal errors the least squares estimator has a narrower asymptotic confidence interval in comparison to the robust estimators. We believe a similar phenomenon will arise in the change-plane problem for fixed $p$ (where again, all the Huber estimators and the LSE will converge at rate $n$), but the limit distributions in the multidimensional case are expected to be multidimensional analogues of compound Poisson processes with extremely involved characterizations. Almost nothing is known about these objects and their study constitutes a highly non-trivial project in its own right. 
%(as is evident from the tables presented in Section \ref{sec:simulation}), as in this case, the least squares estimator is maximum likelihood estimator when the error distribution is normal(?).      
\item In growing dimensions, the robust estimators attain faster rates of convergence than the least squares estimator, in particular attaining the minimax rate which does not depend upon the tail of the error, whilst the rate of convergence of the least squares estimator is dampened by the tail of the error distribution.  
%where the rate of convergence of the robust estimators is free from it and is minimax optimal up to a log factor. 
\end{enumerate}
We now briefly discuss some variants of the problem considered above as well as possible directions for future research.  
\subsection{Binary response model: } A natural variant of the change plane model analyzed in Subsection \ref{sec:findings_multi_dim} is the following binary response model:  
$$
X \sim P,  \ \ \bbP(Y = 1 \mid X) = \alpha_0 \mathds{1}_{X^{\top}d_0 \le 0} + \beta_0 \mathds{1}_{X^{\top}d_0 > 0} \,,
$$
where $0 < \alpha_0 \neq \beta_0 < 1$. One may minimize the squared error loss to estimate the unknown parameters: 
$$
(\hat \alpha, \hat \beta, \hat d) = \argmin_{\alpha, \beta, d} \frac1n \sum_{i=1}^n \left(Y_i - \alpha \mathds{1}_{X_i^{\top}d \le 0} - \beta\mathds{1}_{X_i^{\top}d > 0}\right)^2 = \argmin_{\alpha, \beta, d} \bbP_n f_{\alpha, \beta, d}\,.
$$ 
As in Subsection \ref{sec:findings_multi_dim}, the change plane parameter $d_0$ here is also identified up to its direction and the level parameters $(\alpha_0, \beta_0)$ are identified up to their order, so we assume $\|d_0\| = 1$ and $\alpha_0 < \beta_0$. The loss function $f_{\alpha, \beta, d}$ is uniformly bounded by $1$, hence the techniques used to prove the first part of Theorem \ref{thm:l2} yield: 
\begin{align*}
\left(\sqrt{n} \wedge \frac{n}{p}\left(\log{\frac{n}{p}}\right)^{-1}\right)\left(\hat \alpha - \alpha_0\right) = O_p(1), \  \ \left(\sqrt{n} \wedge \frac{n}{p}\left(\log{\frac{n}{p}}\right)^{-1}\right)\left(\hat \beta - \beta_0\right) & = O_p(1), \\
& \hspace{-25em} \frac{n}{p}\left(\log{\frac{n}{p}}\right)^{-1}\bbP\left(\s(X^{\top}\hat d) \neq \s(X^{\top}d_0)\right)  = O_p(1) \,.
\end{align*} 
Furthermore, the rate obtained above can be shown to be minimax optimal (up to a log factor) by following a similar line argument as in the proof of Theorem \ref{thm:mlb_growing}.

\subsection{More general regression functions: } We have analyzed in this paper a stump based change point model: 
%$$
%Y_i = \alpha_0 \mathds{1}_{X_i \le d_0} + \beta \mathds{1}_{X_i > d_0} + \xi_i
%$$
%and its high dimensional change plane counterpart. 
The model analyzed in this paper can be easily generalized to one where the levels $(\alpha_0, \beta_0)$ on either side of the boundary are replaced by some unknown functions of $X$. As an example, one may fit the following non-parametric model: 
 $$
Y_i = f(X_i) \mathds{1}_{X_i \le d_0} + g(X_i) \mathds{1}_{X_i > d_0} + \xi_i
$$
where both $f, g$ are smooth and $f(d_0) \neq g(d_0)$. One may estimate $f, g, d_0$ using the following HEF: 
$$
(\hat f^k, \hat g^k, \hat d^k) = \argmin_{f, g, d} \frac1n \sum_{i=1}^n \tilde H_k\left(Y_i - f(X_i)\mathds{1}_{X_i \le d} - g(X_i)\mathds{1}_{X_i > d}\right)\,,
$$
with $(f,g)$ restricted to an appropriate class of functions (depending upon the underlying application). This model is well investigated in the literature using the squared error loss [refer], however the properties of the robust estimators (i.e. estimators obtained by minimizing HEF) are still largely unknown and worthy of investigation in the presence of heavy tailed errors.  

\subsection{Smoothed change plane problem: }
The change plane estimators analyzed in Subsection \ref{sec:findings_multi_dim} are NP-hard to compute as  HEF is discontinuous at the change boundary. One may replace the indicator function involved in HEF by a smooth sigmoid function to estimate the unknown parameters as follows: 
$$
(\hat \alpha^k, \hat \beta^k, \hat d^k) = \argmin_{\alpha, \beta, d} \frac1n \sum_i \tilde H_k\left(Y_i - \alpha - (\beta - \alpha)\frac{e^{X_i^{\top}d_0/\sigma_n}}{1 + e^{X_i^{\top}d_0/\sigma_n}}\right) 
$$
for some bandwidth parameter $\sigma_n \to 0$ as $n \to \infty$. The sigmoid function converges to the indicator function as $n \to \infty$ and is differentiable with respect to $(\alpha, \beta, d)$, therefore one may employ gradient descent to estimate the parameters; however, as the loss function is non-convex, there is no guarantee that gradient descent type techniques initiated from a random point on the parameter surface will converge to a global minimum. One way to address this issue is to replace the indicator function by a convex surrogate (i.e. logit function as in logistic regression, exponential function as in adaboost), but as the convex function does not converge to the indicator function, it is unclear whether this method will lead to a consistent estimator of $d_0$. However, such methods merit deeper investigation as they may facilitate efficient computation of the change plane estimator.

\appendix 
\section{Proofs of selected Theorems}
For all proofs below we will assume $\alpha_0 > \beta_0$ for simplicity of presentation. The derivations all go through for the reverse inequality upon minor adjustments of the proofs presented in the paper. 
\subsection{Proof of Theorem \ref{thm:tail_bound_cpp}}
We divide the whole proofs into few supplementary lemmas, whose proofs can be found in Appendix \ref{app:supp_lemmas}. Our first lemma provides a lower bound on the the probability of a random walk staying always positive. We believe that this lemma has been proved before, but we were unable to find a proper source to cite. Hence we will provide our own proof in Appendix \ref{app:supp_lemmas}.  
\begin{lemma}
\label{lem:lower_bound_rw}
Suppose $\{S_i\}_{i \ge 0}$ is a positively drifted random walk (i.e. $S_i = \sum_{j=1}^i (X_i + \mu)$, $\bbE(X) = 0, \mu > 0$) with $S_0 = 0$. Then we have: 
$$
\bbP\left(\max_{1 \le i \le n} S_i < 0\right)  \ge \frac1n \bbP\left(S_n < 0\right)\,.
$$
\end{lemma}
The compound Poisson process is essentially a two sided random walk where are number of steps till time $t$ follows a Poisson process. Therefore, we start by establishing tail bound on the minimizer of the random walk and then relate it to the tail of minimizer of the compound Poisson process. Our next lemma establishes that if the step distribution of a random walk follows a Pareto distribution, then the minimizer of the random walk is also heavy-tailed:  
\begin{lemma}
\label{lem:lower_random_onesided}
Suppose $\xi_1, \xi_2, \dots $ i.i.d. random variables with the following distribution: 
$$
\bbP\left(|\xi| > t\right) = \frac{1}{1 + t^{\gamma}}
$$
and $\bbP(\xi > t) = 1 - \bbP(\xi \le -t)$ for all $t > 0$. Define $X_i = \xi_i + \mu$ for some $\mu > 0$ and a random walk based on $X_i$'s, i.e $S_n = \sum_{i=1}^n X_i$. Suppose $M$ denotes the minimizer of the random walk on $\bbZ^+$. Then we have: 
$$
\bbP\left(M \ge k \right) \ge \frac{c_1c_2p^*}{\gamma} \times \frac{1}{k^{\gamma}} := c_0 k^{-\gamma} \,,
$$
for all $k \ge k_0 := 1 \vee \lceil \mu^{-\gamma/(\gamma - 1)}\rceil$, where: 
\begin{enumerate}
    \item $p^* = \bbP\left(S_i > 0 \ \ \forall \ \ i \in \bbN\right) = \bbP(M = 0)$. 
    \item $c_1 = \frac{1}{2(1 + \mu^{-\gamma})(1 + \mu)^{\gamma}}$. 
    \item $c_2 = \inf_{x \ge 1} \left(1 - \frac{1}{1+x}\right)^{x-1}$. 
\end{enumerate}
\end{lemma}
The previous lemma indicates that the minimizer of the random walk with a heavy tailed step distribution is also heavy tailed. As the compound Poisson process is a two sided process (i.e. supported on entire real line), we next extend our lower bound on the tail of the minimizer of random walk obtained in previous lemma for a two sided random walk in the following lemma:

\begin{lemma}
\label{lem:two_sided_rw_minimizer}
Under the same structure as of Lemma \ref{lem:lower_random_onesided}, we consider a two sided random walk with independent component on the either side. Define by $M_{ts}$ as the minimizer of the two sided random walk and by $M_{os}$ as the minimizer of one-sided random walk. Then we have: 
$$
\bbP\left(|M_{ts}| \ge k \right) \ge 2p^*c_0k^{-\gamma}
$$
for all $k \ge k_0$, where $p^*, k_0, c_0$ are same as defined in Lemma \ref{lem:lower_random_onesided}.  
\end{lemma}
Finally, we translate the lower bound on the tail of the minimizer of the two-sided random walk to the two sided compound Poisson process in the following lemma:

\begin{lemma}
\label{cor:tsCPP}
Consider a two sided independent compound Poisson process with increment independent of the steps. More specifically, let $\{X_i\}_{i \in \bbN}$ be same as defined in Lemma \ref{lem:lower_random_onesided}. Suppose $\{X'_i\}_{i \in \bbN}$ be an independent copy of $\{X_i\}_{i \in \bbN}$. Also suppose $N_1(t)$ and $N_2(t)$ are two independent Poisson process on $\bbR^+$ with some intensity function $\Lambda(t)$. The two sided independent compound Poisson process on $\bbR$ is defined as: 
$$
X(t) = 
\begin{cases}
\sum_{i=1}^{N_1(t)}X_i \,, & \text{if } t > 0 \\
\sum_{i=1}^{N_2(-t)}X'_i \,, & \text{if } t < 0 \\
0 \,, & \text{if } t = 0 \,.
\end{cases}
$$
Let $M$ be the mid-argmin of $X(t)$ over $\bbR$. Then we have for all $x > (k_0 + \gamma + \log{2})/f_X(d_0)$: 
$$
\bbP\left(M_{ts, CPP} > x\right) \ge \frac{c_0}{2f_X^{\gamma}(d_0)}x^{-\gamma} \,,
$$
where $c_0, k_0$ are same constants as defined in Lemma \ref{lem:lower_random_onesided}.  
\end{lemma}
Combining Lemma \ref{lem:lower_random_onesided}, \ref{lem:two_sided_rw_minimizer} and \ref{cor:tsCPP}, we conclude the proof of of lower bound on $F_{\ell_2}$. 
\\\\
\noindent
Now to prove the upper bound for $F_{\ell_1}$ we modify our arguments in the previous lemmas. Note that $F_{\ell_1}$ is the distribution of the minimizer of the following compound Poisson process: 
$$
CPP(t) = \sum_{i=1}^{N_+(t)}\left(\xi^* + \mu_0\right)\mathds{1}_{t \ge 0} +  \sum_{i=1}^{N_-(-t)}\left(\xi^* + \mu_0\right)\mathds{1}_{t < 0}
$$
with $CPP(t) = 0$, $N_+$ and $N_-$ are two independent Poisson processes as before and: 
$$
\xi^* \leftarrow \left\{\left|\xi + (\alpha_0 - \beta_0)\right| - |\xi|\right\} -  \bbE\left[\left|\xi + (\alpha_0 - \beta_0)\right| - |\xi| \right]
$$
with $\mu_0 = \bbE\left[\left|\xi + (\alpha_0 - \beta_0)\right| - |\xi| \right] > 0$. Before going into the details of the proof, we state Hoeffding's inequality bound: 
\begin{align*}
\bbP\left(S_n < 0\right) & = \bbP\left(\sum_{i=1}^n \xi^*_i < -n\mu_0\right) \\
& = \bbP\left(\bar \xi^*_n < -\mu_0\right) \le e^{-\frac{n\mu_0^2}{8(\alpha_0 - \beta_0)^2}} \,.
\end{align*}
Here we will highlight the steps where a modification is needed. First note that, in case of one sided random walk (same situation as in Lemma \ref{lem:lower_random_onesided}) we obtain the upper bound as follows: 
\allowdisplaybreaks
\begin{align}
    P\left(M_{os} \ge k\right) & = \sum_{j \ge k} \bbP\left(M_{os} = j\right) \notag \\
    & = p^* \sum_{j \ge k} \bbP\left(\max_{1 \le i \le j} S_i < 0 \right)  \notag \\
    & \le  p^* \sum_{j \ge k} \bbP\left(S_j < 0\right) \notag \\
    & \le p^* \sum_{j \ge k} e^{-\frac{j\mu_0^2}{8(\alpha_0 - \beta_0)^2}} \notag \\
    & = \frac{p^*}{1 - e^{-\frac{\mu_0^2}{8(\alpha_0 - \beta_0)^2}}} e^{-\frac{k\mu_0^2}{8(\alpha_0 - \beta_0)^2}} 
\end{align}
We now translate the tail bound on the minimizer of the one sided random walk to the minimizer of the two sided random walk as below: 
\begin{align*}
    \bbP\left(M_{ts} = k\right) & = \bbP\left(S_K \le S_i \ \forall \ 0 \le i \le k-1, S_k \le S_i \ \forall \ k+1 \le i < \infty, S_k \le \inf_{j \ge 1}S_{-j}\right) \\
    & \le \bbP\left(S_K \le S_i \ \forall \ 0 \le i \le k-1, S_k \le S_i \ \forall \ k+1 \le i < \infty \right) \\
    & =  \bbP\left(M_{os} = k\right) \,.
\end{align*}
This, along with the upper bound on the tail of the one-sided random walk implies: 
$$
\bbP\left(M_{ts} \ge k\right) \le  \frac{p^*}{1 - e^{-\frac{\mu_0^2}{8(\alpha_0 - \beta_0)^2}}} e^{-\frac{k\mu_0^2}{8(\alpha_0 - \beta_0)^2}} \,.
$$
Next, we translate the bound for the minimizer of a one-sided random walk to a one-sided compound Poisson process with steps $\xi^* + \mu_0$: 
\begin{align*}
    P(M_{os, CPP} > x) & = \sum_{k = 0}^{\infty} \bbP\left(M_{os, CPP} > x \mid N_1(x) = k\right) \bbP\left(N_1(x) = k\right) \\
    & = \sum_{k = 0}^{\infty} \bbP\left(\argmin_{i \ge 0} S_i > k\right) \bbP\left(N_1(x) = k\right) \\
    & \le  \frac{p^*}{1 - e^{-\frac{\mu_0^2}{8(\alpha_0 - \beta_0)^2}}} \sum_{k=0}^{\infty}e^{-\frac{(k+1)\mu_0^2}{8(\alpha_0 - \beta_0)^2}} \bbP\left(N_1(x) = k\right) \\
    & =  \frac{p^*}{1 - e^{-\frac{\mu_0^2}{8(\alpha_0 - \beta_0)^2}}} e^{-\frac{\mu_0^2}{8(\alpha_0 - \beta_0)^2}}\sum_{k=0}^{\infty}e^{-\frac{k\mu_0^2}{8(\alpha_0 - \beta_0)^2}} \frac{e^{-\Lambda(x)}\Lambda(x)^k}{k!} \\
    & =  \frac{p^*}{1 - e^{-\frac{\mu_0^2}{8(\alpha_0 - \beta_0)^2}}} e^{-\frac{\mu_0^2}{8(\alpha_0 - \beta_0)^2}}e^{-\Lambda(x)} \sum_{k=0}^{\infty}\frac{\left(e^{-\frac{\mu_0^2}{8(\alpha_0 - \beta_0)^2}} \Lambda(x)\right)^k}{k!} \\
    & = \frac{p^*}{e^{\frac{\mu_0^2}{8(\alpha_0 - \beta_0)^2}} - 1}\exp{\left(-\Lambda(x)\left(1 - e^{-\frac{\mu_0^2}{8(\alpha_0 - \beta_0)^2}}\right)\right)} \\
    & = \frac{p^*}{e^{\frac{\mu_0^2}{8(\alpha_0 - \beta_0)^2}} - 1}\exp{\left(-xf_X\left(d_0\right)\left(1 - e^{-\frac{\mu_0^2}{8(\alpha_0 - \beta_0)^2}}\right)\right)} \,,\\
\end{align*}
and consequently for the two sided compound Poisson process: 
\begin{align*}
    \bar F_{\ell_1}(x) = \bbP\left(M_{ts, CPP} > x\right) & = \sum_{k=0}^{\infty}\bbP\left(M_{ts, CPP} > x \mid N_1(x) = k\right)\bbP(N_1(x) = k) \\
    & = \sum_{k=0}^{\infty}\bbP\left(M_{ts} > k \right)\bbP(N_1(x) = k)  \\
    & \le \sum_{k=0}^{\infty}\bbP\left(M_{os} > k \right)\bbP(N_1(x) = k) \\
    & = \bbP\left(M_{os, CPP} > x\right) \\
    & \le \frac{p^*}{e^{\frac{\mu_0^2}{8(\alpha_0 - \beta_0)^2}} - 1}\exp{\left(-xf_X\left(d_0\right)\left(1 - e^{-\frac{\mu_0^2}{8(\alpha_0 - \beta_0)^2}}\right)\right)} \,.
\end{align*}
\subsection{Proof of Theorem \ref{thm:parallel_change_point}}
\paragraph{Proof of lower bound: }The following lemma, which is a finite sample analogue of the first conclusion of Theorem \ref{thm:tail_bound_cpp}, is essential to establish the lower bound of Theorem \ref{thm:parallel_change_point}: 
\begin{lemma}
\label{lem:finite_sample_tail_bound}
Suppose, for a fixed $n$, $F_{n, \ell_2}$ denotes the distribution of $n(\hat d^{\ell_2}_i - d_{0, i})$. Then we have for all $2\gamma/f_X(d_0) \le |x| \le (\delta_1 \wedge \delta_2)n$ (for some constants $\delta_1, \delta_2$ independent of $n$ defined explicitly in the proof): 
$$
1 - F_{n, \ell_2}(x) =   \bbP\left(\left| n(\hat d^{\ell_2}_i - d_{0, i})\right| \ge x \right) \ge \frac{c_1c_2(p^*)^2}{\gamma 2^{\gamma + 2}}  \times \left(\frac{f_{X, \max}}{1 - F_X(d_0)}\right)^{-\gamma} \times x^{-\gamma}
$$ 
where $f_{X, \max}$ is the maximum value of the density of $X$ and $c_1, c_2, p^*$ are same as defined in Lemma \ref{lem:lower_random_onesided}. 
\end{lemma}
The proof of the Lemma can be found in Appendix \ref{app:supp_lemmas}. For notational simplicity, set: 
$$
C = \frac{c_1c_2(p^*)^2}{\gamma 2^{\gamma + 2}}  \times \left(\frac{f_{X, \max}}{1 - F_X(d_0)}\right)^{-\gamma} \,.
$$
Using the above lemma we have: 
\begin{align*}
    \bbP\left(\max_{1 \le i \le m}\frac{n}{m^{1/\gamma}}\left|\hat d^{\ell_2}_i - d_{0, i}\right| > t \right) & = 1 - \bbP\left(\max_{1 \le i \le m}\frac{n}{m^{1/\gamma}}\left|\hat d^{\ell_2}_i - d_{0, i}\right| \le t \right) \\
    & = 1 - \left(F_{n, \ell_2}(tm^{1/\gamma})\right)^m \\
     & = 1 - \left(1 -  \bar F_{n, \ell_2}(tm^{1/\gamma})\right)^m \\
    & \ge 1 - \left(1 - C(tm^{1/\gamma})^{-\gamma}\right)^m \\
     & = 1 - \left(1 -Cm^{-1}t^{-\gamma}\right)^m \\
     & \longrightarrow 1 - e^{-Ct^{-\gamma}}
\end{align*}
Note that Lemma \ref{lem:finite_sample_tail_bound} is applicable here as for any fixed $t$, $tm^{1/\gamma} \ll n$ because $m^{1/\gamma} \ll n$ and as $m \uparrow \infty$, $tm^{1/\gamma} \ge  2\gamma/f_X(d_0)$ for all large $m$. This completes the proof.

\paragraph{Proof of upper bound:} The proof of upper bound relies on the following Lemma, which is a analogue of Lemma \ref{lem:finite_sample_tail_bound}, where we establish an upper bound on the finite sample distribution of $n(\hat d_i - d_{i, 0})$ with bounded supported error distribution $\xi$: 
\begin{lemma}
\label{lem:finite_sample_tail_bound_upper}
Let $F_{n, \ell_1}$ denotes the distribution of $\left|n(\hat d^{\ell_1}_i - d_{0, i})\right|$. Then we have for $0 \le |x| \le n\delta_1$ (for some constant $\delta_1$ defined explicitly in the proof): 
$$
1 - F_{n, \ell_1}(x) = \bbP\left(\left| n(\hat d^{\ell_1}_i - d_{0, i})\right| \ge x \right) \le \frac{2e^{-c}}{1 - e^{-c}}  e^{-x\frac{f_X(d_0)}{2}(1 - e^{-c})} \,,
$$ 
where $c = \mu^2/4b^2$, $\mu = \bbE[|\xi + (\alpha_0 - \beta_0)| - |\xi|]$ and $b$ is the range of the random variable $\left(|\xi + (\alpha_0 - \beta_0)| - |\xi|\right) - \mu$. 
\end{lemma}
\noindent
Using the above lemma we have: 
\begin{align*}
\bbP\left(\max_{1 \le i \le m}\frac{n}{\log{m}}\left|\hat d^{\ell_1}_i - d_{0, i}\right| > t \right) & \le \sum_{i=1}^m \bbP\left(n\left|\hat d^{\ell_1}_i - d_{0, i}\right| > t\log{m} \right) \\
 & \le  \frac{2e^{-c}}{1 - e^{-c}} m e^{-t\log{m}\frac{f_X(d_0)}{2}(1 - e^{-c})} \\
 & \le \frac{2e^{-c}}{1 - e^{-c}} e^{-\log{m}\left(t\frac{f_X(d_0)}{2}(1 - e^{-c}) - 1\right)} \,.
\end{align*}
This completes the proof.

\subsection{Proof of Theorem \ref{thm:mlb_parallel}}
To prove the lower bound, we consider a simple model: Assume that, for each problem, the true change point is $0$ (i.e. $d_{i, 0} = 0$ for all $i$), the covariates $X_{i, j}'s$ are all i.i.d. Uniform $(-1, 1)$ and error distribution is normal. We first observe that for any estimator $\hat d_i$ of $d_{i, 0}$ we have: 
$$
\left| \hat d_i - d_{i, 0}\right| \ge \min_{1 \le j \le n} \left|X_{i, j} - d_{i, 0}\right| \,,
$$
i.e. we can't estimate a change point better than its closest order statistic. Note that when $d_{i, 0} = 0$, we have:
$$
\left| \hat d_i\right| \ge \min_{1 \le j \le n} \left|X_{i, j}\right|  =  \min_{1 \le i \le n} U_i 
$$
where $U_1, \dots, U_n \overset{i.i.d.}{\sim} U(0, 1)$. Hence to prove the theorem, all we need to show is that: 
$$
\liminf_{n, m \to \infty} \frac{1}{\log{m}} \bbE\left[\max_{1 \le i \le m} n Z_{i:n}\right] \ge C > 0
$$
where $Z_{i:n}$ are i.i.d with the common distribution being that of the minimum of $n$ uniform $(0, 1)$ random variables. Note that for any $0 \le t \le n$:  
\begin{align*}
\bbP\left(n Z_{i:n} \ge t\right) & = \bbP\left(\min_{1 \le i \le n} U_i \ge \frac{t}{n}\right) \\
& = \left(1 - \frac{t}{n}\right)^n \\
\implies \bbP\left(n Z_{i:n} \le t\right) & = 1 - \left(1 - \frac{t}{n}\right)^n  := F_{nZ_{1:n}}(t) \,.
\end{align*} 
Therefore we have: 
\begin{align}
\frac{1}{\log{m}} \bbE\left[\max_{1 \le i \le m} n Z_{i:n}\right] & = \frac{1}{\log{m}} \int_0^n \bbP\left(\max_{1 \le i \le m} n Z_{i:n} \ge t \right) \ dt \notag \\
& = \frac{1}{\log{m}} \int_0^n \left[1 - \bbP\left(\max_{1 \le i \le m} n Z_{i:n} \le t \right)\right] \ dt \notag \\
& = \frac{1}{\log{m}} \int_0^n \left[1 - F^m_{nZ_{1:n}}(t)\right] \ dt \notag \\
\label{eq:mlb_parallel_1} & = \frac{1}{\log{m}} \int_0^n \left[1 -\left(1 - \left(1 - \frac{t}{n}\right)^n\right)^m\right] \ dt 
\end{align}
Next using the following inequality:
$$
\left(1 + \frac{t}{n}\right)^n \ge e^t\left(1 - \frac{t^2}{n}\right) \ \ \ \forall \ \ \ |t| \le n \,,
$$
we obtain from equation \eqref{eq:mlb_parallel_1}: 
\begin{align}
\frac{1}{\log{m}} \bbE\left[\max_{1 \le i \le m} n Z_{i:n}\right]  & \ge \frac{1}{\log{m}} \int_0^n \left[1 -\left(1 - e^{-t} + e^{-t}\frac{t^2}{n}\right)^m\right] \ dt  \notag \\
& = \frac{1}{\log{m}} \int_0^n \left[1 - \sum_{i=0}^m \dbinom{m}{i} \left(e^{-t}\frac{t^2}{n}\right)^i\left(1 - e^{-t}\right)^{m - i} \right]\notag \\
& = \frac{1}{\log{m}} \int_0^n \left[1 -\left(1 - e^{-t}\right)^{m}\right] \ dt \notag \\
& \qquad \qquad -  \frac{1}{\log{m}} \int_0^n \sum_{i=1}^m \dbinom{m}{i} \left(e^{-t}\frac{t^2}{n}\right)^i\left(1 - e^{-t}\right)^{m - i} \ dt \notag \\
& =  \frac{1}{\log{m}} \int_0^\infty \left[1 -\left(1 - e^{-t}\right)^{m}\right] \ dt \notag \\
& \qquad \qquad - \frac{1}{\log{m}} \int_n^\infty \left[1 -\left(1 - e^{-t}\right)^{m}\right] \ dt \notag \\
& \qquad \qquad \qquad \qquad - \frac{1}{\log{m}} \int_0^n \sum_{i=1}^m \dbinom{m}{i} \left(e^{-t} \frac{t^2}{n}\right)^i\left(1 - e^{-t}\right)^{m - i} \ dt \notag \\
\label{eq:mlb_parallel_2}  & := a_{m, 1}  - a_{m, 2} - a_{m, 3}  
\end{align}
We next show that $a_{m, 2}, a_{m, 3} \longrightarrow 0$ and $a_{m, 1} \longrightarrow 1$ as $m \to \infty$, as long as $n \ge 2m$. We start with $a_{m, 3}$: 
\begin{align*}
a_{m, 3} & = \frac{1}{\log{m}} \int_0^n \sum_{i=1}^m \dbinom{m}{i} \left(e^{-t} \frac{t^2}{n}\right)^i\left(1 - e^{-t}\right)^{m - i} \ dt \\
& \le  \frac{1}{\log{m}} \int_0^n \sum_{i=1}^m \dbinom{m}{i} \left(e^{-t} \frac{t^2}{n}\right)^i \ dt \\
& = \frac{1}{\log{m}}\sum_{i=1}^m \dbinom{m}{i}\int_0^n \left(e^{-t} \frac{t^2}{n}\right)^i \ dt \\
& = \frac{1}{\log{m}}\sum_{i=1}^m  \dbinom{m}{i} \frac{1}{i n^i} \int_0^n t^{2i} \ ie^{-it}  \ dt \\
& \le \frac{1}{\log{m}}\sum_{i=1}^m  \dbinom{m}{i} \frac{1}{i n^i} \int_0^\infty t^{2i} \ ie^{-it}  \ dt \\
& = \frac{1}{\log{m}}\sum_{i=1}^m  \dbinom{m}{i} \frac{1}{i n^i} \frac{(2i)!}{i^{2i}} \\
& =  \frac{1}{\log{m}}\sum_{i=1}^m  \frac{(m-i+1)\cdots m}{i n^i i^i}\frac{(i+1) \cdots (i+i)}{i^i} \\
& \le \frac{1}{\log{m}}\sum_{i=1}^m  \frac{(m-i+1)\cdots m}{i \left(\frac{n}{2}\right)^i i^i} \\
& \le \frac{1}{\log{m}}\sum_{i=1}^m  \left(\frac{2m}{n}\right)^i \frac{1}{i^{i+1}} \le \frac{1}{\log{m}}\sum_{i=1}^m  \frac{1}{i^{i+1}} \longrightarrow 0 \,.
\end{align*}
For the other term $a_{n, 2}$: 
\begin{align*}
\left|a_{n, 2}\right| & =  \left|\frac{1}{\log{m}} \int_n^\infty \left[1 -\left(1 - e^{-t}\right)^{m}\right] \ dt\right| \\
& =  \left|\frac{1}{\log{m}} \int_n^\infty \left[1 - \sum_{i=0}^m \dbinom{m}{i} e^{-it}(-1)^i \right] \ dt\right| \\
& = \left|\frac{1}{\log{m}} \int_n^\infty  \sum_{i=1}^m \dbinom{m}{i} e^{-it}(-1)^{i+1}  \ dt \right| \\
& \le  \frac{1}{\log{m}}  \sum_{i=1}^m \dbinom{m}{i} \int_n^\infty e^{-it}  \ dt \\
& =  \frac{1}{\log{m}}  \sum_{i=1}^m \dbinom{m}{i}\frac{e^{-ni}}{i} \\
& =  \frac{2^m}{\log{m}} \bbE\left[\frac{e^{-nX}}{X}\mathds{1}_{X \ge 1}\right] \hspace{0.2in} [X \sim \text{Bin}(n, p)]\\
& \le \frac{2^m}{\log{m}} \bbE\left[\frac{e^{-2mX}}{X}\mathds{1}_{X \ge 1}\right] \hspace{0.2in} [\because n \ge 2m] \\ 
& \le \frac{2^m}{\log{m}} \bbE\left[e^{-2mX}\right] \hspace{0.2in} [\because n \ge 2m] \\
& = \frac{2^m}{\log{m}} \left(\frac12 e^{-2m} + \frac12\right)^m \\ 
& = \frac{1}{\log{m}} \left(e^{-2m} + 1 \right)^m  \longrightarrow 0\,. 
\end{align*}
Now the calculation of $a_{n, 1}$ is similar by replacing $n$ with $0$. We have: 
\begin{align*}
a_{n, 1} & = \frac{1}{\log{m}} \int_0^\infty \left[1 -\left(1 - e^{-t}\right)^{m}\right] \ dt \\
& = \frac{1}{\log{m}} \int_0^\infty  \sum_{i=1}^m \dbinom{m}{i} e^{-it}(-1)^{i+1}  \ dt \\
& =  \frac{1}{\log{m}}  \sum_{i=1}^m \dbinom{m}{i} \frac{(-1)^{i+1}}{i} \\
& = \frac{1}{\log{m}}  \sum_{i=1}^m \frac{1}{i} \overset{m \to \infty}{\longrightarrow} 1 \,.
\end{align*} 
%\begin{align*}
%a_{n, 1} & =  \frac{1}{\log{m}}  \sum_{i=1}^m \frac{1}{i} \longrightarrow 1 \,.
%\end{align*}
where the last equality follows from the representation of Harmonic number. Therefore from equation \eqref{eq:mlb_parallel_2} we conclude: 
$$
\liminf_{m, n \to \infty} \frac{1}{\log{m}} \bbE\left[\max_{1 \le i \le m} n Z_{i:n}\right] \ge 1 \,.
$$
This concludes the proof.

\subsection{Proof of Theorem \ref{thm:l2}}
\subsubsection{Case 1: $0 \le k < \infty$}
We first establish the rate of convergence of $(\hat \alpha_{\init}, \hat \beta_{\init}, \hat d)$. Towards that direction, we use the following semi-metric over the parameter space $\Theta$: 
$$
\dist(\theta_1, \theta_2) = \sqrt{\left(\alpha_1 - \alpha_2\right)^2 + \left(\beta_1 - \beta_2\right)^2 + \bbP\left(\s(X^{\top}d) \neq \s(X^{\top}d_0)\right)}
$$
The curvature of the population score function $\bbM(\theta)$ around the its value at minimizer $\bbM(\theta_0)$ is obtained via the similar calculation as in in the proof of Theorem \ref{thm:asymptotic_huber} (specifically equation \eqref{eq:pop_centered_huber}). Consider all $\theta \in \Theta$ such that $\dist(\theta, \theta_0) \le \delta$ where $\delta$ is such that $|\alpha_0 - \beta_0| > 2\delta$. For such that $\theta$ we have: 
%$$
%\max\left\{\left| \alpha - \alpha_0\right|, \left|\beta - \beta_0\right|, \left\|d - d_0\right\| \right\} \le \delta
%$$
%for some small $\delta > 0$ such Assumption \ref{assm:wedge} is satisfied and it is ensured that $|\alpha_0 - \beta_0| \ge 2\delta$. 
%Following the same calculation as in equation \eqref{eq:pop_centered_huber} in the proof of Theorem \ref{thm:asymptotic_huber} we have: 
\begin{align}
\bbM(\theta) - \bbM(\theta_0) & = \bbE\left[\left(\tilde H_k\left(\xi_i +  \alpha_0 - \alpha\right) - \tilde H_k(\xi_i)\right)\right]\bbP\left(X^{\top}d \vee X^{\top}d_0 \le 0\right) \notag \\
& \qquad + \bbE\left[\left(\tilde H_k\left(\xi_i +  \alpha_0 - \beta\right) - \tilde H_k(\xi_i)\right)\right]\bbP\left(X^{\top}d_0 < 0 < X^{\top}d\right) \notag \\
& \qquad \qquad +  \bbE\left[\left(\tilde H_k\left(\xi_i +  \beta_0 - \alpha\right) -\tilde H_k(\xi_i)\right)\right]\bbP\left(X^{\top}d < 0 < X^{\top}d_0\right) \notag \\
& \qquad \qquad \qquad +  \bbE\left[\left(\tilde H_k\left(\xi_i +  \beta_0 - \beta\right) - \tilde H_k(\xi_i)\right)\right]\bbP\left(X^{\top}d \wedge X^{\top}d_0 > 0\right) \notag \\
& \ge \frac{C_k}{2}\left[(\alpha_0 - \alpha)^2\bbP\left(X^{\top}d \vee X^{\top}d_0 \le 0\right) ) + (\beta_0 - \beta)^2\bbP\left(X^{\top}d \wedge X^{\top}d_0 > 0\right) \right. \notag \\
 \label{eq:pop_centered_huber_highdim} & \hspace{10em} \left. + \bbP\left(\s(X^{\top}d) \neq \s(X^{\top}d_0)\right)\left\{2(\alpha_0 - \beta_0 - \delta)^2 \right\} \right] \hspace{0.2in} [\text{Lemma }\ref{lem:huber_lower_bound}]\notag \\
 & \ge C_k \left[(\alpha_0 - \alpha)^2 +  (\beta_0 - \beta)^2 + \bbP\left(\s(X^{\top}d) \neq \s(X^{\top}d_0)\right)|\right] \notag \\
 & = C_k \dist^2(\theta, \theta_0) \,.
\end{align}

\paragraph{Consistency: } 
\label{para:consistency_huber}
We use argmin continuous mapping theorem (Theorem 3.2.2 of \cite{vdvw96}) to establish the consistency of the initial estimator. As the parameter space is bounded, our estimates are by default tight. As the process $\bbM(\theta) - \bbM(\theta_0)$ is continuous with respect to $\theta$ and has a clear minima at $\theta = \theta_0$ all we need to show for any compact subset $K \subseteq \Theta$: 
$$
\sup_{\theta \in K} \left|\left(\bbM_n(\theta) - \bbM_n(\theta_0)\right) - \left(\bbM(\theta) - \bbM(\theta_0)\right)\right| = o_p(1) \,.
$$
Consider a collection of functions $\cF = \left\{f_{\theta}: \theta \in \Theta\right\}$ where the individual functions $f_{\theta}(X, \xi)$ is defined as: 
\begin{align*}
f_{\theta}(X, \xi) & = \left(\tilde H_k\left(\xi +  \alpha_0 - \alpha\right) -  \tilde H_k(\xi)\right)\mathds{1}_{X^{\top}d \vee X^{\top}d_0 \le 0} \\
& \qquad \left(\tilde H_k\left(\xi +  \alpha_0 - \beta\right) - \tilde H_k(\xi)\right)\mathds{1}_{X^{\top}d_0 \le 0 < X^{\top}d } \\
& \qquad \qquad +  \left(\tilde H_k\left(\xi +  \beta_0 - \alpha\right) -\tilde H_k(\xi)\right)\mathds{1}_{X^{\top}d \le 0 < X^{\top}d_0} \\
& \qquad \qquad \qquad +   \left(\tilde H_k\left(\xi +  \beta_0 - \beta\right) - \tilde H_k(\xi)\right)\mathds{1}_{X^{\top}d \wedge X^{\top}d_0 > 0} \\
& := \sum_{i=1}^4 g^i_{\theta}(\xi)h^i_{\theta}(X)
\end{align*}
with: 
\begin{align*}
& g^1_{\theta}(\xi) = \left(\tilde H_k\left(\xi +  \alpha_0 - \alpha\right) -  \tilde H_k(\xi)\right), & h^1_{\theta}(X) = \mathds{1}_{X^{\top}d \vee X^{\top}d_0 \le 0} \,, \\
& g^2_{\theta}(\xi) = \left(H_k\left(\xi +  \alpha_0 - \beta\right) - \tilde H_k(\xi)\right), & h^2_{\theta}(X) = \mathds{1}_{X^{\top}d_0 \le 0 < X^{\top}d } \,, \\
& g^3_{\theta}(\xi) =\left(\tilde H_k\left(\xi +  \beta_0 - \alpha\right) -\tilde H_k(\xi)\right), & h^3_{\theta}(X) = \mathds{1}_{X^{\top}d \le 0 < X^{\top}d_0} \,, \\
& g^4_{\theta}(\xi) =\left(\tilde H_k\left(\xi +  \beta_0 - \beta\right) - \tilde H_k(\xi)\right), & h^4_{\theta}(X) =\mathds{1}_{X^{\top}d \wedge X^{\top}d_0 > 0} \,.
\end{align*}
As the Huber function $H_k$ is Lipschitz with Lipschitz constant $k$, our criterion function $\tilde H_k$ is Lipschitz with Lipschitz constant $(k+1)$. As our parameter space is compact, the functions $\{g^i_{\theta}, h^i_{\theta}\}_{1 \le i \le 4}$ are uniformly bounded, and has constant envelope, say $F$. That the functions $\{g^i_{\theta}\}_{\theta \in \Theta}$ for $i = 1, 2, 3, 4$ has finite VC dimension $v$ (i.e. does not grow with $n$ or $p$) is immediate. On the other hands, as all the $p$-dimensional hyperplanes passing through origin has VC dimension $p$. Hence the functions $\{h^i_{\theta}\}_{\theta \in \Theta}$ has VC dimension $p$. Define $\cF_{g, i} = \{g^i_{\theta}: \theta \in \Theta\}$ for $1 \le i \le 4$ and $\cF_{h, i} = \{h_{\theta}^i: \theta \in \Theta\}$ for $1 \le i \le 4$. Combining these we obtain: 
\begin{align*}
\sup_{Q} N(\eps\|F\|_{Q, 1}, \cF, L_1(Q)) & \le \sup_{Q} N\left(\eps\|F\|_{Q, 1}, \sum_{i=1}^4 \cF_{g, i}\cF_{h, i}, L_1(Q)\right) \\
& \le \Pi_{i=1}^4 \sup_{Q} N(\eps\|F\|_{Q, 1}, \cF_{g, i}\cF_{h, i}, L_1(Q)) \\
& \le  \Pi_{i=1}^4 KVC(\cF_{g, i})VC(\cF_{h, i})(16e)^{VC(\cF_{g, i}) + VC(\cF_{h, i})} \\
& \qquad \qquad \qquad \qquad \qquad \qquad \times \left(\frac{1}{\eps}\right)^{(VC(\cF_{g, i}) + VC(\cF_{h, i})- 2)} \\
& \le  \Pi_{i=1}^4 Kvp(16e)^{v+p}\left(\frac{1}{\eps}\right)^{(v+p- 2)} \\
&  = K^4 (16e)^{4(v+p)}\left(\frac{1}{\eps}\right)^{4(v+p-2)} \\
& = K^4 \left(\frac{16e}{\eps}\right)^{4(v+p)}
\end{align*}
This along with the fact $p/n \to 0$ implies: 
$$
\frac1n \log{\left(\sup_{Q} N(\eps\|F\|_{Q, 1}, \cF, L_1(Q)))\right)} \to 0 \,.
$$
Therefore $\cF$ is Glivenko-Cantelli class of functions and using Theorem 2.4.3 of \cite{vdvw96} we conclude that: 
\begin{align*}
\sup_{\theta \in K} \left|\left(\bbM_n(\theta) - \bbM_n(\theta_0)\right) - \left(\bbM(\theta) - \bbM(\theta_0)\right)\right| & = \left\|\left(\bbP_n - P\right)\right\|_{\cF} = o_p(1) \,.
\end{align*}
This establishes the consistency of $(\hat \alpha_{\init}, \hat \beta_{\init}, \hat d)$.

\paragraph{Rate of convergence of initial estimators: }So far we have established the quadratic curvature of $\bbM(\theta)$ around its unique minimizer $\theta_0$ and also the consistency of $\hat \theta_{\init} = (\hat \alpha_{\init}, \hat \beta_{\init}, \hat d)$. In this section we show that: 
$$
\sqrt{\frac{n}{p}}\left(\log{\frac{n}{p}}\right)^{-\frac12}\dist(\hat \theta_{\init}, \theta_0) = O_p(1)
$$
which (along with Assumption \ref{assm:wedge}) implies: 
$$
\frac{n}{p}\left(\log{\frac{n}{p}}\right)^{-1}\left[\left(\hat \alpha_{\init} - \alpha_0\right)^2 + \left(\hat \beta_{\init} - \beta_0\right)^2 + \left\|\hat d - d_0\right\| \right] = O_p(1) \,.
$$
To establish the rate, all is left to do is to find a bound on the modulus of continuity of the empirical process $\bbM_n(\theta)$, i.e we need to find $\phi_n(\delta)$ such that: 
\begin{equation}
\label{eq:moc_highdim}
\bbE\left[\sup_{d(\theta, \theta_0) \le \delta} \left| (\bbM_n - \bbM)(\theta_0) - (\bbM_n- \bbM)(\theta) \right|\right] \le \frac{\phi_n(\delta)}{\sqrt{n}} \,.
\end{equation}
Towards that end, define a local collection of functions $\cF_{\delta} = \left\{f_{\theta}: d(\theta, \theta_0) \le \delta\right\}$. Note that when $d(\theta, \theta_0) \le \delta$ we have: 
$$
\max\left\{|\alpha - \alpha_0|, |\beta - \beta_0|, \sqrt{\bbP\left(\s(X^{\top}d) \neq \s(X^{\top}d_0)\right)}\right\} \le \delta \,.
$$ 
For any such $\theta$ we have: 
\begin{align*}
\bbE\left[f_{\theta}(X, \xi)^2\right] & = \bbE\left(\tilde H_k\left(\xi +  \alpha_0 - \alpha\right) -  \tilde H_k(\xi)\right)^2\bbP\left(X^{\top}d \vee X^{\top}d_0 \le 0\right) \\
& \qquad \bbE\left(\tilde H_k\left(\xi +  \alpha_0 - \beta\right) - \tilde H_k(\xi)\right)^2\bbP\left(X^{\top}d_0 \le 0 < X^{\top}d\right) \\
& \qquad \qquad +  \bbE\left(\tilde H_k\left(\xi +  \beta_0 - \alpha\right) -\tilde H_k(\xi)\right)^2\bbP\left(X^{\top}d \le 0 < X^{\top}d_0\right) \\
& \qquad \qquad \qquad +  \bbE\left(\tilde H_k\left(\xi +  \beta_0 - \beta\right) - \tilde H_k(\xi)\right)^2\bbP\left(X^{\top}d \wedge X^{\top}d_0 > 0\right) \\
& \lesssim C_k\left[\delta^2 + \bbP\left(\s(X^{\top}d) \neq \s(X^{\top}d_0)\right)\right] \\
& \lesssim C_k \delta^2 
\end{align*}
Hence applying Theorem 8.7 of \cite{sen2018gentle} we conclude: 
$$
\bbE\left[\sup_{d(\theta, \theta_0) \le \delta} \left| (\bbM_n - \bbM)(\theta_0) - (\bbM_n- \bbM)(\theta) \right|\right] \lesssim  \sqrt{\frac{p}{n}}\delta\sqrt{\log{\frac{1}{\delta}}} \vee \frac{p}{n}\log{\frac{1}{\delta}}
$$
Therefore a valid choice of $\phi_n$ in equation \eqref{eq:moc_highdim} is: 
$$
\phi_n(\delta) = \sqrt{p}\delta\sqrt{\log{\frac{1}{\delta}}} \vee \frac{p}{\sqrt{n}}\log{\frac{1}{\delta}} \,.
$$
Using this $\phi_n$ in Theorem 3.4.1 of \cite{vdvw96} we conclude that: 
$$
\sqrt{\frac{n}{p}}\left(\log{\frac{n}{p}}\right)^{-\frac12}d(\hat \theta_{\init}, \theta_0) = O_p(1) \,.
$$
Finally, as the function class under consideration $\cF = \{f_{\theta}: \theta \in \Omega \times S^{p-1}\}$ is uniformly bounded, an application of Theorem 2 of \cite{massart2006risk} yields:   
\begin{equation}
\label{eq:concentration_bound}
\bbP\left(\frac{n}{p}\left(\log{\frac{n}{p}}\right)^{-1} d(\hat \theta_{\init}, \theta_0) \ge t\right) \le C_k e^{-c_k t}
\end{equation}
for some constants $C_k, c_k > 0$ which depends on $k$. This in particular implies that: 
\begin{equation}
\label{eq:exp_bound}
\bbE\left[\frac{n}{p}\left(\log{\frac{n}{p}}\right)^{-1} \bbP\left(\s(X^{\top}\hat d) \neq \s(X^{\top}d_0)\right)\right] \le C_k
\end{equation}
for some constant $C_k > 0$ depends on $k$.

\paragraph{Rate of convergence of the final estimators: }
\label{para:rate_final}
 We now present the proof that the rate of convergence of the final estimator. The proof for $\hat \alpha$ and $\hat \beta$ are similar and therefore we only the present the proof for $\hat \alpha$. Before delving into the technical details, we introduce some notation: 
\begin{align*}
    \bbM_n(\alpha, d) & = \frac{1}{n}\sum_{i=1}^n \tilde H_k(Y_i - \alpha) \mathds{1}_{X_i^{\top}d \le 0} \\
    \bbR_n(\alpha, d_1, d_2) & = \frac{1}{n}\sum_{i=1}^n \tilde H_k(Y_i - \alpha) \left(\mathds{1}_{X_i^{\top}d_1 \le 0} - \mathds{1}_{X_i^{\top}d_2 \le 0}\right) \\
    \bbM(\alpha, d) & = \bbE[\bbM_n(\alpha, d)] = \bbE\left[\tilde H_k(Y - \alpha) \mathds{1}_{X^{\top}d \le 0}\right] \,.
\end{align*}
From Lemma \ref{lem:huber_lower_bound} we have for all $|\alpha - \alpha_0| \le \eta$ (for some small enough $\eta > 0$): 
$$
\bbM(\alpha, d_0) - \bbM(\alpha_0, d_0) \ge C_k (\alpha - \alpha_0)^2 
$$
In terms of the processes introduced above, we can write our final estimator $\hat \alpha$ as: 
\begin{align}
    \hat \alpha & = \argmin_{\alpha} \frac{1}{n}\sum_{i=1}^n \tilde H_k(Y_i - \alpha) \mathds{1}_{X_i^{\top}\hat d \le 0} \notag \\
    & = \argmin_{\alpha} \bbM_n(\alpha, \hat d) \notag \\
    \label{eq:acc_1} & = \argmin_{\alpha} \left[\bbM_n(\alpha, d_0) + \bbR_n(\alpha, \hat d, d_0)\right]
\end{align}
Consistency of the above estimator follows from the similar calculation as of its previous incarnation, hence we skip it here for brevity.
The remainder term $\bbR_n$ can be bounded as: 
\begin{align*}
\sup_{\alpha: |\alpha - \alpha_0| \le \delta} \left|\bbR_n(\alpha, \hat d, d_0) - \bbR_n(\alpha_0, \hat d, d_0)\right| & \le \frac{k\delta}{n} \sum_{i=1}^n \mathds{1}_{\s(X_i^{\top}\hat d) \neq \s(X_i^{\top}d_0)} \\
& = k\delta \bbP_nf_{\hat d} \\
& = k\delta \left[(\bbP_n - P)f_{\hat d} + Pf_{\hat d}\right]
\end{align*}
Fix $\eps > 0$. Previously, we have established that: 
$$
Pf_{\hat d} = O_p\left(\frac{p}{n}\log{\frac{n}{p}}\right) \,.
$$
Therefore we can find $t_0$ such that for all $t > t_0$: 
$$
\bbP\left(Pf_{\hat d} > t_0\frac{p}{n}\log{\frac{n}{p}}\right) \le \eps \,.
$$
Next we bound the fluctuation of the centered empirical process. \begin{align*}
    & \bbP\left(\frac{n}{p}\left(\log{\frac{n}{p}}\right)^{-1}(\bbP_n - P)f_{\hat d} > t\right) \\
    & \le \bbP\left(\frac{n}{p}\left(\log{\frac{n}{p}}\right)^{-1}(\bbP_n - P)f_{\hat d} > t, Pf_{\hat d} \le t_0 \frac{p}{n}\log{\frac{n}{p}}\right) + \bbP\left(Pf_{\hat d} > t_0 \frac{p}{n}\log{\frac{n}{p}} \right)  \\
    & \le \bbP\left(\sup_{d: Pf_d \le t_0 \frac{p}{n}\log{\frac{n}{p}} }\frac{n}{p}\left(\log{\frac{n}{p}}\right)^{-1}\left\|(\bbP_n - P)f_d\right\| > t\right) + \eps \\
    & \le \frac{n}{p}\left(\log{\frac{n}{p}}\right)^{-1} \times \frac{1}{t}\bbE\left[\sup_{d: Pf_d \le t_0 \frac{p}{n}\log{\frac{n}{p}} }\left\|(\bbP_n - P)f_d\right\|\right] + \eps \\
    & \le \frac{n}{p}\left(\log{\frac{n}{p}}\right)^{-1} \times \frac{\sqrt{t_0}\frac{p}{n}\sqrt{\log{\frac{n}{p}}\log{\frac{n}{t_0p\log{(n/p)}}}} \vee \frac{p}{n}\log{\frac{n}{t_0p\log{(n/p)}}}}{t} + \eps \hspace{0.2in} [\text{Theorem 8.7 of \cite{sen2018gentle}}]  \\
    & = \frac{\sqrt{t_0}}{t}\sqrt{\frac{\log{(n/p)} - \log{\left(t_0\log{(n/p)}\right)}}{\log{(n/p)}}} \ \vee \frac{1}{t}\frac{\log{(n/p)} - \log{\left(t_0\log{(n/p)}\right)}}{\log{(n/p)}} + \eps
\end{align*}
Therefore we have: 
$$
\limsup_{n \to \infty} \bbP\left(\frac{n}{p}\left(\log{\frac{n}{p}}\right)^{-1}(\bbP_n - P)f_{\hat d} > t\right) \le \frac{\sqrt{t_0}}{t} \vee \frac{1}{t} + \eps \,,
$$
which implies that for any fixed $\eps > 0$, we can $t$ large enough to ensure that: 
$$
\limsup_{n \to \infty} \bbP\left(\frac{n}{p}\left(\log{\frac{n}{p}}\right)^{-1}(\bbP_n - P)f_{\hat d} > t\right) \le 2 \eps
$$
Hence we have: 
$$
\sup_{|\alpha - \alpha_0| \le \delta} \left|\bbR_n(\alpha, \hat d, d_0) - \bbR_n(\alpha_0, \hat d, d_0)\right| = \delta \times O_p\left(\frac{p}{n}\log{\frac{n}{p}}\right) \,.
$$
Again fix $\eps > 0$ and choose $t_0 > 0$ such that: 
$$
\limsup_{n \to \infty}\bbP\left(\frac{n}{p}\left(\log{\frac{n}{p}}\right)^{-1} \sup_{|\alpha - \alpha_0| \le \delta} \left|\bbR_n(\alpha, \hat d, d_0) - \bbR_n(\alpha_0, \hat d, d_0)\right| \ge \delta t_0\right) \le \eps \,.
$$
Also note that this $t_0$ only depends on $\eps$, not $\delta$. Henceforth define $r_n = \sqrt{n} \wedge (n/p)(\log{(n/p)})^{-1}$, the desired rate of convergence for the second stage estimator for notational simplicity. Further, by an application of Lemma 2.14.1 of \cite{vdvw96} we have for any $\delta > 0$: 
\begin{align}
    \label{eq:acc_moc}\bbE\left[\sup_{|\alpha - \alpha_0| \le \delta}\left|(\bbM_n - \bbM)(\alpha, d_0) - (\bbM_n - \bbM)(\alpha_0, d_0) \right|\right] \lesssim \frac{\delta}{\sqrt{n}} \,.
\end{align}
Using a shelling type of argument, we have for any $t > 0$: 
\begin{align*}
    & \bbP\left(r_n |\hat \alpha - \alpha_0| > t\right) \\
    & \le \bbP\left(r_n |\hat \alpha - \alpha_0| > t, |\hat \alpha - \alpha_0| \le \eta\right) + \bbP\left(|\hat \alpha - \alpha_0| > \eta\right) \\
    & \le  \bbP\left(r_n |\hat \alpha - \alpha_0| > t, |\hat \alpha - \alpha_0| \le \eta\right) + \eps \\
    & \le \bbP\left(\sup_{\alpha: tr_n^{-1} < |\alpha - \alpha_0| \le \eta}\left\{\bbM_n(\alpha_0, d_0) + \bbR_n(\alpha_0, \hat d, d_0) - \bbM_n(\alpha, d_0) - \bbR_n(\alpha, \hat d, d_0)\right\} \ge 0 \right) + \eps\\
    & \le \bbP\left(\sup_{\alpha: tr_n^{-1} < |\alpha - \alpha_0| \le \eta}\left\{\bbM_n(\alpha_0, d_0) + \bbR_n(\alpha_0, \hat d, d_0) - \bbM_n(\alpha, d_0) - \bbR_n(\alpha, \hat d, d_0)\right\} \ge 0, \right. \\
    & \hspace{20em} \left. \bbP_nf_{\hat d} \le t_0 \frac{p}{n}\log{\frac{n}{p}} \right) + \bbP\left(\bbP_nf_{\hat d} > t_0 \frac{p}{n}\log{\frac{n}{p}}\right) + \eps\\
    & \le \bbP\left(\sup_{\alpha: tr_n^{-1} < |\alpha - \alpha_0| \le \eta}\left\{\bbM_n(\alpha_0, d_0) + \bbR_n(\alpha_0, \hat d, d_0) - \bbM_n(\alpha, d_0) - \bbR_n(\alpha, \hat d, d_0)\right\} \ge 0, \right. \\
    & \hspace{25em} \left. \bbP_nf_{\hat d} \le t_0 \frac{p}{n}\log{\frac{n}{p}} \right) + 2\eps \\
    & \le \sum_{j=1}^{\log_2{(\eta r_n/t)}}\bbP\left(\sup_{\alpha: 2^{j-1}tr_n^{-1} < |\alpha - \alpha_0| \le 2^jtr_n^{-1}}\left\{\bbM_n(\alpha_0, d_0) + \bbR_n(\alpha_0, \hat d, d_0) \right. \right. \\
    & \qquad \qquad \qquad \qquad \left. \left. - \bbM_n(\alpha, d_0) - \bbR_n(\alpha, \hat d, d_0)\right\} \ge 0, \bbP_nf_{\hat d} \le t_0 \frac{p}{n}\log{\frac{n}{p}} \right) + 2\eps \\
    & \le \sum_{j=1}^{\log_2{(\eta r_n/t)}} \bbP\left(\sup_{\alpha: 2^{j-1}tr_n^{-1} < |\alpha - \alpha_0| \le 2^jtr_n^{-1}}\left\{\bbM_n(\alpha_0, d_0) - \bbM_n(\alpha, d_0)\right\} \right. \\
    & \qquad \qquad \qquad \qquad \left. + \sup_{2^{j-1}tr_n^{-1} < |\alpha - \alpha_0| \le 2^jtr_n^{-1}}\left|\bbR_n(\alpha, \hat d, d_0) - \bbR_n(\alpha_0, \hat d, d_0)\right|  \ge 0 , \bbP_nf_{\hat d} \le t_0 \frac{p}{n}\log{\frac{n}{p}}\right) + 2\eps \\
    & \le \sum_{j=1}^{\log_2{(\eta r_n/t)}} \bbP\left(\sup_{\alpha: 2^{j-1}tr_n^{-1} < |\alpha - \alpha_0| \le 2^jtr_n^{-1}}\left\{\bbM_n(\alpha_0, d_0) - \bbM_n(\alpha, d_0)\right\} + 2^jtr_n^{-1} t_0\frac{p}{n}\log{\frac{n}{p}} \ge 0 \right) + 2\eps  \\
    & \le \sum_{j= 1}^{\log_2{(\eta r_n/t)}}\bbP\left(\sup_{\alpha: 2^{j-1}tr_n^{-1} < |\alpha - \alpha_0| \le 2^jtr_n^{-1}}\left\{\bbM_n(\alpha_0, d_0) - \bbM_n(\alpha, d_0)\right\} + 2^jtr_n^{-1} t_0\frac{p}{n}\log{\frac{n}{p}} \ge 0 \right) + 2\eps \\
    & \le   \sum_{j= 1}^{\log_2{(\eta r_n/t)}}\bbP\left(\sup_{\alpha: 2^{j-1}tr_n^{-1} < |\alpha - \alpha_0| \le 2^jtr_n^{-1}}\left\{(\bbM_n - \bbM)(\alpha_0, d_0) - (\bbM_n - \bbM)(\alpha, d_0)\right\} \right. \\
    & \qquad \qquad \qquad \qquad \qquad \left. + 2^jtr_n^{-1} t_0\frac{p}{n}\log{\frac{n}{p}} \ge \inf_{\alpha: 2^{j-1}tr_n^{-1} < |\alpha - \alpha_0| \le 2^jtr_n^{-1}}\left(\bbM(\alpha) - \bbM(\alpha_0)\right) \right) + 2\eps \\
    & \le  \sum_{j= 1}^{\log_2{(\eta r_n/t)}} \frac{\bbE\left[\sup_{2^{j-1}tr_n^{-1} < |\alpha - \alpha_0| \le 2^jtr_n^{-1}}\left\{(\bbM_n - \bbM)(\alpha_0, d_0) - (\bbM_n - \bbM)(\alpha, d_0)\right\}\right] + 2^jtr_n^{-1} t_0\frac{p}{n}\log{\frac{n}{p}}}{\inf_{\alpha: 2^{j-1}tr_n^{-1} < |\alpha - \alpha_0| \le 2^jtr_n^{-1} }\left(\bbM(\alpha) - \bbM(\alpha_0)\right) } + 2\eps \\
    & \le \sum_{j= 1}^{\log_2{(\eta r_n/t)}} \frac{\bbE\left[\sup_{|\alpha - \alpha_0| \le 2^jtr_n^{-1}}\left\{(\bbM_n - \bbM)(\alpha_0, d_0) - (\bbM_n - \bbM)(\alpha, d_0)\right\}\right] + 2^jtr_n^{-1} t_0\frac{p}{n}\log{\frac{n}{p}}}{\inf_{\alpha: |\alpha - \alpha_0| \ge 2^{j-1}tr_n^{-1} }\left(\bbM(\alpha) - \bbM(\alpha_0)\right) } + 2\eps \\
    & \le \sum_{j= 1}^{\log_2{(\eta r_n/t)}} \frac{\bbE\left[\sup_{|\alpha - \alpha_0| \le 2^jtr_n^{-1}}\left\{(\bbM_n - \bbM)(\alpha_0, d_0) - (\bbM_n - \bbM)(\alpha, d_0)\right\}\right] + 2^jt t_0r_n^{-2}}{\inf_{\alpha: |\alpha - \alpha_0| \ge 2^{j-1}tr_n^{-1} }\left(\bbM(\alpha) - \bbM(\alpha_0)\right) } + 2\eps\\
    & \le \sum_{j= 1}^{\log_2{(\eta r_n/t)}} \frac{\bbE\left[\sup_{|\alpha - \alpha_0| \le 2^jtr_n^{-1}}\left\{(\bbM_n - \bbM)(\alpha_0, d_0) - (\bbM_n - \bbM)(\alpha, d_0)\right\}\right] + 2^jt t_0r_n^{-2}}{2^{2(j-1)}t^2r_n^{-2}} + 2\eps\\
    & \le \sum_{j= 1}^{\log_2{(\eta r_n/t)}} \frac{\bbE\left[\sup_{|\alpha - \alpha_0| \le 2^jtr_n^{-1}}\left\{(\bbM_n - \bbM)(\alpha_0, d_0) - (\bbM_n - \bbM)(\alpha, d_0)\right\}\right] }{2^{2(j-1)}t^2r_n^{-2}} + \frac{t_0}
    {t} + 2\eps \\
     & \le \sum_{j= 1}^{\log_2{(\eta r_n/t)}} \frac{2^j t r_n^{-1}}{\sqrt{n}2^{2(j-1)}t^2r_n^{-2}} + \frac{t_0}
    {t} + 2\eps \hspace{0.2in} [\text{From equation \eqref{eq:acc_moc}}]\\
    & \le \frac{1}{t} + \frac{t_0}{t} + 2 \eps \,.
\end{align*}
Taking $t$ large enough we conclude: 
$$
\limsup_{n \to \infty} \bbP\left(r_n |\hat \alpha - \alpha_0| > t\right) \le 3\eps
$$
This completes the proof.

\subsubsection{Case 2: $k = \infty$, i.e. squared error loss}
The proof for $k = \infty$ is similar to that of $0 \le k < \infty$, the only difference is that the collection of functions $\cF$ defined in the proof of the previous part is no longer bounded. Hence we need to modify some parts of the proof carefully to take care of that. 

\paragraph{Consistency: }
\label{para:consistency_l2} 
Consider the same function class $\cF$ as in paragraph \ref{para:consistency_huber}. Note that now any individual function $f_{\theta}$ is: 
\begin{align*}
f_{\theta}(X, \xi) & = \left(\xi(\alpha_0 - \alpha) + \frac12 (\alpha_0 - \alpha)^2 \right)\mathds{1}_{X^{\top}d \vee X^{\top}d_0 \le 0} \\
& \qquad \left(\xi(\alpha_0 - \beta) + \frac12 (\alpha_0 - \beta)^2\right)\mathds{1}_{X^{\top}d_0 \le 0 < X^{\top}d } \\
& \qquad \qquad +  \left(\xi(\beta_0 - \alpha) + \frac12 (\beta_0 - \alpha)^2\right)\mathds{1}_{X^{\top}d \le 0 < X^{\top}d_0} \\
& \qquad \qquad \qquad +   \left(\xi(\beta_0 - \beta) + \frac12(\beta_0 - \beta)^2\right)\mathds{1}_{X^{\top}d \wedge X^{\top}d_0 > 0}
\end{align*}
The envelope function $F$ of $\cF$ is as follows:
\begin{align*}
\sup_{\theta \in \Theta} \left|f_{\theta}(X, \xi)\right| & \le \sup_{(\alpha, \beta) \in \Omega} \left[|\xi|\max\left\{\left|\alpha - \alpha_0\right|, \left|\alpha - \beta_0\right|, \left|\beta - \alpha_0\right|, \left|\beta - \beta_0\right|\right\} \right. \\
& \qquad \qquad \qquad + \left. \frac12 \left(\max\left\{\left|\alpha - \alpha_0\right|, \left|\alpha - \beta_0\right|, \left|\beta - \alpha_0\right|, \left|\beta - \beta_0\right|\right\}\right)^2\right] \\
& \le C|\xi| + \frac{C^2}{2} := F(X, \xi) 
\end{align*} 
The envelope function is integrable and following same analysis as of paragraph \ref{para:consistency_huber} we conclude: 
$$
\sup_{Q} N\left(\eps \left\|F\right\|_{Q, 1}, \cF, L_1(Q)\right) \le K^4\left(\frac{16e}{\eps}\right)^{v + p} \,.
$$
Hence $\cF$ is a Glivenko-Cantelli class of functions and consistency follows from Theorem 2.4.3 of \cite{vdvw96}.

\paragraph{Rate of convergence of the initial estimate: }
To control the modulus of continuity, we can no longer apply Theorem 8.7 of \cite{sen2018gentle} directly here as the functions are not uniformly bounded. Here we use the following modified version of Theorem 1 of \cite{han2019convergence}: 
\begin{proposition}
\label{prop:multiplier_ineq}
Suppose $\{\xi_1, \dots, \xi_n\}$ are independent of random variables of $\{X_1, \dots, X_n\}$ and moreover $\{X_1, \dots, X_n\}$ are permutation invariant. Assume further that there exists a non-decreasing concave function $\varphi_n: \bbR_+ \to \bbR_+$ with $\varphi_n(0) = 0$ and constant $b_n > 0$ such that for $1 \le k \le n$: 
$$
\bbE\left\| \sum_{i=1}^k \eps_i f(X_i)\right\|_\cF \le \varphi_n(k) + b_n
$$ 
for some i.i.d Rademacher random variables $\eps_1, \dots, \eps_n$. Then we have: 
$$
\bbE\left\|\sum_{i=1}^n \xi_i f(X_i) \right\|_\cF \le 4 \int_0^{\infty} \varphi_n\left(\sum_{i=1}^n \bbP\left(\left|\xi_i\right| > t\right) \right) \ dt + 2b_n\bbE\left[\max_{1 \le i \le n} \left|\xi_i\right|\right]\,.
$$
\end{proposition} 
\noindent
The proof of this proposition can be found in the Appendix \ref{app:supp_lemmas}. To apply the above proposition, we define 
\begin{align*}
f_{\theta}(X_i, \xi_i) & = \xi_if_{\theta, 1}(X_i) + f_{\theta, 2}(X_i)
\end{align*}
where: 
\begin{align*}
f_{\theta, 1}(X_i) & = (\alpha_0 - \alpha)\mathds{1}_{X^{\top}d \vee X^{\top}d_0 \le 0}  + (\alpha_0 - \beta)\mathds{1}_{X^{\top}d_0 \le 0 < X^{\top}d } \\
& \qquad \qquad + (\beta_0 - \alpha) \mathds{1}_{X^{\top}d \le 0 < X^{\top}d_0} + (\beta_0 - \beta)\mathds{1}_{X^{\top}d \wedge X^{\top}d_0 > 0} \\
f_{\theta_2}(X_i) & = \frac12(\alpha_0 - \alpha)^2\mathds{1}_{X^{\top}d \vee X^{\top}d_0 \le 0}  + \frac12(\alpha_0 - \beta)^2\mathds{1}_{X^{\top}d_0 \le 0 < X^{\top}d } \\
& \qquad \qquad + \frac12(\beta_0 - \alpha)^2 \mathds{1}_{X^{\top}d \le 0 < X^{\top}d_0} + \frac12 (\beta_0 - \beta)^2\mathds{1}_{X^{\top}d \wedge X^{\top}d_0 > 0} 
\end{align*}
Both the collections $\cF_1 = \{f_{\theta, 1}: d(\theta, \theta_0) \le \delta\}$ and $\cF_2 = \{f_{\theta, 2}: d(\theta, \theta_0) \le \delta\}$ are uniformly bounded with VC dimension of the order $p$. It is also immediate that $Pf_{\theta, j}^2 \lesssim \delta^2$ for all $\theta: d(\theta, \theta_0) \le \delta$, for $j \in \{1, 2\}$. Hence we have from Theorem 8.7 of \cite{sen2018gentle} for any $1 \le k \le n$ and $\eps_1, \dots. \eps_n$ i.i.d Rademacher random variables:  
\begin{align}
\label{eq:non_multiplier} \bbE\left[\sup_{d(\theta, \theta_0) \le \delta}\left|\sum_{i=1}^k \eps_i f_{\theta, j}(X_i)\right|\right] & \le L\left(\delta\sqrt{k}\sqrt{p\log{\frac{AU}{\delta}}} + pU\log{\frac{AU}{\delta}}\right) \\
& := \varphi_n(k) + b_n \notag 
\end{align}
for some constants $L > 0, A> e^2$ and $U$ is the uniform bounds on the individual functions and $\varphi_n(k) = L\delta\sqrt{k}\sqrt{p\log{\left(AU/\delta\right)}}$ and $b_n = pU\log{(AU/\delta)}$. Therefore using Proposition \ref{prop:multiplier_ineq}: 
\begin{align}
\label{eq:moc_multiplier_1}
\bbE\left\| \sum_{i=1}^k \xi_i f(X_i)\right\|_{\cF_1} & \le 4 \int_0^{\infty} \varphi_n\left(n \bbP\left(\left|\xi_1\right| > t\right) \right) \ dt + 2b_n\bbE\left[\max_{1 \le i \le n} \left|\xi_i\right|\right] \notag \\
& \le 4L\delta \sqrt{n} \sqrt{p\log{\left(\frac{AU}{\delta}\right)}}\int_0^{\infty} \sqrt{\bbP(|\xi_1| > t)} \ dt + 2pU\log{\left(\frac{AU}{\delta}\right)}\bbE\left[\max_{1 \le i \le n} \left|\xi_i\right|\right] \notag \\ 
& = 4L\left\|\xi\right\|_{2, 1} \delta \sqrt{n} \sqrt{p\log{\left(\frac{AU}{\delta}\right)}} + 2pU\log{\left(\frac{AU}{\delta}\right)}\bbE\left[\max_{1 \le i \le n} \left|\xi_i\right|\right] 
\end{align}
For the collection $\cF_2$ we can directly use equation \eqref{eq:non_multiplier} for $k = n$ we obtain: 
\begin{equation}
\label{eq:non_multiplier_2}
\bbE\left[\left\|\left(\bbP_n - P\right)f_{\theta, 2}\right\|_{\cF_2}\right] \le L\left(\frac{\delta}{\sqrt{n}}\sqrt{p\log{\frac{AU}{\delta}}} + \frac{pU}{n}\log{\frac{AU}{\delta}}\right)
\end{equation}
Therefore combining equation \eqref{eq:moc_multiplier_1} and \eqref{eq:non_multiplier_2} we conclude: 
\begin{align}
\label{eq:moc_final_1}
\bbE\left[\sup_{d(\theta, \theta_0) \le \delta} \left|\left(\bbP_n - P\right)f_{\theta}\right|\right] & \le L(4\left\|\xi\right\|_{2, 1} + 1)\frac{\delta}{\sqrt{n}}\sqrt{p\log{\frac{AU}{\delta}}} \notag \\
& \qquad \qquad \qquad +  \frac{pU}{n}\log{\frac{AU}{\delta}}\left(1 + 2\bbE\left[\max_{1 \le i \le n} \left|\xi_i\right|\right] \right)
\end{align}
Ignoring constants (as they won't effect the rate of convergence) we can take $\phi_n(\delta)$ is Theorem 3.4.1  of \cite{vdvw96} as: 
$$
\phi_n(\delta) = \delta \sqrt{p\log{\frac{1}{\delta}}} \vee \frac{p}{\sqrt{n}}\log{\frac{1}{\delta}} \bbE\left[\max_{1 \le i \le n} \left|\xi_i\right|\right] \,.
$$
Finally solving the equation the equation $r_n^2\phi_n(1/r_n) \le \sqrt{n}$ we conclude the rate of convergence.

\paragraph{Rate of convergence of the final estimators: }
\label{para:rate_final_l2}
The calculation is exactly same as in Paragraph \ref{para:rate_final} and hence skipped. 
%\newpage
%\vskip 0.2in
%\bibliography{sample}

\subsection{Proof of Theorem \ref{thm:asymptotic_huber}}
To establish the rate of convergence and the asymptotic distribution of the change point estimators obtained via Huber loss, we first need to establish a curvature of the population loss function around its unique minimizer. The following lemma is imperative to that end: 

\begin{lemma}
\label{lem:huber_lower_bound}
If $\xi$ follows a symmetric distribution around the origin with with continuous density $f_{\xi}$ satisfying $f_{\xi}(0) > 0$, then for any $k > 0, |\mu| < 2k$, we have: 
$$
\bbE\left[\tilde H_k(\xi + \mu) - \tilde H_k(\xi)\right] \ge \frac{\mu^2}{2}\bbP\left(-k \le \xi \le k- \mu \right) \ge \frac{\mu^2}{2}\bbP\left(-k \le \xi \le 0 \right) \,.
$$
For $k = 0$, if we choose $\delta$ such that for all $|x| \le \delta$, $f_\xi(x) \ge f_\xi(0)/2$, then we have for $|\mu| \le \delta$: 
$$
\bbE\left[\tilde H_k(\xi + \mu) - \tilde H_k(\xi)\right] \ge \frac{\mu^2}{2}f_{\xi}(0) \,.
$$
\end{lemma}
\noindent
The proof of the above lemma can be found in Appendix \ref{app:supp_lemmas}. Now set $\delta > 0$ such that $\beta_0 + \delta < \alpha_0$. Define the empirical stochastic process $\bbM_n(\theta)$ as: 
\allowdisplaybreaks
\begin{align*}
\bbM_n(\theta) \equiv \bbM_n(\alpha, \beta, d) & =  \frac1n \sum_{i=1}^n \tilde H_k\left(Y_i - \alpha\mathds{1}_{X_i \le d} - \beta \mathds{1}_{X_i > d}\right) \\
& = \frac1n \sum_{i=1}^n \tilde H_k\left(\xi_i +  \alpha_0\mathds{1}_{X_i \le d_0} + \beta_0 \mathds{1}_{X_i > d_0} - \alpha\mathds{1}_{X_i \le d} - \beta \mathds{1}_{X_i > d}\right) \\
& =   \frac1n\sum_{i=1}^n \tilde H_k\left(\xi_i +  \alpha_0 - \alpha\right)\mathds{1}_{X_i \le d_0 \wedge d} +  \frac1n \sum_{i=1}^n \tilde H_k\left(\xi_i +  \alpha_0 - \beta\right)\mathds{1}_{d < X_i \le d_0} \\
& \qquad \qquad +  \frac1n \sum_{i=1}^n \tilde H_k\left(\xi_i +  \beta_0 - \alpha\right)\mathds{1}_{d_0 < X_i \le d} +  \frac1n \sum_{i=1}^n \tilde H_k\left(\xi_i +  \beta_0 - \beta\right)\mathds{1}_{X_i > d \vee d_0}
\end{align*}
This implies the centred empirical stochastic process is: 
\begin{align*}
\bbM_n(\theta) - \bbM_n(\theta_0) & = \bbE\left[\bbM_n(\theta) - \bbM_n(\theta_0)\right] \\
& = \frac1n\sum_{i=1}^n \left(\tilde H_k\left(\xi_i +  \alpha_0 - \alpha\right) - \tilde H_k(\xi_i)\right)\mathds{1}_{X_i \le d_0 \wedge d} \\
& \qquad +  \frac1n \sum_{i=1}^n \left(\tilde H_k\left(\xi_i +  \alpha_0 - \beta\right) - \tilde H_k(\xi_i)\right)\mathds{1}_{d < X_i \le d_0} \\
& \qquad \qquad +  \frac1n \sum_{i=1}^n \left(\tilde H_k\left(\xi_i +  \beta_0 - \alpha\right) -\tilde H_k(\xi_i)\right)\mathds{1}_{d_0 < X_i \le d} \\
& \qquad \qquad \qquad +  \frac1n \sum_{i=1}^n \left(\tilde H_k\left(\xi_i +  \beta_0 - \beta\right) - \tilde H_k(\xi_i)\right)\mathds{1}_{X_i > d \vee d_0}
\end{align*}
and the corresponding population deterministic process: 
\begin{align}
\bbM(\theta) - \bbM(\theta_0) & = \bbE\left[\left(\tilde H_k\left(\xi_i +  \alpha_0 - \alpha\right) - \tilde H_k(\xi_i)\right)\right]\bbP\left(X \le d \wedge d_0\right) \notag \\
& \qquad + \bbE\left[\left(\tilde H_k\left(\xi_i +  \alpha_0 - \beta\right) - \tilde H_k(\xi_i)\right)\right]\bbP\left(d < X < d_0\right) \notag \\
& \qquad \qquad +  \bbE\left[\left(\tilde H_k\left(\xi_i +  \beta_0 - \alpha\right) -\tilde H_k(\xi_i)\right)\right]\bbP\left(d_0 < X < d\right) \notag \\
& \qquad \qquad \qquad +  \bbE\left[\left(\tilde H_k\left(\xi_i +  \beta_0 - \beta\right) - \tilde H_k(\xi_i)\right)\right]\bbP\left(X > d \vee d_0\right) \notag \\
& \ge \frac{C_k}{2}\left[(\alpha_0 - \alpha)^2 \bbP\left(X \le d \wedge d_0\right) + (\beta_0 - \beta)^2\bbP\left(X > d \vee d_0\right) \right. \notag \\
 \label{eq:pop_centered_huber} & \hspace{10em} \left. + |d - d_0|\left\{2(\alpha_0 - \beta_0 - \delta)^2 \right\} \right] \notag \\
 & \ge C_k \left[(\alpha_0 - \alpha)^2 +  (\beta_0 - \beta)^2 + |d - d_0|\right]
\end{align}
for all $|\alpha - \alpha_0| \le \delta, |\beta - \beta_0| \le \delta$, where the penultimate inequality follows from Lemma \ref{lem:huber_lower_bound}. Also note that the definition of $C_k$ is different in the last two lines, but as they are constant, we refrain ourselves from using different notations in each line. 

\paragraph{Consistency: } 
\label{para:consistency}
We have established that $\bbM(\theta)$ has local quadratic curvature with respect to $(\alpha, \beta)$ in a $\delta-$ neighbourhood around the truth. Now to establish the rate of convergence, we first need to establish the consistency of our estimator. To that end, we use Theorem 3.2.2 of \cite{vdvw96}. That the process $\bbM(\theta) - \bbM(\theta_0)$ is continuous with respect to $\theta$ and has a clear minima at $\theta = \theta_0$ is immediate from the definition. Also the tightness of the minimizer $\hat \theta = (\hat \alpha, \hat \beta, \hat d)$ follows directly from our assumption of the compact parameter space $\Theta$. Therefore, all we need to show is that for any compact subset $K \subseteq \Theta$: 
$$
\sup_{\theta \in K} \left|\left(\bbM_n(\theta) - \bbM_n(\theta_0)\right) - \left(\bbM(\theta) - \bbM(\theta_0)\right)\right| = o_p(1) \,.
$$
Towards that direction, define the function $f_{\theta}(X, \xi)$ as: 
\begin{align*}
f_{\theta}(X, \xi) & = \left(\tilde H_k\left(\xi +  \alpha_0 - \alpha\right) -  \tilde H_k(\xi)\right)\mathds{1}_{X \le d_0 \wedge d} \\
& \qquad \left(\tilde H_k\left(\xi +  \alpha_0 - \beta\right) - \tilde H_k(\xi)\right)\mathds{1}_{d < X \le d_0} \\
& \qquad \qquad +  \left(\tilde H_k\left(\xi +  \beta_0 - \alpha\right) -\tilde H_k(\xi)\right)\mathds{1}_{d_0 < X \le d} \\
& \qquad \qquad \qquad +   \left(\tilde H_k\left(\xi +  \beta_0 - \beta\right) - \tilde H_k(\xi)\right)\mathds{1}_{X > d \vee d_0}
\end{align*}
It is immediate from the definition of $f_{\theta}(X, \xi)$ that the collection of functions: 
$$
\cF_K = \left\{f_{\theta}: \theta \in K\right\}
$$
has finite VC dimension. Furthermore, as $K$ is compact, there exist $c$ such that: 
$$
\max_{\theta \in K}\left\{\left|\alpha\right|, \left|\beta\right|, \left|d\right| \right\} \le c \,.
$$
Note that the Huber function $H_k$ is Lipschitz with the Lipschitz constant being $k$. Therefore we have for any $\mu > 0 $: 
\begin{align}
\label{eq:consistency_eq1}  \left|\tilde H_k(\xi + \mu) - \tilde H_k(\xi)\right| &= \frac{k+1}{k}\left|H_k(\xi + \mu) -  H_k(\xi)\right| \le k|\mu| \,.
\end{align}
%\begin{align}
%\left|\tilde H_k(\xi + \mu) - \tilde H_k(\xi)\right| & = \frac{k+1}{k}\left[\left|\frac{(\xi + \mu)^2 - \xi^2}{2}\mathds{1}_{|\xi| \vee |\xi + \mu| \le k}\right| \right. \notag \\
%& \qquad + \left. \left|\left\{k\left(|\xi + \mu| - \frac{k}{2}\right) - \frac{\xi^2}{2}\right\}\mathds{1}_{|\xi| \le k, |\xi + \mu| > k}\right| \right. \notag \\
%& \qquad \qquad + \left. \left|\left\{\frac{(\xi + \mu)^2}{2} - k\left(|\xi| - \frac{k}{2}\right) \right\}\mathds{1}_{|\xi| > k, |\xi + \mu| \le k}\right| \right. \notag \\
%& \qquad \qquad \qquad + \left. \left|\left\{k\left(|\xi + \mu| - |\xi|\right)\right\}\mathds{1}_{|\xi| \wedge |\xi + \mu| > k}\right| \right] \notag \\
%& \le \frac{k+1}{k}\left[2\mu k + \frac{\mu^2}{2} + \left|\left\{k\left(\xi + \mu - \frac{k}{2}\right) - \frac{\xi^2}{2}\right\}\mathds{1}_{(-k \vee k - \mu) \le \xi \le k}\right| \right. \notag \\
%& \qquad + \left. \left|\left\{\frac{(\xi + \mu)^2}{2} + k\left(\xi + \frac{k}{2}\right)\right\}\mathds{1}_{-k - \mu \le \xi \le (-k \wedge k-\mu)}\right|\right] \notag \\
%& \le \frac{k+1}{k}\left[ 2\mu k + \frac{\mu^2}{2} + \left|\left\{\mu k - \frac12 (\xi - k)^2\right\}\mathds{1}_{(-k \vee k - \mu) \le \xi \le k}\right| \right. \notag \\
%\label{eq:consistency_eq1} & \qquad \left. + \left|\left\{\frac12 (\xi + \mu + k)^2 - \mu k\right\}\mathds{1}_{-k - \mu \le \xi \le (-k \wedge k-\mu)}\right|\right] \notag \\ 
%& \le \frac{k+1}{k}\left[4\mu k + \frac12\mu^2 + \left(\mu \wedge 2k\right)^2\right] \,.
%\end{align}
This implies that the function of $\cF$ are uniformly bounded. Hence using Glivenko-Cantelli theorem (e.g. see Theorem 2.8.1 of \cite{vdvw96}) we conclude that: 
\begin{align*}
\sup_{\theta \in K} \left|\left(\bbM_n(\theta) - \bbM_n(\theta_0)\right) - \left(\bbM(\theta) - \bbM(\theta_0)\right)\right| & = \left\|\left(\bbP_n - P\right)\right\|_{\cF} = o_p(1) \,.
\end{align*}
This establishes the consistency.

\paragraph{Tightness upon proper scaling: } 
We next show that: 
$$
\max\left\{\sqrt{n}\left(\hat \alpha - \alpha_0\right), \sqrt{n}\left(\hat \beta - \beta_0\right), n\left(\hat d - d_0\right)\right\} = O_p(1) \,.
$$
Here we apply Theorem 3.2.5 of \cite{vdvw96}. Define a semi-metric on $\Theta$ as: 
$$
d\left(\theta_1, \theta_2\right) = \sqrt{\left(\alpha_1 - \alpha_2\right)^2 + \left(\beta_1 - \beta_2\right)^2 + \left|d_1 - d_2\right|}
$$
From \eqref{eq:pop_centered_huber} we have $\bbM(\theta) - \bbM(\theta_0) \ge C_k d^2(\theta, \theta_0)$. To establish asymptotic equicontinuity of the process we need to bound: 
\allowdisplaybreaks
\begin{align*}
& \bbE\left[\sup_{d(\theta, \theta_0) \le \delta} \left|\bbM_n(\theta) - \bbM_n(\theta_0) -\left(\bbM(\theta) - \bbM(\theta_0)\right)\right| \right] \\
& = \bbE\left[\sup_{d(\theta, \theta_0) \le \delta} \left|\left(\bbP_n - P\right)f_{\theta}\right| \right] \\
& = \bbE\left[\left\|\bbP_n - P\right\|_{\cF_{\delta}}\right]
\end{align*}
where we define the collection $\cF_{\delta}$ as $\cF_{\delta} = \left\{f_{\theta}: d(\theta, \theta_0) \le \delta\right\}$. The envelope function of $\cF_{\delta}$ is defined as:
\begin{align*}
&  \sup_{\theta: d(\theta, \theta_0) \le \delta} \left|\left(\tilde H_k\left(\xi +  \alpha_0 - \alpha\right) - \tilde H_k(\xi)\right)\right|\mathds{1}_{X \le d_0 \wedge d} \\
& \qquad \qquad \qquad \left|\left(\tilde H_k\left(\xi +  \alpha_0 - \beta\right) - \tilde H_k(\xi)\right)\right|\mathds{1}_{d < X \le d_0} \\
& \qquad \qquad \qquad \qquad +  \left|\left(\tilde H_k\left(\xi +  \beta_0 - \alpha\right) -\tilde H_k(\xi)\right)\right|\mathds{1}_{d_0 < X \le d} \\
& \qquad \qquad \qquad \qquad \qquad +   \left|\left(\tilde H_k\left(\xi +  \beta_0 - \beta\right) - \tilde H_k(\xi)\right)\right|\mathds{1}_{X > d \vee d_0} \\
& \le C_k\left(2\delta + \mathds{1}_{d < X \le d_0} + \mathds{1}_{d_0 < X \le d} \right) \\
& \le 2C_k\left(2\delta \vee \mathds{1}_{d < X \le d_0} + \mathds{1}_{d_0 < X \le d} \right) := F_{\delta}(X, \xi) 
\end{align*}
Hence the $L_2(P)$ norm of the envelope function: 
$$
\sqrt{PF_{\delta}^2} \le 2C_k \left(2\delta \vee \sqrt{\bbP\left(d_0 < X < d\right) + \bbP\left(d < X < d_0\right)}\right) \le 4C_k\delta  := \phi_n(\delta)\,.
$$
Hence an application of Theorem 3.2.5 of \cite{vdvw96} yields $\sqrt{n} \ d\left(\hat \theta, \theta_0\right) = O_p(1)$, which completes the proof.

\paragraph{Asymptotic distribution: }
In the final paragraph we establish the asymptotic distribution of $\sqrt{n}\left(\hat \alpha - \alpha_0\right), \sqrt{n}\left(\hat \beta - \beta_0\right)$ and $n\left(\hat d - d_0\right)$. Towards that end, we largely follow the approach of Subsection 14.5.1 of \cite{kosorok2007introduction}. For any $\bh := (h_1, h_2, h_3) \in \bbR^3$ define a paramter vector $\theta_{n, \bh} = \alpha_0 + \frac{h_1}{\sqrt{n}}, \beta_0 + \frac{h_2}{\sqrt{n}}, d_0 + \frac{h_3}{n}$. Define a stochastic process $\bbQ_n$ on $\bbR^3$ as: 
\begin{align*}
\bbQ_n(h_1, h_2, h_3) & = n \times \bbP_n \left(f_{\theta_{n, \bh}} - f_{\theta_0}\right) \\
& =\sum_{i=1}^n \left(\tilde H_k\left(\xi_i +  \frac{h_1}{\sqrt{n}}\right) - \tilde H_k(\xi_i)\right)\mathds{1}_{X_i \le d_0 \wedge d} \\
& \qquad +   \sum_{i=1}^n \left(\tilde H_k\left(\xi_i +  \alpha_0 - \beta_0 - \frac{h_2}{\sqrt{n}}\right) - \tilde H_k(\xi_i)\right)\mathds{1}_{d_0 + \frac{h_3}{n} < X_i \le d_0} \\
& \qquad \qquad +  \sum_{i=1}^n \left(\tilde H_k\left(\xi_i +  \beta_0 - \alpha_0 - \frac{h_1}{\sqrt{n}}\right) -\tilde H_k(\xi_i)\right)\mathds{1}_{d_0 < X_i \le d_0 + \frac{h_3}{n}} \\
& \qquad \qquad \qquad +   \sum_{i=1}^n \left(\tilde H_k\left(\xi_i +\frac{h_2}{\sqrt{n}}\right) - \tilde H_k(\xi_i)\right)\mathds{1}_{X_i > d \vee d_0} \\
& := \bbQ_{n, 1}(\bh) + \bbQ_{n, 2}(\bh) + \bbQ_{n, 3}(\bh) + \bbQ_{n, 4}(\bh) \,.
\end{align*}
It is immediate from the definition of $\bbQ_n(\bh)$ that: 
$$
\widehat{\bh}_n := \left(\sqrt{n}\left(\hat \alpha - \alpha_0\right), \sqrt{n}\left(\hat \beta - \beta_0\right), n\left(\hat d - d_0\right)\right) = \sargmin_{\bh \in \bbR^3}\bbQ_n(\bh) \,.
$$ 
We next show that there exist a stochastic process $\bbQ$ on $\bbR^3$ such that for any compact rectangle $\bbI = I_1 \times I_2 \times I_3 \subset \bbR^3$:
$$
\bbQ_n\vert_{\bbI} \overset{\mathscr{L}}{\implies} \bbQ\vert_{\bbI} \,.
$$
where the process $\bbQ$ is defined as: 
\begin{align*}
\bbQ(\bh) & = \left(h_1\sigma_k \sqrt{F_X(d_0)} \times Z_1 + \frac{h_1^2}{2} \mu_kF_X(d_0)\right) \\
& \hspace{5em} + \left(h_2 \sigma_k \sqrt{\bar F_X(d_0)} \times Z_2 +   \frac{h_2^2}{2}\mu_k  \bar F_X(d_0)\right) \\
&  \hspace{10em} + \CPP\left(\tilde H_k\left(\xi_i + (\alpha_0 - \beta_0)\right) - \tilde H_k(\xi_i), f_X(\theta_0)\right) \,.
\end{align*}
with $Z_1, Z_2 \overset{i.i.d}{\sim} \cN(0, 1)$ and CPP is (as defined in the main paper) compound Poisson process. Note that the stochastic process $\bbQ_n$ is continuous with respect to its first two co-ordinates and cadlag (right continuous with left limit) with respect to its third co-ordinate. Hence to establish the convergence of $\left\{\bbQ_n \vert_{\bbI}\right\}_{n \in \bbN}$ we need to use Skorohod topology. We mainly use Theorem 13.5 of \cite{billingsley2013convergence} to establish the convergence result. Towards that end, define: 
\begin{align*}
\tilde \xi_{i ,h_1}  & =\left(\tilde H_k\left(\xi_i +  \frac{h_1}{\sqrt{n}}\right) - \tilde H_k(\xi_i)\right)- \bbE\left[\left(\tilde H_k\left(\xi_i +  \frac{h_1}{\sqrt{n}}\right) - \tilde H_k(\xi_i)\right)\right] \\
\tilde \xi_{i ,h_2}  & =  \left(\tilde H_k\left(\xi_i +\frac{h_2}{\sqrt{n}}\right) - \tilde H_k(\xi_i)\right) - \bbE\left[ \left(\tilde H_k\left(\xi_i +\frac{h_2}{\sqrt{n}}\right) - \tilde H_k(\xi_i)\right)\right]
\end{align*}
and another stochastic process $\tilde \bbQ_n(\bh)$ as: 
\begin{align*}
\tilde \bbQ_n(\bh) & = \sum_{i=1}^n \tilde \xi_{i, h_1} \mathds{1}_{X_i \le d_0 \wedge d_0 + \frac{h_3}{n}} \\
& \qquad +   \sum_{i=1}^n \left(\tilde H_k\left(\xi_i +  (\alpha_0 - \beta_0)\right) - \tilde H_k(\xi_i)\right)\mathds{1}_{d_0 + \frac{h_3}{n} < X_i \le d_0}\\
& \qquad \qquad +   \sum_{i=1}^n \left(\tilde H_k\left(\xi_i +  (\beta_0 - \alpha_0 )\right) - \tilde H_k(\xi_i)\right)\mathds{1}_{d_0  < X_i \le d_0 + \frac{h_3}{n}} \\
& \qquad \qquad \qquad +  \sum_{i=1}^n \tilde \xi_{i, h_2} \mathds{1}_{X_i > d_0 \vee d_0 + \frac{h_3}{n}} \\
& := \tilde \bbQ_{n, +}(\bh)\mathds{1}_{h_3 \ge 0} + \tilde \bbQ_{n, -}(\bh)\mathds{1}_{h_3 < 0} \,.
\end{align*}
Hence we have the following decomposition: 
\begin{equation}
\label{eq:decomposition}
 \bbQ_n(\bh) = \tilde \bbQ_n(\bh) + \mathfrak{E}_n(\bh) + \mathfrak{R}_n(\bh)
\end{equation}
where: 
\begin{align*}
\mathfrak{R}_n(\bh) & =  \sum_{i=1}^n \left(\tilde H_k\left(\eps_i +  (\alpha_0 - \beta_0) - \frac{h_2}{\sqrt{n}}\right)  - \tilde H_k\left(\eps_i +  (\alpha_0 - \beta_0)\right)\right)\mathds{1}_{d_0 + \frac{h_3}{n} < X_i \le d_0} \\
& \qquad \qquad +    \sum_{i=1}^n \left(\tilde H_k\left(\eps_i +  (\beta_0 - \alpha_0) - \frac{h_1}{\sqrt{n}}\right)  - \tilde H_k\left(\eps_i +  (\beta_0 - \alpha_0)\right)\right)\mathds{1}_{d_0 < X_i \le d_0 + \frac{h_3}{n}} 
\end{align*}
and 
\begin{align*}
\mathfrak{E}_n(\bh) & = \bbE\left[\left(\tilde H_k\left(\xi_i + \frac{h_1}{\sqrt{n}}\right) - \tilde H_k\left(\xi_i\right)\right)\right]\sum_{i=1}^n\mathds{1}_{X_i \le d_0 \wedge d_0 + \frac{h_3}{n}} \\
& \qquad \qquad +  \bbE\left[\left(\tilde H_k\left(\xi_i + \frac{h_2}{\sqrt{n}}\right) - \tilde H_k\left(\xi_i\right)\right)\right]\sum_{i=1}^n\mathds{1}_{X_i > d_0 \vee d_0 + \frac{h_3}{n}} 
\end{align*}
We next show $\mathfrak{R}_n(\bh)$ is $o_p(1)$ uniformly over a compact set: 
\begin{align*}
\left|\mathfrak{R}_n(\bh) \right| & = \left|\sum_{i=1}^n \left(\tilde H_k\left(\xi_i +  (\alpha_0 - \beta_0) - \frac{h_2}{\sqrt{n}}\right) - \tilde H_k\left(\xi_i +  (\alpha_0 - \beta_0)\right)\right)\mathds{1}_{d_0 + \frac{h_3}{n} < X_i \le d_0} \right. \\
& \qquad \qquad + \left.  \sum_{i=1}^n \left(\tilde H_k\left(\xi_i +  ( \beta_0 - \alpha_0) - \frac{h_1}{\sqrt{n}}\right)  - \tilde H_k\left(\xi_i +  (\beta_0 - \alpha_0)\right)\right)\mathds{1}_{d_0 < X_i \le d_0 + \frac{h_3}{n}} \right| \\
& \le \sum_{i=1}^n \left|\tilde H_k\left(\xi_i +  (\alpha_0 - \beta_0) - \frac{h_2}{\sqrt{n}}\right) - \tilde H_k\left(\xi_i +  (\alpha_0 - \beta_0)\right)\right|\mathds{1}_{d_0 + \frac{h_3}{n} < X_i \le d_0} \\
& \qquad \qquad +    \sum_{i=1}^n \left|\tilde H_k\left(\xi_i +  ( \beta_0 - \alpha_0) - \frac{h_1}{\sqrt{n}}\right)  - \tilde H_k\left(\xi_i +  (\beta_0 - \alpha_0)\right)\right|\mathds{1}_{d_0 < X_i \le d_0 + \frac{h_3}{n}} \\
& \le \frac{(k+1)h_2}{\sqrt{n}} \sum_{i=1}^n\mathds{1}_{d_0 + \frac{h_3}{n} < X_i \le d_0}  + \frac{(k+1)h_1}{\sqrt{n}}  \sum_{i=1}^n \mathds{1}_{d_0 < X_i \le d_0 + \frac{h_3}{n}} 
\end{align*}
Now suppose $\bh \in \bbI$. There $|h_1| \vee |h_2| \vee |h_3| \le K$ for some $K > 0$. Hence we have: 
\begin{align*}
\sup_{\bh \in \bbI} \left|\mathfrak{R}_n(\bh) \right| & \le \frac{(k+1)K}{\sqrt{n}} \left[\sum_{i=1}^n\mathds{1}_{d_0 - \frac{K}{n} < X_i \le d_0} + \sum_{i=1}^n \mathds{1}_{d_0 < X_i \le d_0 + \frac{K}{n}}\right] \\
& =  \frac{(k+1)K}{\sqrt{n}}\sum_{i=1}^n\mathds{1}_{d_0 - \frac{K}{n} < X_i \le d_0 + \frac{K}{n}}
\end{align*}
as we know: 
$$
\sum_{i=1}^n\mathds{1}_{d_0 - \frac{K}{n} < X_i \le d_0 + \frac{K}{n}} \overset{\mathscr{L}}{\implies} \text{Pois}\left(Kf_X(d_0)\right)
$$
we conclude that: 
\begin{equation}
\label{eq:remainder_conv}
\sup_{\bh \in \bbI} \left|\mathfrak{R}_n(\bh) \right|  = O_p(n^{-1/2}) = o_p(1) \,.
\end{equation}
We now establish the convergence of $\mathfrak{E}(\bh)$: 
\begin{align*}
\mathfrak{E}_n(\bh) & = \bbE\left[\left(\tilde H_k\left(\xi_i + \frac{h_1}{\sqrt{n}}\right) - \tilde H_k(\xi_i)\right)\right]\sum_{i=1}^n\mathds{1}_{X_i \le d_0 \wedge d_0 + \frac{h_3}{n}} \\
& \qquad \qquad +  \bbE\left[\left(\tilde H_k\left(\xi_i + \frac{h_2}{\sqrt{n}}\right) - \tilde H_k(\xi_i)\right)\right]\sum_{i=1}^n\mathds{1}_{X_i > d_0 \vee d_0 + \frac{h_3}{n}}  \\
& =  \frac{k+1}{k}\left\{\frac{h_1^2}{2n}\bbP\left( -k \le \xi \le k -  \frac{h_1}{\sqrt{n}}\right)  + \bbE\left[\left( \frac{h_1}{\sqrt{n}} k - \frac{h_1}{\sqrt{n}} \xi - \frac12 \left(\xi - k\right)^2\right) \mathds{1}_{k -  \frac{h_1}{\sqrt{n}} \le \xi \le k}\right]  \right. \\
& \qquad \qquad \qquad \qquad + \left. \bbE\left[\frac12 \left(\xi + \frac{h_1}{\sqrt{n}} + k\right)^2 \mathds{1}_{-k -  \frac{h_1}{\sqrt{n}} \le \xi \le -k}\right] \right\}\sum_{i=1}^n\mathds{1}_{X_i \le d_0 \wedge d_0 + \frac{h_3}{n}} \\
& \qquad \qquad +  \frac{k+1}{k}\left\{\frac{h_2^2}{2n}\bbP\left( -k \le \xi \le k -  \frac{h_2}{\sqrt{n}}\right)  + \bbE\left[\left( \frac{h_2}{\sqrt{n}} k -  \frac{h_2}{\sqrt{n}} \xi - \frac12 \left(\xi - k\right)^2\right) \mathds{1}_{k -  \frac{h_2}{\sqrt{n}} \le \xi \le k}\right] \right.   \\
& \qquad \qquad \qquad \qquad + \left. \bbE\left[\frac12 \left(\xi +  \frac{h_2}{\sqrt{n}} + k\right)^2 \mathds{1}_{-k -  \frac{h_2}{\sqrt{n}} \le \xi \le -k}\right] \right\}\sum_{i=1}^n\mathds{1}_{X_i > d_0 \vee d_0 + \frac{h_3}{n}} \\
& = n \frac{k+1}{k}\left\{\frac{h_1^2}{2n}\bbP\left( -k \le \xi \le k -  \frac{h_1}{\sqrt{n}}\right)  + \bbE\left[\left( \frac{h_1}{\sqrt{n}} k - \frac{h_1}{\sqrt{n}} \xi - \frac12 \left(\xi - k\right)^2\right) \mathds{1}_{k -  \frac{h_1}{\sqrt{n}} \le \xi \le k}\right]  \right. \\
& \qquad \qquad \qquad \qquad + \left. \bbE\left[\frac12 \left(\xi + \frac{h_1}{\sqrt{n}} + k\right)^2 \mathds{1}_{-k -  \frac{h_1}{\sqrt{n}} \le \xi \le -k}\right] \right\}\frac1n \sum_{i=1}^n\mathds{1}_{X_i \le d_0 \wedge d_0 + \frac{h_3}{n}} \\
& \qquad \qquad + n \frac{k+1}{k}\left\{\frac{h_2^2}{2n}\bbP\left( -k \le \xi \le k -  \frac{h_2}{\sqrt{n}}\right)  + \bbE\left[\left( \frac{h_2}{\sqrt{n}} k -  \frac{h_2}{\sqrt{n}} \xi - \frac12 \left(\xi - k\right)^2\right) \mathds{1}_{k -  \frac{h_2}{\sqrt{n}} \le \xi \le k}\right] \right.   \\
& \qquad \qquad \qquad \qquad + \left. \bbE\left[\frac12 \left(\xi +  \frac{h_2}{\sqrt{n}} + k\right)^2 \mathds{1}_{-k -  \frac{h_2}{\sqrt{n}} \le \xi \le -k}\right] \right\}\frac1n \sum_{i=1}^n\mathds{1}_{X_i > d_0 \vee d_0 + \frac{h_3}{n}} \\
& =  \frac{k+1}{k}\left\{\frac{h_1^2}{2}\bbP\left( -k \le \xi \le k -  \frac{h_1}{\sqrt{n}}\right)  + n \bbE\left[\left( \frac{h_1}{\sqrt{n}} k - \frac{h_1}{\sqrt{n}} \xi - \frac12 \left(\xi - k\right)^2\right) \mathds{1}_{k -  \frac{h_1}{\sqrt{n}} \le \xi \le k}\right]  \right. \\
& \qquad \qquad \qquad \qquad + \left. n \bbE\left[\frac12 \left(\xi + \frac{h_1}{\sqrt{n}} + k\right)^2 \mathds{1}_{-k -  \frac{h_1}{\sqrt{n}} \le \xi \le -k}\right] \right\}\frac1n \sum_{i=1}^n\mathds{1}_{X_i \le d_0 \wedge d_0 + \frac{h_3}{n}} \\
& \qquad \qquad +  \frac{k+1}{k}\left\{\frac{h_2^2}{2}\bbP\left( -k \le \xi \le k -  \frac{h_2}{\sqrt{n}}\right)  + n \bbE\left[\left( \frac{h_2}{\sqrt{n}} k -  \frac{h_2}{\sqrt{n}} \xi - \frac12 \left(\xi - k\right)^2\right) \mathds{1}_{k -  \frac{h_2}{\sqrt{n}} \le \xi \le k}\right] \right.   \\
& \qquad \qquad \qquad \qquad + \left. n \bbE\left[\frac12 \left(\xi +  \frac{h_2}{\sqrt{n}} + k\right)^2 \mathds{1}_{-k -  \frac{h_2}{\sqrt{n}} \le \xi \le -k}\right] \right\}\frac1n \sum_{i=1}^n\mathds{1}_{X_i > d_0 \vee d_0 + \frac{h_3}{n}}\end{align*}
From strong law of large numbers we have: 
\begin{align}
\label{eq:en_conv_1} \frac1n \sum_{i=1}^n\mathds{1}_{X_i \le d_0 \wedge d_0 + \frac{h_3}{n}} & \overset{p}{\longrightarrow} \bbP(X \le d_0) \,, \\
\label{eq:en_conv_2}  \frac1n \sum_{i=1}^n\mathds{1}_{X_i > d_0 \vee d_0 + \frac{h_3}{n}} & \overset{p}{\longrightarrow} \bbP(X > d_0) \,, 
\end{align}
For the other terms in the expectation: 
\begin{align}
\label{eq:en_conv_3}  \frac{h_1^2}{2} \bbP\left( -k \le \xi \le k -  \frac{h_1}{\sqrt{n}}\right) & \overset{}{\longrightarrow} \frac{h_1^2}{2} \bbP\left( -k \le \eps \le k \right) \,, \\
\label{eq:en_conv_4}  \frac{h_2^2}{2} \bbP\left( -k \le \xi \le k -  \frac{h_2}{\sqrt{n}}\right) & \overset{}{\longrightarrow} \frac{h_2^2}{2} \bbP\left( -k \le \eps \le k \right)  \,.
\end{align}
For the other terms: 
\begin{align}
& n\bbE\left[\left( \frac{h_2}{\sqrt{n}} k -  \frac{h_2}{\sqrt{n}} \xi - \frac12 \left(\xi - k\right)^2\right) \mathds{1}_{k -  \frac{h_2}{\sqrt{n}} \le \xi \le k}\right] \notag \\
& = n\int_{k -  \frac{h_2}{\sqrt{n}}}^k \left( \frac{h_2}{\sqrt{n}} k -  \frac{h_2}{\sqrt{n}} x - \frac12 \left(x - k\right)^2\right)f_{\xi}(x) \ dx \notag \\
& = n\left[\frac{h_2}{\sqrt{n}}\int_{k -  \frac{h_2}{\sqrt{n}}}^k (k -  x) f_{\xi}(x) \ dx - \frac12 \int_{k -  \frac{h_2}{\sqrt{n}}}^k \left(x - k\right)^2 \ f_{\xi}(x)\ dx\right] \notag \\
& = n \left[\frac{h_2}{\sqrt{n}}\int_{0}^{\frac{h_2}{\sqrt{n}}} z f_{\xi}(z-k) \ dx - \frac12 \int_{-  \frac{h_2}{\sqrt{n}}}^0 z^2 \ f_{\xi}(z+k)\ dx\right] \notag \\
& \le Cn \left[\frac{h_2}{\sqrt{n}}\int_{0}^{\frac{h_2}{\sqrt{n}}} z \ dz +\frac12 \int_{-  \frac{h_2}{\sqrt{n}}}^0 z^2 \ dz \right] \hspace{0.2in} [C = \max_x f_{\xi}(x)] \notag \\ 
\label{eq:en_conv_5}  & = n \times O(n^{-3/2}) = o(1)\,.
\end{align}
and 
\begin{align}
& n\bbE\left[\frac12 \left(\xi +  \frac{h_2}{\sqrt{n}} + k\right)^2 \mathds{1}_{-k -  \frac{h_2}{\sqrt{n}} \le \xi \le -k}\right] \notag \\
& = n\int_{-k -  \frac{h_2}{\sqrt{n}}}^{-k} \frac12 \left(x +  \frac{h_2}{\sqrt{n}} + k\right)^2 \ f_{\xi}(x) \ dx \notag \\
\label{eq:en_conv_6}  &  = \frac{h_2^3}{6n^{1/2}}f_{\xi}(k) + o(1) = o(1) \,.
\end{align}
Similar calculation holds for the terms involving $h_1$. Hence we conclude combining equations \eqref{eq:en_conv_1} - \eqref{eq:en_conv_6} we conclude: 
\begin{align}
\label{eq:mean_conv}
\mathfrak{E}_n(\bh) & \overset{P}{\longrightarrow}   \frac{k+1}{k}\left[\frac{h_1^2}{2} \bbP\left( -k \le \eps \le k \right)F_X(d_0) + \frac{h_2^2}{2} \bbP\left( -k \le \eps \le k \right)\bar F_X(d_0)\right]\,.
\end{align}
Finally we show the weak convergence of $\tilde \bbQ_n(\bh)$ to $\bbQ(\bh)$, for which we need to show that the collection $\{\tilde \bbQ_n(\bh)\}_{n \in \bbN}$ is tight with respect to appropriate topology and every finite dimensional projection of $\tilde \bbQ_n(\bh)$ converges to that of $ \bbQ(\bh)$. We embark on by showing that for any fixed $\bh$, $\tilde \bbQ_n(\bh)$ converges to $ \bbQ(\bh)$ in distribution. 
Towards that direction, fix $h_3 > 0$:
\begin{align*}
\bbE\left[\exp{\left(it\tilde \bbQ_n(\bh)\right)}\right] & = \bbE\left[\bbE\left[\exp{\left(it\tilde \bbQ_n(\bh)\right)} \mid X_1, \dots, X_n \right]\right] \\
& := \bbE\left[\bbE_{\bX}\left[\exp{\left(it\tilde \bbQ_n(\bh)\right)}\right]\right]
\end{align*} 
We start with analyzing the inner expectation $\bbE_{\bX}\left[\exp{\left(it\tilde \bbQ_n(\bh)\right)}\right]$: 
\begin{align}
\bbE_{\bX}\left[\exp{\left(it\tilde \bbQ_n(\bh)\right)}\right] & = \bbE_{\bX}\left[\exp\left\{\left(it\sum_{i=1}^n \tilde \xi_{i, h_1} \mathds{1}_{X_i \le d_0} \right. \right. \right. \notag \\
& \qquad \qquad  \qquad   \qquad  \left. \left. \left. +   \sum_{i=1}^n \left(\tilde H_k\left(\xi_i + (\alpha_0 - \beta_0)\right) - \tilde H_k(\xi_i)\right)\mathds{1}_{d_0 < X_i \le d_0 + \frac{h_3}{n}} \right. \right. \right. \notag \\
& \hspace{20em} \left. \left. \left. +  \sum_{i=1}^n \tilde \xi_{i, h_2} \mathds{1}_{X_i > d_0 + \frac{h_3}{n}}\right)\right\}\right] \notag \\
& = \left(\phi_{\tilde \xi_{h_1}}(t)\right)^{\sum_{i=1}^n \mathds{1}_{X_i \le d_0}} \times \left(\phi_{\tilde \xi_{h_2}}(t)\right)^{\sum_{i=1}^n \mathds{1}_{X_i > d_0 + \frac{h_3}{n}}} \notag \\
& \qquad \qquad \qquad \qquad \times \left(\phi_{\tilde H_k\left(\xi_i + (\alpha_0 - \beta_0)\right) - \tilde H_k(\xi_i)}(t)\right)^{\sum_{i=1}^n \mathds{1}_{d_0 < X_i \le d_0 + \frac{h_3}{n}} } \notag \\
& = \left(\phi_{\tilde \xi_{h_1}}(t)\right)^{n \times \frac{1}{n}\sum_{i=1}^n \mathds{1}_{X_i \le d_0}} \times \left(\phi_{\tilde \xi_{h_2}}(t)\right)^{n \times \frac{1}{n} \sum_{i=1}^n \mathds{1}_{X_i > d_0 + \frac{h_3}{n}}} \notag \\
\label{eq:cf_1} & \qquad \qquad \qquad \qquad \times \left(\phi_{\tilde H_k\left(\xi_i + (\alpha_0 - \beta_0)\right) - \tilde H_k(\xi_i)}(t)\right)^{\sum_{i=1}^n \mathds{1}_{d_0 < X_i \le d_0 + \frac{h_3}{n}} }
\end{align}
To show the convergence of the characteristic function of $\tilde \xi_{h_1}$ (and similarly for $\tilde \xi_{h_2}$) we first note that the variance of $\tilde \xi_{h_1}$ for $h_1 > 0$ is: 
\begin{align*}
\var\left(\tilde \xi\right) & =\var\left(\tilde H_k\left(\xi_i +  \frac{h_1}{\sqrt{n}}\right) - \tilde H_k(\xi_i)\right) \\
& = \bbE\left[\left(\tilde H_k\left(\xi_i +  \frac{h_1}{\sqrt{n}}\right) - \tilde H_k(\xi_i)\right)^2\right] + O(n^{-2})  \hspace{0.2in} [O(n^{-2}) \text{ follows from the analysis of } \mathfrak{E}_n(\bh)]\\
& = \left( \frac{k+1}{k}\right)^2 \left\{\bbE\left[\left(\frac{h_1}{\sqrt{n}}\xi + \frac{h_1^2}{2n}\right)^2\mathds{1}_{-k \le \xi \le k - \frac{h_1}{\sqrt{n}}}\right] \right. \\
& \qquad + \left. \bbE\left[\left(\frac12 \left(\xi + \frac{h_1}{\sqrt{n}} + k\right)^2 - \frac{h_1}{\sqrt{n}}k\right)^2\mathds{1}_{-k - \frac{h_1}{\sqrt{n}} \le \xi \le -k}\right] \right. \\
& \qquad \qquad + \left. \bbE\left[\left(-\frac12 \left(\xi - k\right)^2 + k\frac{h_1}{\sqrt{n}}\right)^2\mathds{1}_{k-\frac{h_1}{\sqrt{n}} \le \xi \le k}\right] \right. \\
& \qquad \qquad \qquad + \left. \frac{k^2h_1^2}{n}\left\{\bbP\left(\xi > k\right) + \bbP\left(\xi < -k - \frac{h_1}{\sqrt{n}}\right)\right\} + O(n^{-2}) \right\}\\
& =\left( \frac{k+1}{k}\right)^2\left( \frac{h_1^2}{n}\bbE\left[\xi^2\mathds{1}_{-k \le \xi \le k}\right] + \frac{k^2h_1^2}{n}\left\{\bbP\left(\xi > k\right) + \bbP\left(\xi < -k\right)\right\} \right) + o(n^{-1}) \\
& = \left( \frac{k+1}{k}\right)^2\left(\frac{h_1^2}{n}\bbE\left[\xi^2\mathds{1}_{-k \le \xi \le k}\right] + \frac{2k^2h_1^2}{n}\bbP\left(\xi > k\right)\right) + o(n^{-1}) \\
& := \frac{h_1^2 \sigma_k^2}{n} + o(n^{-1}) \,.
\end{align*}
where the variance parameter $\sigma_k^2$ is defined as: 
$$
 \left( \frac{k+1}{k}\right)^2\left(\bbE\left[\xi^2\mathds{1}_{-k \le \xi \le k}\right] + 2k^2\bbP\left(\xi > k\right)\right) \,.
$$
Similar calculation holds for $h_1 < 0$ and for $h_2$. Hence going back to equation \eqref{eq:cf_1} we have: 
\begin{align*}
\bbE_{\bX}\left[\exp{\left(it\tilde \bbQ_n(\bh)\right)}\right] & = \left(1 - \frac{t^2}{2}\frac{h_1^2\sigma_k^2}{n} + o(n^{-1})\right)^{n \times \frac{1}{n}\sum_{i=1}^n \mathds{1}_{X_i \le d_0}} \\
& \qquad \qquad \times\left(\phi_{\tilde H_k\left(\xi_i + (\alpha_0 - \beta_0)\right) - \tilde H_k(\xi_i)}(t)\right)^{\sum_{i=1}^n \mathds{1}_{d_0 < X_i \le d_0 + \frac{h_3}{n}} } \notag \\
\label{eq:cf_1} & \qquad \qquad \qquad \qquad \times \left(1 - \frac{t^2}{2}\frac{h_2^2\sigma_k^2}{n} + o(n^{-1})\right)^{n \times \frac{1}{n}\sum_{i=1}^n \mathds{1}_{X_i >  d_0 + \frac{h_3}{n}} }
\end{align*}
As
\begin{align*}
\frac{1}{n}\sum_{i=1}^n \mathds{1}_{X_i \le d_0} & \overset{a.s.}{\longrightarrow} F_X(d_0) \\
\frac{1}{n}\sum_{i=1}^n \mathds{1}_{X_i >  d_0 + \frac{h_3}{n}} & \overset{a.s.}{\longrightarrow} \bar F_X(d_0)
\end{align*}
we conclude: 
\begin{align*}
 \left(1 - \frac{t^2}{2}\frac{h_1^2\sigma_k^2}{n} + o(n^{-1})\right)^{n \times \frac{1}{n}\sum_{i=1}^n \mathds{1}_{X_i \le d_0}} & \overset{a.s.}{\longrightarrow} e^{-\frac{t^2h_1^2\sigma_k^2}{2}F_X(d_0)} \\
  \left(1 - \frac{t^2}{2}\frac{h_2^2\sigma_k^2}{n} + o(n^{-1})\right)^{n \times \frac{1}{n}\sum_{i=1}^n \mathds{1}_{X_i > d_0 + \frac{h_3}{n}} } & \overset{a.s.}{\longrightarrow} e^{-\frac{t^2h_2^2\sigma_k^2}{2}\bar F_X(d_0)}
\end{align*}
Also we know: 
$$
\sum_{i=1}^n \mathds{1}_{d_0 < X_i \le d_0 + \frac{h_3}{n}} \overset{\mathscr{L}}{\implies} \text{Pois}\left(f_X(d_0)h_3\right) 
$$
which further implies: 
$$
\left(\phi_{\tilde H_k\left(\xi_i + (\alpha_0 - \beta_0)\right) - \tilde H_k(\xi_i)}(t)\right)^{\sum_{i=1}^n \mathds{1}_{d_0 < X_i \le d_0 + \frac{h_3}{n}} } \overset{\mathscr{L}}{\implies} \left(\phi_{\tilde H_k\left(\xi_i + (\alpha_0 - \beta_0)\right) - \tilde H_k(\xi_i)}(t)\right)^{\text{Pois}\left(f_X(\theta_0)h_3\right) } 
$$
Hence combining these we conclude: 
\begin{align*}
\bbE_{\bX}\left[\exp{\left(it\tilde \bbQ_n(\bh)\right)}\right] & = \left(\phi_{\tilde \xi_{h_1}}(t)\right)^{n \times \frac{1}{n}\sum_{i=1}^n \mathds{1}_{X_i \le d_0}} \times \left(\phi_{\tilde \xi_{h_2}}(t)\right)^{n \times \frac{1}{n} \sum_{i=1}^n \mathds{1}_{X_i > d_0 + \frac{h_3}{n}}} \notag \\
 & \qquad \qquad \qquad \qquad \times \left(\phi_{\tilde H_k\left(\xi_i + (\alpha_0 - \beta_0)\right) - \tilde H_k(\xi_i)}(t)\right)^{\sum_{i=1}^n \mathds{1}_{d_0 < X_i \le d_0 + \frac{h_3}{n}} }\\
& \overset{\mathscr{L}}{\implies}  e^{-\frac{t^2h_1^2\sigma_k^2}{2}F_X(d_0)} \times  e^{-\frac{t^2h_2^2\sigma_k^2}{2}\bar F_X(d_0)} \\
& \qquad \qquad \qquad \qquad \times \left(\phi_{\tilde H_k\left(\xi_i + (\alpha_0 - \beta_0)\right) - \tilde H_k(\xi_i)}(t)\right)^{\text{Pois}\left(f_X(d_0)h_3\right) } 
\end{align*}
Applying DCT and taking expectation on the both side we conclude: 
\begin{align*}
\bbE\left[\exp{\left(it\tilde \bbQ_n(\bh)\right)}\right] & \longrightarrow e^{-\frac{t^2h_1^2\sigma_k^2}{2}F_X(d_0)} \times  e^{-\frac{t^2h_2^2\sigma_k^2}{2}\bar F_X(d_0)} \\
& \qquad \qquad \qquad \qquad \times \bbE\left[\left(\phi_{\tilde H_k\left(\xi_i + (\alpha_0 - \beta_0)\right) - \tilde H_k(\xi_i)}(t)\right)^{\text{Pois}\left(f_X(d_0)h_3\right) } \right] \,.
\end{align*}
This concludes that $\tilde \bbQ_n(\bh) \overset{\mathscr{L}}{\implies}  \bbQ(\bh)$. The proof of the fact that for any finite collection $(\bh_1, \dots, \bh_l)$: 
\begin{equation}
\label{eq:finite_dim_tilde}
\left(\tilde \bbQ_n(\bh_1), \dots, \tilde \bbQ_n(\bh_l)\right)  \overset{\mathscr{L}}{\implies}  \left( \bbQ(\bh_1), \dots,   \bbQ(\bh_l)\right)
\end{equation}
is similar (same analysis of characteristic function) and hence skipped for brevity. Interested readers can take a look at the proof of Lemma 3.2 of \cite{lan2009change} or the proof of Theorem 5 of \cite{kosorok2007inference} for more details of this type of calculations. We next establish the tightness of the process. Define another process $\dbtilde \bbQ_n(\bh)$ as: 
\allowdisplaybreaks
\begin{align}
\dbtilde \bbQ_n(\bh) & = \sum_{i=1}^n \tilde \xi_{i, h_1} \mathds{1}_{X_i \le d_0} \notag \\
& \qquad +  \sum_{i=1}^n \left(H_k\left(\xi_i + (\alpha_0 - \beta_0)\right) - \tilde H_k(\xi_i)\right)\mathds{1}_{d_0 + \frac{h_3}{n} < X_i \le d_0} \notag \\
& \qquad \qquad +   \sum_{i=1}^n \left(H_k\left(\xi_i + (\alpha_0 - \beta_0)\right) - \tilde H_k(\xi_i)\right)\mathds{1}_{d_0 < X_i \le d_0 + \frac{h_3}{n}} \notag \\
& \qquad \qquad \qquad +  \sum_{i=1}^n \tilde \xi_{i, h_2} \mathds{1}_{X_i > d_0}  \notag \\
\label{eq:dbtilde_q} & := \dbtilde \bbQ_{n, 1}(\bh) + \dbtilde \bbQ_{n, 2}(\bh) + \dbtilde \bbQ_{n, 3}(\bh) + \dbtilde \bbQ_{n, 4}(\bh)
\end{align}
We now show that $\dbtilde \bbQ_n(\bh)$ uniformly approximate $\tilde \bbQ_n(\bh)$ over compact sets: 
\begin{align}
& \bbE\left[\sup_{\bh \in \bbI} \left|\dbtilde \bbQ_n(\bh) - \tilde \bbQ_n(\bh)\right|\right] \notag \\
& \le \bbE\left[\sup_{\bh \in \bbI}  \left\{\left|\sum_{i=1}^n \tilde \xi_{i, h_1} \left[ \mathds{1}_{X_i \le d_0} -  \mathds{1}_{X_i \le d_0 \wedge d_0 + \frac{h_3}{n}}\right]\right| + \left|\sum_{i=1}^n \tilde \xi_{i, h_2} \left[ \mathds{1}_{X_i > d_0} -  \mathds{1}_{X_i > d_0 \vee d_0 + \frac{h_3}{n}}\right]\right|\right\}\right] \notag \\
& \le 2 \bbE\left[\sup_{\bh \in \bbI} \left|\sum_{i=1}^n \tilde \xi_{i, h_1} \left[ \mathds{1}_{X_i \le d_0} -  \mathds{1}_{X_i \le d_0 \wedge d_0 + \frac{h_3}{n}}\right]\right|\right] \notag \\
& \hspace{10em} + 2 \bbE\left[\sup_{\bh \in \bbI} \left|\sum_{i=1}^n \tilde \xi_{i, h_2} \left[ \mathds{1}_{X_i \le d_0} -  \mathds{1}_{X_i > d_0 \vee d_0 + \frac{h_3}{n}}\right]\right|\right] \notag \\
& \le 2 \bbE\left[\sum_{i=1}^n |\tilde \xi_{i, h_1}| \mathds{1}_{d_0 - \frac{K}{n} < X_i \le d_0}\right] + 2 \bbE\left[\sum_{i=1}^n |\tilde \xi_{i, h_2}|\mathds{1}_{d_0 < X_i \le d_0 + \frac{K}{n} } \right] \notag \\
& \le 2n\bbE\left[|\tilde \xi_{h_1}|\right]\bbP\left( d_0 - \frac{K}{n} < X \le d_0\right) +2n\bbE\left[|\tilde \xi_{h_2}|\right]\bbP\left( d_0 < X \le d_0 + \frac{K}{n} \right)  \notag \\
& \le 2n \times \left[\sqrt{\var\left(\tilde \xi_{h_1}\right)} \times \left(\frac{K}{n}f_X(d_0) + o(n^{-1})\right) + \sqrt{\var\left(\tilde \xi_{h_2}\right)} \times \left(\frac{K}{n}f_X(d_0) + o(n^{-1})\right)\right] \notag \\
& = 2n \times \left[\left(\frac{h_1\sigma_k}{\sqrt{n}} + o(n^{-1/2})\right) \times \left(\frac{K}{n}f_X(d_0) + o(n^{-1})\right) \right. \notag \\
& \qquad \qquad \qquad \left. + \left(\frac{h_2\sigma_k}{\sqrt{n}} + o(n^{-1/2})\right)  \times \left(\frac{K}{n}f_X(d_0) + o(n^{-1})\right)\right] \notag \\
\label{eq:tightness_1}  & = O(n^{-1/2}) = o(1) \,.
\end{align}
Hence from equation \eqref{eq:tightness_1} we conclude: 
\begin{equation}
\label{eq:dbtilde_approx}
\sup_{\bh \in \bbI} \left|\dbtilde \bbQ_n(\bh) - \tilde \bbQ_n(\bh)\right| = o_p(1) \,.
\end{equation}
Therefore it is immediate from equation \eqref{eq:finite_dim_tilde} that: 
\begin{equation}
\label{eq:finite_dim_dbtilde}
\left(\dbtilde \bbQ_n(\bh_1), \dots, \dbtilde \bbQ_n(\bh_l)\right)  \overset{\mathscr{L}}{\implies}  \left( \bbQ(\bh_1), \dots,  \bbQ(\bh_l)\right) \,.
\end{equation}
Next, we show that tightness $\left\{\dbtilde \bbQ_n(\bh)\right\}_{n \in \bbN} \,.$ As evident from equation \eqref{eq:dbtilde_q}, it is enough to show tightness of $\left\{\dbtilde \bbQ_{n, i}(\bh)\right\}_{n \in \bbN}$ for $i = 1, 2, 3, 4$. For $i= 1$, the process $\dbtilde \bbQ_{n, 1}(\bh)$ only depends on $h_1$ and hence have continuous paths. Therefore to establish tightness, it is enough to show: 
$$
\lim_{\delta \downarrow 0} \limsup_{n \to \infty} \bbE\left[\sup_{\substack{|h_{1, 1}| \vee |h_{1, 2}| \le K \\ |h_{1,1} - h_{1, 2}| \le \delta}}\sum_{i=1}^n  \left|\tilde\xi_{i, h_{1,1}} - \tilde \xi_{i, h_{1,2}}\right|\mathds{1}_{X_i \le d_0}\right] = 0
$$
Towards that end, define a collection of  functions: 
$$
\cF_{1, \delta} = \left\{f_{h_{1, 1}, h_{1, 2}}: \left|h_{1,1}\right| \vee \left|h_{1,2}\right| \le K, |h_{1,1} - h_{1,2}| \le \delta\right\} \,,
$$ where the individual functions $f_{h_{1,1}, h_{1,2}}$ is defined as: 
\begin{align*}
f_{h_{1, 1}, h_{1, 2}}(X, \eps) & = \left\{\left(\tilde H_k\left(\xi_i +  \frac{h_1}{\sqrt{n}}\right) - \tilde H_k\left(\xi_i + \frac{h_2}{\sqrt{n}}\right)\right) \right. \\
& \qquad \qquad \qquad \left. - \bbE\left[\left(\tilde H_k\left(\xi_i +  \frac{h_1}{\sqrt{n}}\right) - \tilde H_k\left(\xi_i+  \frac{h_2}{\sqrt{n}}\right)\right)\right]\right\} \mathds{1}_{X \le d_0} \,.
\end{align*}
Clearly $\cF_{1,\delta}$ has finite VC dimension. Also from equation \eqref{eq:consistency_eq1} we have: 
\begin{align*}
& \left|\tilde H_k\left(\xi_i +  \frac{h_{1,1}}{\sqrt{n}}\right) - \tilde H_k\left(\xi_i + \frac{h_{1,2}}{\sqrt{n}}\right)\right| \\
& \qquad \le \frac{k+1}{k}\left[4k\frac{\left|h_{1,1} - h_{1,2}\right|}{\sqrt{n}} + \frac12\frac{\left(h_{1,1} - h_{1,2}\right)^2}{n} + \left(\frac{|h_{1,1}- h_{1,2}|}{\sqrt{n}} \wedge 2k\right)^2\right] \le C_k\frac{\delta}{\sqrt{n}} \,.
\end{align*}
Therefore, the envelope function can be taken as: 
$$
F_{1, \delta}(X, \eps) = \frac{2C_k\delta}{\sqrt{n}} \,.
$$
Hence using Lemma 2.14.1 of \cite{vdvw96} we conclude: 
$$
 \limsup_{n \to \infty} \bbE\left[\sup_{\substack{|h_{1,1}| \vee |h_{1,2}| \le K \\ |h_{1,1} - h_{1,2}| \le \delta}}\sum_{i=1}^n  \left|\tilde\xi_{i, h_{1,1}} - \tilde \xi_{i, h_{1,2}}\right|\mathds{1}_{X_i \le d_0}\right]   \lesssim \delta
$$
which established tightness of $\dbtilde \bbQ_{n, 1}(\bh)$. The proof of tightness of $\dbtilde \bbQ_{n, 4}(\bh)$ is similar and hence skipped. Finally we show the tightness of $\dbtilde \bbQ_{n, 23}(\bh) = \dbtilde \bbQ_{n, 2}(\bh) + \dbtilde \bbQ_{n, 3}(\bh)$. As these terms only depend on $h_3$ which has cadlag paths, we use equation (13.14) of Theorem 13.5 from \cite{billingsley2013convergence} with $\beta = 1/2$, $\alpha = 1$ and $F(x) = Cx$ for some constant $C$. Fix $h_{3, 1} < 0 < h_{3, 2} < h_{3,3}$. The other cases (i.e. say $0 < h_{3, 1} < h_{3, 2} < h_{3,3}$) are similar and hence skipped.
\begin{align*}
& \bbE\left[\left|\dbtilde \bbQ_{n, 23}(h_{3,1}) - \dbtilde \bbQ_{n, 23}(h_{3,2})\right|\left|\dbtilde \bbQ_{n, 23}(h_{3, 2}) - \dbtilde \bbQ_{n, 23}(h_{3,3})\right|\right] \\
& = \bbE\left[\left|\sum_{i=1}^n \left(\tilde H_k\left(\xi_i + (\alpha_0 - \beta_0)\right) - \tilde H_k(\xi_i)\right)\left[\mathds{1}_{d_0 + \frac{h_{3,1}}{n} \le X_i < d_0} - \mathds{1}_{d_0 \le X_i < d_0 + \frac{h_{3,2}}{n}}\right]\right| \right. \\
& \qquad \qquad \qquad \left. \times \left|\sum_{i=1}^n\left(\tilde H_k\left(\xi_i + (\alpha_0 - \beta_0)\right) - \tilde H_k(\xi_i)\right)\left[\mathds{1}_{d_0 \le X_i < d_0 + \frac{h_{3,2}}{n}} - \mathds{1}_{d_0 \le X_i < d_0 + \frac{h_{3,3}}{n}}\right]\right|\right] \\
& = \bbE\left[\sum_{i=1}^n \left|\tilde H_k\left(\xi_i + (\alpha_0 - \beta_0)\right) - \tilde H_k(\xi_i)\right|\left[\mathds{1}_{d_0 + \frac{h_{3,1}}{n} \le X_i < d_0} + \mathds{1}_{d_0 \le X_i < d_0 + \frac{h_{3,2}}{n}}\right] \right. \\
& \qquad \qquad \qquad \left. \times \sum_{i=1}^n\left|\tilde H_k\left(\xi_i + (\alpha_0 - \beta_0)\right) - \tilde H_k(\xi_i)\right|\left[\mathds{1}_{d_0 \le X_i < d_0 + \frac{h_{3,2}}{n}} + \mathds{1}_{d_0 \le X_i < d_0 + \frac{h_{3,3}}{n}}\right]\right] \\
& \le \frac{k+1}{k}\left(4(\alpha_0 - \beta_0) k + \frac12(\alpha_0 - \beta_0)^2 + \left((\alpha_0 - \beta_0) \wedge 2k\right)^2\right)  \times \\
& \qquad \qquad \qquad \bbE\left[\left(\sum_{i=1}^n\mathds{1}_{d_0 + \frac{h_{3,1}}{\sqrt{n}} \le X_i < d_0 + \frac{h_{3,2}}{n}}\right) \times \left(\sum_{i=1}^n\mathds{1}_{d_0 + \frac{h_{3,2}}{n} \le X_i < d_0 + \frac{h_{3,3}}{n}}\right)\right]\\
& = C_{k, \theta_0} \sum_{i \neq j} \bbE\left[\mathds{1}_{d_0 + \frac{h_{3,1}}{n} \le X_i < d_0 + \frac{h_{3,2}}{n}} \times \mathds{1}_{d_0 + \frac{h_{3,2}}{n} \le X_j < d_0 + \frac{h_{3,3}}{n}}\right] \\
& \le C_{k, \theta_0}  \times n^2 \times \bbP\left(d_0 + \frac{h_{3,1}}{n} \le X < d_0 + \frac{h_{3,2}}{n}\right) \times \bbP\left(d_0 + \frac{h_{3,2}}{n} \le X < d_0 + \frac{h_{3,3}}{n}\right) \\
& \le C_{k, \theta_0} C^2 \times (h_{3,2} - h_{3, 1}) \times (h_{3,3} - h_{3,2}) \hspace{0.2in} [C = \max_x f_X(x)]\\
& \le C_{k, \theta_0} C^2 \times (h_{3,3} - h_{3,1})^2
\end{align*}
This completes the proof of tightness of $\dbtilde \bbQ_n$. Hence using equation \eqref{eq:finite_dim_dbtilde} we conclude: 
$$
\dbtilde \bbQ_n \vert_{\bbI} \overset{\mathscr{L}}{\implies} \bbQ\vert_{\bbI} \,. 
$$
which, along with equation \eqref{eq:dbtilde_approx} implies: 
\begin{equation}
\label{eq:tilde_convergence}
\tilde \bbQ_n \vert_{\bbI} \overset{\mathscr{L}}{\implies} \bbQ\vert_{\bbI} \,.
\end{equation}
Finally, from the decomposition in equation \eqref{eq:decomposition} and combining our conclusions from equation \eqref{eq:remainder_conv}, \eqref{eq:mean_conv} and equation \eqref{eq:tilde_convergence} we have: 
\begin{equation}
\label{eq:final_conv}
\bbQ_n\vert_{\bbI} \overset{\mathscr{L}}{\implies} \bbQ\vert_{\bbI} \,.
\end{equation}
Next and last step is to invoke argmin continuity mapping theorem to say that: 
$$
\widehat \bh_n = \left(\sqrt{n}\left(\hat \alpha - \alpha_0\right), \sqrt{n}\left(\hat \beta - \beta_0\right), n\left(\hat d - d_0\right)\right) \overset{\mathscr{L}}{\implies} \sargmin_{\bh \in \bbR^3}\bbQ(\bh)  \,.
$$
This will complete the proof. Following the proof of Lemma 3.2 of \cite{lan2009change}, all we need to establish the joint asymptotic tightness of $\left\{\left(\bbQ_n(\bh), \bbJ_n(\bh)\right)\right\}_{n \in \bbN}$ where $\bbJ_n(\bh)$ is the jump process corresponding to $\bbQ_n(\bh)$, i.e.  
$$
\bbJ_n(\bh) = \text{sign}(h_3) \sum_{i=1}^n \left[\mathds{1}_{X_i \le d_0 + \frac{h_3}{n}} - \mathds{1}_{X_i \le d_0}\right] \,.
$$
As we have already established the tightness of $\{\bbQ_n(\bh)\}$, we only need to establish the tightness of $\{\bbJ_n(\bh)\}$. The proof is very similar to the proof of Lemma 3.2 of \cite{lan2009change} and skipped here for the brevity.

\subsection{Proofs of Theorem \ref{cor:l1_loss} and \ref{cor:l2_loss}}
The proofs of Theorem \ref{cor:l1_loss} and \ref{cor:l2_loss} are similar to that of Theorem \ref{thm:asymptotic_huber} by replacing $\tilde H_k$ with $\ell_1$ loss function and $\ell_2$ loss function respectively and hence skipped.

%Now we bound the tail of the binomial distribution we use the the following Lemma: 
%
%\begin{lemma}[Lemma 4.7.2 of Ash's Information Theory]
%\label{lem:binomial_tail}
%Suppose $X \sim Bin(n, p)$. Then for any $k \le np$ we have: 
%$$
%\bbP(X \le k) \ge \frac{1}{\sqrt{8k(1 - k/n)}}\exp{\left(-nD_{KL}\left(\frac{k}{n} \vert\vert p\right)\right)}
%$$
%where $D_{KL}(a \vert \vert b)$ is the Kullback-Liebler divergence between two Bernoulli's i.e. 
%$$
%D_{KL}(a \vert \vert b) = a \log{\left(\frac{a}{b}\right)} + (1-a) \log{\left(\frac{1-a}{1-b}\right)} \,.
%$$
%\end{lemma}

\subsection{Proof of Theorem \ref{thm:mlb_growing}}
Here we assume $\alpha_0 = 0, \beta_0 = 1$ is known and derive the bound on the $\hat d$. When $\alpha_0, \beta_0$ is not known, then the problem becomes harder and rate of convergence obviously can not be faster. Our proof is based on the proof of Theorem 5 of \cite{massart2006risk}.  Consider $\mathscr{A}$ to be the set of all half-spaces, i.e.
$$
\mathscr{A} = \left\{x^{\top}d > 0\right\}_{d \in S^{p-1}} \,.
$$
Our model is $X \sim P$ and: 
$$
Y = \mathds{1}_{X^{\top}d > 0} + \xi
$$
where $\xi \sim \cN(0, 1)$. Now the class of hyperplanes $\mathscr{A}$ has VC dimension $p$. From the properties of the hyperplane we know that, given any $N> p$ (not to be confused with sample size $n$) there exist $x_1, \dots, x_N \in \bbR^p$ such that $\mathscr{A}$ shatters all subsets of $\left\{x_1, \dots, x_N\right\}$ with cardinality $k \le \lfloor p/2 \rfloor := V$ (e.g. see \cite{edelsbrunner2012algorithms}). Define $\Theta_{N, V}$ to be the collection of all such $d \in S^{p-1}$ which shatters all subsets of length $V$ of $\left\{x_1, \dots, x_N\right\}$.  Hence we have: 
$$
\left\{\left(\mathds{1}_{x_1^{\top}d > 0}, \dots, \mathds{1}_{x_N^{\top}d > 0}\right)\right\}_{d \in \Theta_{N, V}} := B
$$
where: 
$$
B = \left\{0, 1\right\}_{N, V} = \left\{b \in \{0, 1\}^N: \sum_{i=1}^N b_i = V\right\} \,.
$$
Henceforth for any $d \in \Theta_{N, V}$ we denote by $b_{d}$ to be the corresponding unique $b \in B$. Define $\mu$ to be the uniform measure on $\{x_1, \dots, x_N\}$ and for any $d \in \Theta_{N, V}$ define: 
$$
Y = \mathds{1}_{X^{\top}d > 0} + \xi \,.
$$
The loss function we use here is the squared error loss defined as: 
\begin{align*}
    \ell\left(d, d_0\right) & = \bbE\left[(Y - \mathds{1}_{X^{\top}d > 0})^2 - (Y - \mathds{1}_{X^{\top}d_0 > 0})^2\right] \\
    & = \bbE_X\left[\left| \mathds{1}_{X^{\top}d > 0} - \mathds{1}_{X^{\top}d_0 > 0}\right|\right] \\
    & = \left\|\mathds{1}_{X^{\top}d > 0} - \mathds{1}_{X^{\top}d_0 > 0}\right\|_{L_1(P)} \,.
\end{align*}
The minimax risk is defined as: 
\begin{align}
\cR_n & = \inf_{\hat d} \sup_{d \in S^{p-1}}\bbE_d\left[\ell(\hat d, d)\right] \notag \\
& \ge \inf_{\hat d} \sup_{d \in \Theta_{N, V}}\bbE_d\left[\ell(\hat d, d)\right] \notag \\
& = \inf_{\hat d} \sup_{d \in \Theta_{N, V}}\bbE_d\left[\left\|\mathds{1}_{X^{\top}\hat d > 0} - \mathds{1}_{X^{\top}d > 0}\right\|_{L_1(\mu)}\right] \notag \\
\label{eq:lower_bound_1} & = \frac{1}{N}\inf_{\hat d} \sup_{d \in \Theta_{N, V}}\bbE_d\left[\sum_{i=1}^N\left|\mathds{1}_{x_i^{\top}\hat d > 0} - \mathds{1}_{x_i^{\top}d > 0}\right|_{L_1(\mu)}\right]
\end{align}
Now for any $\hat d \in S^{p-1}$, define $\hat d_{\new}$ as: 
$$
\hat d_{\new} = \argmin_{d_* \in \Theta_{N, V}} \left\|\mathds{1}_{X^{\top}\hat d > 0} - \mathds{1}_{X^{\top}d_* > 0}\right\|_{L_1(\mu)} \,.
$$
Then we have for any $d \in \Theta_{N, V}$:
\begin{align*}
    \left\|\mathds{1}_{X^{\top}\hat d_\new > 0} - \mathds{1}_{X^{\top}d > 0}\right\|_{L_1(\mu)} & =  \left\|\mathds{1}_{X^{\top}\hat d_\new > 0} -  \mathds{1}_{X^{\top}\hat d > 0} +  \mathds{1}_{X^{\top}\hat d > 0} - \mathds{1}_{X^{\top}d > 0}\right\|_{L_1(\mu)} \\
    & \le \left\|\mathds{1}_{X^{\top}\hat d_\new > 0} -  \mathds{1}_{X^{\top}\hat d > 0}\right\|_{L_1(\mu)} + \left\|\mathds{1}_{X^{\top}\hat d > 0} - \mathds{1}_{X^{\top}d > 0}\right\|_{L_1(\mu)} \\
    & \le 2\left\|\mathds{1}_{X^{\top}\hat d > 0} - \mathds{1}_{X^{\top}d > 0}\right\|_{L_1(\mu)} \,.
\end{align*}
Putting this bound in equation \eqref{eq:lower_bound_1} we obtain: 
\begin{align}
    \cR_n & \ge \frac{1}{N}\inf_{\hat d} \sup_{d \in \Theta_{N, V}}\bbE_d\left[\sum_{i=1}^N\left|\mathds{1}_{x_i^{\top}\hat d > 0} - \mathds{1}_{x_i^{\top}d > 0}\right|_{L_1(\mu)}\right] \notag \\
    & \ge \frac{1}{2N}\inf_{\hat d} \sup_{d \in \Theta_{N, V}}\bbE_d\left[\sum_{i=1}^N\left|\mathds{1}_{x_i^{\top}\hat d_\new > 0} - \mathds{1}_{x_i^{\top}d > 0}\right|_{L_1(\mu)}\right] \notag \\
    \label{eq:lower_bound_2}& = \frac{1}{2N} \inf_{\hat d \in \Theta_{N,V}} \sup_{d \in \Theta_{N, V}}\bbE_d\left[\sum_{i=1}^N\left|\mathds{1}_{x_i^{\top}\hat d > 0} - \mathds{1}_{x_i^{\top}d > 0}\right|_{L_1(\mu)}\right]
\end{align}
Next note that: 
$$
\sum_{i=1}^N\left|\mathds{1}_{x_i^{\top}\hat d > 0} - \mathds{1}_{x_i^{\top}d > 0}\right| = d_{H}\left(b_{\hat d}, b_{d}\right)
$$
where $d_H$ is the Hamming distance. As $\Theta_{N,V}$ has a bijection with $B$ we can write equation \eqref{eq:lower_bound_2} as: 
\begin{equation}
    \label{eq:lower_bound_3}
    \cR_n \ge \frac{1}{2N} \inf_{\hat b} \sup_{b \in B}\bbE_b\left[d_{H}\left(\hat b, b\right)\right] \ge \frac{1}{2N} \inf_{\hat b \in  \cD} \sup_{b \in \cD}\bbE_b\left[d_{H}\left(\hat b, b\right)\right]
\end{equation}
for any subset $\cD \subseteq B$. We now choose $\cD$ carefully. Note that for any $N \ge 4V$ (i.e. $N \ge 2p$), we can choose $\cD$ such that (Lemma 8 of \cite{reynaud2003adaptive}): 
\begin{enumerate}
    \item $d_H\left(b, b'\right) > \frac{V}{2}$ for all $b \neq b' \in \cD$. 
    \item $\log{|\cD|} \ge \rho V\log{\left(\frac{N}{V}\right)}$ with $\rho = 0.233$. 
\end{enumerate}
Using this we modify equation \eqref{eq:lower_bound_3} as follows: 
\begin{align}
    \cR_n & \ge \frac{1}{2N} \inf_{\hat b \in  \cD} \sup_{b \in \cD}\bbE_b\left[d_{H}\left(\hat b, b\right)\right] \notag \\
     & = \frac{1}{2N} \inf_{\hat b \in  \cD} \sup_{b \in \cD}\bbE_b\left[d_{H}\left(\hat b, b\right)\mathds{1}_{\hat b \neq b}\right] \notag \\
    & \ge \frac{V}{4N} \inf_{\hat b \in  \cD} \sup_{b \in \cD}\bbP_b\left(\hat b \neq b\right) \hspace{0.2in} [\text{From point 1. above}] \notag \\
    & = \frac{V}{4N} \inf_{\hat b \in  \cD} \sup_{b \in \cD}\left(1 - \bbP_b\left(\hat b = b\right)\right) \notag \\
    \label{eq:lower_bound_4}& = \frac{V}{4N} \inf_{\hat b \in  \cD}\left(1 - \min_{b \in \cD}\bbP_b\left(\hat b = b\right)\right)
\end{align}
Now to further bound the above equation, we use the following lemma (see \cite{birge2005new}): 
\begin{lemma}
\label{lem:kl_bound}
Let $m \ge 1$, $(P_i)_{0 \le i \le m}$ be a family of probability distributions and $(A_i)_{0 \le i \le m}$ be a family of disjoint events. Let $a = \min_{0 \le i \le m} P_i(A_i)$. Then setting: 
$$
\bar{\mathscr{K}} = \frac1m \sum_{i=1}^m \mathscr{K}(P_i, P_0)
$$ 
where $\mathscr{K}$ is the Kullback-Liebler divergence, we have: 
$$
a \le \alpha \vee \left(\frac{\bar{\mathscr{K}}}{\log{(1+m)}}\right) \,.
$$
where $\alpha = 0.71$. 
\end{lemma}
We now use this bound in equation \eqref{eq:lower_bound_4}. Fix any $b_0 \in \cD$. Define $A_i = \left\{\hat b = b_i\right\}$ for all $b_i \in \cD$ which are disjoint events. Hence using Lemma \ref{lem:kl_bound} we obtain: 
\begin{equation}
    \label{eq:lower_bound_5}
    \min_{b \in \cD}\bbP_b\left(\hat b = b\right) \le \alpha \vee \left(\frac{\bar{\mathscr{K}}}{\log{|\cD|}}\right) 
\end{equation}
where: 
$$
\bar{\mathscr{K}} = \frac{1}{|\cD| - 1}\sum_{b \in \cD, b \neq b_0} \mathscr{K}\left(P_b^{\otimes n}, P_{b_0}^{\otimes n}\right) = \frac{n}{|\cD| - 1}\sum_{b \in \cD, b \neq b_0} \mathscr{K}\left(P_b, P_{b_0}\right) \,.
$$
Now note that for any $d_1, d_2$ we have: 
\begin{align*}
    \mathscr{K}\left(P_{d_1}, P_{d_2}\right) & = \bbE_X\left[\mathscr{K}\left(P_{d_1}(Y \mid X), P_{d_2}(Y \mid X)\right)\right] \\
    &= \bbE_X\left[\left(\mathds{1}_{X^{\top}d_1 > 0} - \mathds{1}_{X^{\top}d_2 > 0}\right)^2\right] \\
    & = \bbE_X\left[\left|\mathds{1}_{X^{\top}d_1 > 0} - \mathds{1}_{X^{\top}d_2 > 0}\right|\right] \\
    & = \frac{1}{N}\sum_{i=1}^N \left|\mathds{1}_{x_i^{\top}d_1 > 0} - \mathds{1}_{x_i^{\top}d_2 > 0}\right| \\
    & = \frac1N d_H(b_{d_1}, b_{d_2}) \,.
\end{align*}
Also by definition for any two $b, b' \in \cD$ we have $d_H(b, b') \le 2V$. Hence we have: 
$$
\bar{\mathscr{K}} \le \frac{2Vn}{N} \,.
$$
This bound along with equation \eqref{eq:lower_bound_5} modifies the bound of equation \eqref{eq:lower_bound_4} as: 
\begin{align}
    \cR_n & \ge \frac{V}{4N} \inf_{\hat b \in  \cD}\left(1 - \min_{b \in \cD}\bbP_b\left(\hat b = b\right)\right) \notag \\
    & \ge \frac{V}{4N} \left(1 - \left(\alpha \vee \frac{2Vn}{N\log{|\cD|}}\right)\right) \notag \\
    \label{eq:lower_bound_6} & = \frac{V(1-\alpha)}{4N}
\end{align}
when, 
\begin{align*}
    & \alpha \ge \frac{2Vn}{N\log{|\cD|}} \,.
\end{align*}
which holds if: 
\begin{align*}
    & \alpha > \frac{2Vn}{N\rho V\log{(N/V)}}
\end{align*}
i.e. if: 
\begin{equation}
    \label{eq:lower_bound_7}
    N\log{\left(\frac{N}{V}\right)} \ge \frac{2n}{\alpha \rho} \,.
\end{equation}
which is satisfied, if for example: 
$$
N = \left\lfloor \frac{4n}{\rho \alpha \left(1 + \log{\left(\frac{n}{V}\right)}\right)}\right\rfloor. 
$$
Using this in equation \eqref{eq:lower_bound_6} we conclude:
$$
\cR_n \gtrsim \frac{V}{n}\left(1 + \log{\left(\frac{n}{V}\right)}\right) \asymp \frac{p}{n}\left(1 + \log{\left(\frac{n}{p}\right)}\right)\,.
$$
as $V = \lfloor p/2 \rfloor$. We finally need to verify $N > 4V$ and that it satisfies equation \eqref{eq:lower_bound_7}. The first one is obviously true for all large $n$ as $n/p \to \infty$. For the second one, lets forget the $\lfloor \cdot \rfloor$ in the definition of $N$ for the time being as it will not affect asymptotically. Then: 
\begin{align*}
    & N\log{\left(\frac{N}{V}\right)} \ge \frac{2n}{\alpha \rho} \\
    \iff & \frac{4n}{\rho \alpha \left(1 + \log{\left(\frac{n}{V}\right)}\right)}\log{\left( \frac{4n}{V\rho \alpha \left(1 + \log{\left(\frac{n}{V}\right)}\right)}\right)} \ge \frac{2n}{\alpha \rho} \\
    \iff & \frac{1}{\left(1 + \log{\left(\frac{n}{V}\right)}\right)}\log{\left( \frac{4n}{V\rho \alpha \left(1 + \log{\left(\frac{n}{V}\right)}\right)}\right)} \ge \frac12 \\
    \iff & \frac{1}{\left(1 + \log{\left(\frac{n}{V}\right)}\right)}\left[\log{\left(\frac{n}{V}\right)} + \log{\left(\frac{4}{\rho \alpha}\right)} - \log{\left(1 + \log{\left(\frac{n}{V}\right)}\right)}\right] \ge \frac{1}{2}
\end{align*}
As $n/V \to \infty$, LHS converges to 1 and eventually $> 1/2$. Therefore the choice of $N$ is valid for all  large $n$.

\subsection{Proof of Theorem \ref{thm:huber_high}}
We first assume that Assumption \ref{assm:wedge} holds globally and our parameter space $\Omega = \Omega_{\alpha} \times \Omega_\beta$ for $(\alpha_0, \beta_0)$ is such that: 
$$
\min_{\alpha \in \Omega_\alpha} |\alpha - \beta_0| \wedge \min_{\beta \in \Omega_\beta} |\beta - \alpha_0|  \ge \delta > 0 \,.
$$
This is just to avoid the issue of consistency. One can relax this assumption with an additionally showing that the estimators are consistent. To prove Theorem \ref{thm:huber_high} we use Theorem A.3 of \cite{mukherjee2019non}. To match our notation with that theorem, here our loss function $\gamma(\theta, \cdot)$ is: 
\begin{align*}
\gamma(\theta, (X, \xi)) & = \tilde H_k\left(\xi + \alpha_0 \mathds{1}_{X^{\top}d_0 \le 0} + \beta_0 \mathds{1}_{X^{\top}d_0 > 0} -  \alpha \mathds{1}_{X^{\top}d \le 0} - \beta \mathds{1}_{X^{\top}d > 0} \right) - \tilde H_k(\xi) \\
 & = \left(\tilde H_k\left(\xi +  \alpha_0 - \alpha\right) -  \tilde H_k(\xi)\right)\mathds{1}_{X^{\top}d \vee X^{\top}d_0 \le 0} \\
& \qquad \left(\tilde H_k\left(\xi +  \alpha_0 - \beta\right) - \tilde H_k(\xi)\right)\mathds{1}_{X^{\top}d_0 \le 0 < X^{\top}d } \\
& \qquad \qquad +  \left(\tilde H_k\left(\xi +  \beta_0 - \alpha\right) -\tilde H_k(\xi)\right)\mathds{1}_{X^{\top}d \le 0 < X^{\top}d_0} \\
& \qquad \qquad \qquad +   \left(\tilde H_k\left(\xi +  \beta_0 - \beta\right) - \tilde H_k(\xi)\right)\mathds{1}_{X^{\top}d \wedge X^{\top}d_0 > 0} 
\end{align*}
It is immediate from the definition that $\gamma(\theta_0, (X, \xi)) = 0$. Also from equation \eqref{eq:consistency_eq1} we know $\gamma(\theta, \cdot)$ is uniform bounded by some constant only depending on $k$ and the width of $\Omega$. Following similar arguments used in the proof of Theorem \ref{thm:l2}, we obtain: 
\begin{align}
\label{eq:loss_dist_lower_bound}
\ell(\theta, \theta_0) & = \bbE[\gamma(\theta, (X, \xi))] \ge  c_k \ \dist^2(\theta, \theta_0) \,,
\end{align}
for some constant $c_k$ (independent of $n$), where the $\dist$ function is:
$$
\dist(\theta, \theta_0) = \sqrt{(\alpha - \alpha_0)^2 + (\beta - \beta_0)^2 +\bbP\left(\s(X^{\top}d) \neq \s(X^{\top}d_0)\right)} \,.
$$
Moreover, from equation \eqref{eq:consistency_eq1} we have: 
\begin{align*}
\var\left(\gamma(\theta, (X, \xi))\right) & \le \bbE[\gamma^2(\theta, (X, \xi))] \le C_k \ \dist^2(\theta, \theta_0) \,.
\end{align*}
Hence this semi-metric $\dist$ satisfies conditions of Theorem A.3 of \cite{mukherjee2019non} with $\omega(x) = x$. Next we need to bound the modulus of continuity: 
\begin{align*}
\sqrt{n}\bbE\left[\sup_{\substack{\theta: f_\theta \in \cF_m \\ \dist(\theta, \theta_m) \le \eps}}\left|\left(\bbP_n - P\right)\left(\gamma(\theta, (X, \xi)) - \gamma(\theta_m, (X, \xi))\right)\right| \right] 
\end{align*} 
Note that another application of equation \eqref{eq:consistency_eq1} yields: 
$$
\sup_{\substack{\theta: f_\theta \in \cF_m \\ \dist(\theta, \theta_m) \le \eps}} \bbE\left[\left(\gamma(\theta, (X, \xi)) - \gamma(\theta_m, (X, \xi))\right)^2\right] \lesssim \eps^2 \,.
$$
Hence applying Theorem 8.7 of \cite{sen2018gentle} we have: 
\begin{align}
\label{eq:moc_highdim_1}
\sqrt{n}\bbE\left[\sup_{\substack{\theta: f_\theta \in \cF_m \\ \dist(\theta, \theta_m) \le \eps}}\left\|\left(\bbP_n - P\right)\left(\gamma(\theta, (X, \xi)) - \gamma(\theta_m, (X, \xi))\right)\right\| \right] & \lesssim \eps\sqrt{V_m\log{\frac{1}{\eps}}} \vee \frac{V_m}{\sqrt{n}}\log{\frac{1}{\eps}} \notag \\
& := \psi_m(\eps) \,.
\end{align} 
So a value of $\eps_m$ that satisfies $\sqrt{n}\eps_m^2 \ge \phi_m(\eps_m)$ is: 
$$
\eps_m = \frac{V_m}{n}\log{\frac{n}{V_m}} \,.
$$
Therefore we can take $x_m = V_m\log{(n/V_m)}$ and as $\omega(x) = x$ we have $b(n) = 1$ (see Theorem A.3 of \cite{mukherjee2019non} for the exact expression of $\omega(x), b(n)$). Note that, as we are assuming $(s\log{p})/n \to 0$ (Assumption \ref{assm:sparsity}), it is sufficient to search among the models with $1 \le m \le C\lfloor (n/\log{p}) \rfloor$ for any constant $C$. We take $C = 1/4$ here. With these choices, the value of $\Sigma$ (as defined in Theorem A.3 of \cite{mukherjee2019non}) is: 
\begin{align*}
\Sigma & = \sum_{i=1}^{\frac14 \lfloor \frac{n}{\log{p}} \rfloor} e^{-V_i\log{\frac{n}{V_i}}} \\
& \le \lim_{n \to \infty} \sum_{i=1}^{\frac14 \lfloor \frac{n}{\log{p}} \rfloor} e^{-V_i\log{\frac{n}{V_i}}} < \infty \,.
\end{align*}
Hence, an application Theorem A.3 of \cite{mukherjee2019non} yields: 
$$
\bbP\left(\ell(\hat \theta, \theta_0) > C\pen(s) + t \frac{C_1}{n}\right) \le \Sigma e^{-t} \,.
$$
This along with the value of $\pen(s)$ from equation \eqref{eq:penalty_huber} of the main paper and the lower bound equation \eqref{eq:loss_dist_lower_bound} completes the first part of the proof, i.e. 
$$
\dist^2\left((\hat \alpha_{\init}, \hat \beta_{\init}, \hat d), (\alpha_0, \beta_0, d_0)\right) = O_p\left(\frac{V_s}{n}\log{\frac{n}{V_s}}\right) \,.
$$  
The acceleration of the rate of $\hat \alpha, \hat \beta$ via replacing $d_0$ by $\hat d$ in the model equation is exactly same as that of Theorem \ref{thm:l2} and hence skipped.

\subsection{Proof of Theorem \ref{thm:l2_high}}
Here again we assume for technical simplicity that the wedge assumption (Assumption \ref{assm:wedge}) is valid on entire $S^{p-1}$, although all our arguments can be extended to the case where the assumption is valid only locally along with a separate argument for the consistency of the estimator. We use the same notations as of Theorem \ref{thm:huber_high} through out the proof. Define 
$$
y^2_{m} = 2\kappa \frac{V_m\|\xi\|_{n, 2}}{n}\log{\frac{n}{V_m}}
$$ 
for all $m \in \mathcal{M}$, where $V_m$ is the VC dimension of model $m$. For any such model $m$ we obtain $\hat \theta_m$ as: 
\begin{align*}
\hat \theta_m & = \argmin_{\theta \in \Omega \times S^{p-1}_m} \bbP_n \left\{(Y - \alpha\mathds{1}_{X^{\top}d \le 0} - \beta\mathds{1}_{X^{\top}d > 0})^2 - \xi^2 \right\} \\
& = \argmin_{\theta \in \Omega \times S^{p-1}_m} \bbP_n \left\{\xi \left(\alpha_0 \mathds{1}_{X^{\top}d_0 \le 0} + \beta_0\mathds{1}_{X^{\top}d_0 > 0} - \alpha \mathds{1}_{X^{\top}d \le 0}  - \beta \mathds{1}_{X^{\top}d > 0}\right) \right. \\
& \qquad \qquad \qquad \qquad \qquad \left. + \frac12 \left(\alpha_0 \mathds{1}_{X^{\top}d_0 \le 0} + \beta_0\mathds{1}_{X^{\top}d_0 > 0} - \alpha \mathds{1}_{X^{\top}d \le 0}  - \beta \mathds{1}_{X^{\top}d > 0}\right)^2 \right\} \\
& := \argmin_{\theta \in \Omega \times S^{p-1}_m} \bbP_n f_{\theta} \\
& := \argmin_{\theta \in \Omega \times S^{p-1}_m} \bbP_n \left(f_{\theta, 1} + f_{\theta, 2}\right)
\end{align*}
where the functions $f_{d, 1}, f_{d, 2}$ are defined as: 
\begin{align*}
    f_{\theta, 1} & = \xi \left(\alpha_0 \mathds{1}_{X^{\top}d_0 \le 0} + \beta_0\mathds{1}_{X^{\top}d_0 > 0} - \alpha \mathds{1}_{X^{\top}d \le 0}  - \beta \mathds{1}_{X^{\top}d > 0}\right) \,,\\
    f_{\theta, 2} & = \frac12 \left(\alpha_0 \mathds{1}_{X^{\top}d_0 \le 0} + \beta_0\mathds{1}_{X^{\top}d_0 > 0} - \alpha \mathds{1}_{X^{\top}d \le 0}  - \beta \mathds{1}_{X^{\top}d > 0}\right)^2 \,.
\end{align*}
The loss function used here is: 
$$
\ell(\theta, \theta_0) = \bbP f_{\theta} \gtrsim \dist^2(\theta, \theta_0) \,,
$$
for all $\theta \in \Omega \in S^{p-1}$ via the global version of Assumption \ref{assm:wedge}. 
From the definition of $\hat m$ we have: 
\begin{align*}
    \bbP_n f_{\hat \theta_{\hat m}} + \text{pen}(\hat m) & \le \bbP_n f_{\hat \theta_{s_0}} + \text{pen}(s_0) \\
    & \le \bbP_n f_{\theta_{s_0}} + \text{pen}(s_0) := \text{pen}(s_0) \,.
\end{align*}
Using this we can bound the loss function: 
\allowdisplaybreaks
\begin{align*}
\ell(\hat \theta_{\hat m}, \theta_0) & = \bbP f_{\hat \theta_{\hat m}} \notag \\
& = \left(\bbP - \bbP_n\right)f_{\hat \theta_{\hat m}} + \bbP_n f_{\hat \theta_{\hat m}} \\
& =  \left(\bbP - \bbP_n\right)f_{\hat \theta_{\hat m}} + \bbP_n f_{\hat \theta_{\hat m}} + \pen(\hat m) - \pen(\hat m) \\
& \le \left(\bbP - \bbP_n\right)f_{\hat \theta_{\hat m}} + \pen(s_0) - \pen(\hat m) \\
& = \frac{\left(\bbP - \bbP_n\right)f_{\hat \theta_{\hat m}}}{\ell(\hat \theta_{\hat m}, \theta_0) + y_{\hat m}^2}\left(\ell(\hat \theta_{\hat m}, \theta_0) + y_{\hat m}^2\right) + \pen(s_0) - \pen(\hat m) \\
& \le \sup_{\theta \in \Omega \times S^{p-1}_{\hat m}}\frac{\left|\left(\bbP - \bbP_n\right)f_{\theta}\right|}{\ell(\theta, \theta_0) + y_{\hat m}^2}\left(\ell(\hat \theta_{\hat m}, d_0) + y_{\hat m}^2\right) + \pen(s_0) - \pen(\hat m) 
\end{align*}
For the rest of the calculation, define: 
\begin{align*}
\Gamma_m & = \sup_{\theta \in \Omega \times S^{p-1}_{ m}}\frac{\left|\left(\bbP - \bbP_n\right)f_{\theta}\right|}{\ell(\theta, \theta_0) + y_{m}^2} \\
& =  \sup_{\theta \in \Omega \times S^{p-1}_{ m}}\frac{\left|\left(\bbP - \bbP_n\right)(f_{\theta, 1} + f_{\theta, 2})\right|}{\ell(\theta, \theta_0) + y_{m}^2} \\
& \le  \sup_{\theta \in \Omega \times S^{p-1}_{ m}}\frac{\left|\left(\bbP - \bbP_n\right)f_{\theta, 1}\right|}{\ell(\theta, \theta_0) + y_{m}^2} +  \sup_{\theta \in \Omega \times S^{p-1}_{ m}}\frac{\left|\left(\bbP - \bbP_n\right)f_{\theta, 2}\right|}{\ell(\theta, \theta_0) + y_{m}^2} \\
& := \Gamma_{m, 1} + \Gamma_{m, 2} \,.
\end{align*}
Next we try to bound $\Gamma_{\hat m}$. More specifically, we bound $\Gamma_m$ for all $m$ and then use a union bound to bound $\Gamma_{\hat m}$. Note that, as the function class under the consideration of $\Gamma_{m, 2}$ is bounded we can use similar as of the proof of Theorem A.3 of \cite{mukherjee2019non} (i.e. applying Talagrand's inequality and then bound the expectation and variance) to conclude: 
\begin{equation}
\label{eq:exp_bound_1}
\bbP\left(\Gamma_{\hat m, 2} \ge 1/4\right) = o(1) \,.
\end{equation}
For $\Gamma_{m, 1}$, we first decompose it as follows: 
\begin{align*}
\Gamma_{m ,1} & \le \sup_{d \in S^{p-1}_{m}}\frac{\left|\left(\bbP - \bbP_n\right)\left(f_{\theta, 1} - f_{\theta_m, 1}\right)\right|}{\ell(\theta, \theta_0) + y_{m}^2} + \sup_{\theta \in \Omega \times S^{p-1}_{m}} \frac{\left|\left(\bbP - \bbP_n\right)f_{\theta_m, 1}\right|}{\ell(\theta, \theta_0) + y_{m}^2} \\
& \le \sup_{\theta \in \Omega \times S^{p-1}_{m}}\frac{\left|\left(\bbP - \bbP_n\right)\left(f_{\theta, 1} - f_{\theta_m, 1}\right)\right|}{\ell(\theta, \theta_0) + y_{m}^2} + \frac{\left|\left(\bbP - \bbP_n\right)f_{\theta_m, 1}\right|}{\ell(\theta_m, \theta_0) + y_{m}^2} \\
& = \Gamma_{m, 11} + \Gamma_{m, 12} \,.
\end{align*}
Bounding $\bbE[\Gamma_{m, 12}]$ is straight-forward: 
\begin{align*}
    \bbE[\Gamma_{m, 12}] & \le \frac{\sqrt{\var\left(f_{\theta_m ,1}\right)}}{\sqrt{n}\left(\ell(\theta_m, \theta_0) + y_{m}^2\right)} \\
    & \le \frac{\|\xi\|_2 \dist(\theta_m, \theta_0)}{\sqrt{n}\left(\dist^2(\theta_m, \theta_0) + y_{m}^2\right)} \\
    & \le \frac{\|\xi\|_2}{\sqrt{n}} \ \sup_{x \ge 0} \frac{x}{x^2 + y_m^2}  \le \frac{\|\xi\|_2}{2} \frac{1}{\sqrt{n}y_m} \,.
\end{align*}
To bound $\bbE[\Gamma_{m, 1}]$ we use the maximal inequality of the weighted empirical process (see Lemma A.5 of \cite{massart2006risk}), which is a variant of peeling argument. First of all note that, by symmetrization and applying Theorem 8.7 of \cite{sen2018gentle}, we have for any $1 \le k \le n$: 
\begin{align*}
    & \bbE\left[\sup_{\substack{\theta \in \Omega \times S^{p-1}_m \\ \dist(\theta, \theta_m) \le \eps}}\left|\sum_{i=1}^k \xi_i \left(\alpha_m \mathds{1}_{X_i^{\top}d_m \le 0} + \beta_m\mathds{1}_{X_i^{\top}d_m > 0} - \alpha \mathds{1}_{X_i^{\top}d \le 0}  - \beta \mathds{1}_{X_i^{\top}d > 0}\right)\right|\right] \\
    & \qquad \qquad \qquad \lesssim \left(\sigma_{\eps}\sqrt{kV_m\log{\frac{1}{\sigma_{\eps}}}} \vee V_m\log{\frac{1}{\sigma_\eps}}\right) \,.
\end{align*}
where the \emph{wimpy variance} $\sigma_\eps^2$ is defined as: 
\begin{align*}
    \sigma_\eps^2 = \sup_{\substack{\theta: f_\theta \in \cF_m \\ \dist(\theta, \theta_m) \le \eps}} \sigma_\xi^2 \ \bbE\left[\left(\alpha_m \mathds{1}_{X_i^{\top}d_m \le 0} + \beta_m\mathds{1}_{X_i^{\top}d_m > 0} - \alpha \mathds{1}_{X_i^{\top}d \le 0}  - \beta \mathds{1}_{X_i^{\top}d > 0}\right)^2\right] \lesssim \eps^2  \,.
\end{align*}
This, along with Proposition \ref{prop:multiplier_ineq} implies: 
\begin{align*}
& \bbE\left[\sup_{\substack{\theta: f_\theta \in \cF_m \\ \dist(\theta, \theta_m) \le \eps}} \left|\left(\bbP - \bbP_n\right)\left(f_{\theta, 1} - f_{\theta_m, 1}\right)\right|\right] \\
    & = \bbE\left[\sup_{\substack{\theta: f_\theta \in \cF_m \\ \dist(\theta, \theta_m) \le \eps}} \left|\sum_{i=1}^k \xi_i\left(\alpha_m \mathds{1}_{X_i^{\top}d_m \le 0} + \beta_m\mathds{1}_{X_i^{\top}d_m > 0} - \alpha \mathds{1}_{X_i^{\top}d \le 0}  - \beta \mathds{1}_{X_i^{\top}d > 0}\right) \right|\right] \\
    & \qquad \lesssim \left\|\xi\right\|_{2, 1} \eps \sqrt{\frac{V_m}{n}\log{\left(\frac{1}{\eps}\right)}} + 2
    \frac{V_m}{n}\log{\left(\frac{1}{\eps}\right)}\bbE\left[\max_{1 \le i \le n} \left|\xi_i\right|\right] := \frac{\psi_m(\eps)}{\sqrt{n}} \,.
\end{align*}
An application of Lemma A.5 of \cite{massart2006risk} yields: 
$$
\bbE[\Gamma_{m, 1}] \lesssim \frac{\psi_m(2\sqrt{2} y_m)}{\sqrt{n} y_m^2}
$$
which, in turn, yields: 
\begin{align}
E(\Gamma_{m, 1})  & \le \frac{4\phi_{m}(2\sqrt{2}y_{m})}{\sqrt{n}y^2_{m}} + \frac{\|\xi\|_2}{2\sqrt{n}y_{m}} \notag\\
& \le \frac{4\phi_{m}(2\sqrt{2}\epsilon_{m})}{\sqrt{n}y_{m}\epsilon_{m}} + \frac{\|\xi\|_2}{2\sqrt{n}y_{m}}\sqrt{\frac{y_{m}^2}{y^2_{m}}}\notag\\
& \le \frac{8\sqrt{2}\phi_{m}(\epsilon_{m})}{\sqrt{n}y_{m}\epsilon_{m}} + \frac{\|\xi\|_2}{2\sqrt{n}y_{m}}\sqrt{\frac{\phi^2_{m}\left(y_{m}\right)}{y^2_{m}}}\notag\\
& \le \frac{8\sqrt{2}\epsilon_{m}}{y_{m}} + \frac{2}{\sqrt{n}y_{m}}\sqrt{\frac{\phi^2_{m}\left(\epsilon_{m}\right)}{\epsilon^2_{m}}}\notag\\
& \le \frac{8}{\sqrt{\kappa}} + \frac{\|\xi\|_2}{2\sqrt{n}y_{m}}\frac{\phi_{m}\left(\epsilon_{m}\right)}{\epsilon_{m}}\notag \notag\\ 
\label{bound-on-expectation} & \le \frac{8}{\sqrt{\kappa}} + \frac{\|\xi\|_2\epsilon_{m}}{2y_{m}} \le \frac{8+\sqrt{2}}{\sqrt{\kappa}}
\end{align}
which can be made arbitrarily small by making $\kappa$ arbitrarily large. Next, we bound the fluctuation of $V_{m , 1}$ around its mean using Chebychev inequality: 
$$
\bbP\left(\left|\Gamma_{m, 1} - \bbE[\Gamma_{m, 1}]\right| \ge t\right) \le \frac{\var(
\Gamma_{m, 1})}{t^2} \,.
$$
To bound the variance we use Theorem 11.17 along with Theorem 11.1 of \cite{boucheron2013concentration}. To match with their notation for the ease of the readers, we have: 
\begin{align*}
X_{i, \theta} & = \frac{1}{n}\left[\frac{\xi_i\left(\alpha_m \mathds{1}_{X_i^{\top}d_m \le 0} + \beta_m\mathds{1}_{X_i^{\top}d_m > 0} - \alpha \mathds{1}_{X_i^{\top}d \le 0}  - \beta \mathds{1}_{X_i^{\top}d > 0}\right)}{\ell(\theta, \theta_0) + y_m^2} \right. \\
& \qquad \qquad \qquad \left. - \bbE\left(\frac{\xi_i\left(\alpha_m \mathds{1}_{X_i^{\top}d_m \le 0} + \beta_m\mathds{1}_{X_i^{\top}d_m > 0} - \alpha \mathds{1}_{X_i^{\top}d \le 0}  - \beta \mathds{1}_{X_i^{\top}d > 0}\right)}{\ell(\theta, \theta_0) + y_m^2}\right)\right] \\
& = \frac{1}{n}\left[\frac{\xi_i\left(\alpha_m \mathds{1}_{X_i^{\top}d_m \le 0} + \beta_m\mathds{1}_{X_i^{\top}d_m > 0} - \alpha \mathds{1}_{X_i^{\top}d \le 0}  - \beta \mathds{1}_{X_i^{\top}d > 0}\right)}{\ell(\theta, \theta_0) + y_m^2}\right] \hspace{0.2in} [\text{As } \bbE(\xi)= 0]
\end{align*}
That $X_{i, \theta}$ is symmetric follows from the symmetry of $\xi$. We define $M$ as: 
$$
\max_{1 \le i \le n} \sup_{\theta \in \Omega \times S^{p-1}_m} X_{i, \theta}^2 \lesssim \frac{\max_{1 \le i \le n} \xi_i^2}{n^2 y_m^4} := M
$$
and the wimpy variance: 
\begin{align*}
   \sup_{\theta \in \Omega \times S^{p-1}_m} \sum_{i=1}^n \bbE(X_{i, \theta}^2) & \le \frac{2\sigma_\xi^2}{n} \frac{\ell(\theta_m, \theta_0)}{(\ell(\theta_m, \theta_0) + y_m^2)^2} \\
    & \le \frac{2\sigma_\xi^2}{n} \sup_{x \ge 0}\frac{x}{(x + y_m^2)^2} \\
    & \le \frac{2\sigma_\xi^2}{4ny_m^2} := \sigma_m^2 \,.
\end{align*}
An application of Theorem 11.17  and Theorem 11.1 of \cite{boucheron2013concentration} yields: 
\begin{align}
    \label{eq:var_bound_1} \var(\Gamma_{m, 1}) \le \sigma_m^2 + 64\sqrt{\bbE[M_m]}\bbE[\Gamma_{m, 1}] + 18^2 \bbE[M_m] \,.
\end{align}
Note that we set $V_m = m(\log p)^{1+\delta}$ (which is slightly larger than the VC dimension) and choose $y_m^2$ as: 
$$
y_m^2 = 2\frac{V_m \|\xi\|_{2, n}}{n}\log{\frac{n}{V_m}} = 2\pen(m) \,.
$$
As per Assumption \ref{assm:sparsity_l2}, we confine the model selection in $1 \le m \le (1/4)\lfloor n/(\log{p})^2 \rfloor$. To facilitate the union, we next show that $\sum_{i=1}^{\cM} \var(V_{m, 1}) \to 0$ as $n \to \infty$. We bound each terms on RHS of equation \eqref{eq:var_bound_1}: 
\begin{align*}
   \sum_{i=1}^{\cM}  \sigma_m^2 & =  \frac{\sigma_\xi^2}{4}\sum_{i=1}^{\cM} \frac{1}{n y_m^2} \\
   & = \frac{\sigma_\xi^2}{4\|\xi\|_{n, 2}} \sum_{i=1}^{n/(4(\log p)^2)}  \frac{1}{V_m \log{\frac{n}{V_m }}} \\
   & \le \frac{\sigma_\xi^2}{4 \log{4}\|\xi\|_{n, 2}} \sum_{i=1}^{n/4(\log p)^2} \frac{1}{V_m } \\
   & \le \frac{\sigma_\xi^2}{4 \log{4}(\log{p})^{\delta/2}\|\xi\|_{n, 2}} \sum_{i=1}^{n/(4(\log p)^2)} \frac{1}{m(\log{m})^{1 + \delta/2}} \\
   & \le \frac{\sigma_\xi^2}{4 \log{4}(\log{p})^{\delta/2}\|\xi\|_{n, 2}} \sum_{i=1}^{\infty} \frac{1}{m(\log{m})^{1 + \delta/2}} \\
   & \longrightarrow 0 \text{ as } n \to \infty. 
\end{align*}
Now for the second summand: 
\begin{align*}
    64 \sum_{i=1}^{\cM} \sqrt{\bbE[M_m]}\bbE[\Gamma_{m, 1}] & \le 16 \sqrt{\bbE[M_m]} \\
    & \lesssim 8 \sum_{i=1}^{\cM} \sqrt{ \frac{\bbE[\max_{1 \le i \le n} \xi_i^2]}{n^2 y_m^4} } \\
    & = 8 \sum_{i=1}^{\cM} \frac{\|\xi\|_{n, 2}}{ny_m^2} \\
    & \le \frac{8}{4 \log{4}(\log{p})^{\delta/2}} \sum_{i=1}^{\infty} \frac{1}{m(\log{m})^{1 + \delta/2}} \\
   & \longrightarrow 0 \text{ as } n \to \infty.
\end{align*}
And similarly for the third summand: 
\begin{align*}
    18^2\sum_{i=1}^{\cM}  \bbE[M_m] & = 18^2 \sum_{i=1}^{\cM} \frac{\bbE[\max_{1 \le i \le n} \xi_i^2]}{n^2 y_m^4} \\
    & \le \frac{18^2}{4\log{4}(\log{p})^{\delta}} \sum_{i=1}^{\infty} \frac{1}{m^2(\log{m})^{2 + \delta}} \\
   & \longrightarrow 0 \text{ as } n \to \infty.
\end{align*}
Hence taking $t = 1/8$ and using the fact that $\bbE[\Gamma_{m, 1}] \le 1/8$ for all $m$ for our choice of $\kappa$, we have: 
\begin{equation}
\label{eq:exp_bound_2}
\bbP\left(\Gamma_{\hat m, 1} > \frac14\right) \longrightarrow 0 \,.
\end{equation}
Therefore, combining equation \eqref{eq:exp_bound_1} and \eqref{eq:exp_bound_2} we conclude:  
$$
\bbP\left(\Gamma_{\hat m} > \frac12\right) \longrightarrow 0 
$$
Hence, on its complement event, we have 
\begin{align*}
\ell(\hat \theta_{\hat m}, \theta_0) & \le \frac12 \left(\ell(\hat \theta_{\hat m}, d_0) + y_{\hat m}^2\right) + \pen(s_0) - \pen(\hat m) \\
& = \frac12 \ell(\hat \theta_{\hat m}, d_0)  + \pen(s_0) \,,
\end{align*}
which further implies, 
$$
\ell(\hat \theta_{\hat m}, \theta_0) \le 2\pen(s_0) \,.
$$
This, along with equation \eqref{eq:loss_dist_lower_bound} indicates: 
$$
\dist^2\left((\hat \alpha_{\init}, \hat \beta_{\init}, \hat d), (\alpha_0, \beta_0, d_0)\right) = O_p\left( \frac{s_0 (\log{p})^{(1+\delta)}\|\xi\|_{n, 2}}{n}\left(\log{\frac{n}{s_0\log{p}}}\right)\right) \,.
$$  
The boosting of the rate of $\hat \alpha, \hat \beta$ by replacing $d_0$ by $\hat d$ in the model equation is exactly same as that of Theorem \ref{thm:l2} and hence skipped.

\subsection{Proof of Theorem \ref{thm:mlb_high}}
For the proof of this theorem we follow the techniques of proof of Theorem 2.18 of \cite{mukherjee2019non}. Recall Fano's inequality: if $\Theta \subseteq S^{p-1}$ is a finite $2\eps$ packing set, i.e. for any two $d_i, d_j \in \Theta$, we have $\|d_i - d_j\| \ge 2\eps$ with $|\Theta| < \infty$, then based on n i.i.d. observations $z_1, \dots, z_n$ we have the following minimax lower bound in estimating $d_0$: 
$$
\inf_{\hat d} \sup_{\cP_d} \bbE\left[\left\|\hat d - d\right\|^2\right] \ge \eps^2 \left(1 - \frac{\frac{n}{M^2}\sum_{i, j: d_i, d_j \in \Theta} KL(P_{d_i} || P_{d_j}) + \log{2}}{\log{(|\Theta|-1)}}\right)
$$
Next recall Gilbert-Varshamov Lemma: if $d_H$ is the Hamming distance, i.e. $d_H(x, y) = \sum_{i=1}^d \mathds{1}(x_i \neq y_i)$ with $d$ being the ambient dimension. Then given any $v$ with $1 \le v \le p/8$, we can find $\omega_1, \dots, \omega_M \in \{0, 1\}^p$ which satisfy the following: 
\begin{enumerate}[a)]
\item $d_H(\omega_i, \omega_j) \ge \frac{v}{2} \ \ \forall \ i \neq j \ \ \in \{1, \dots, m\}$. 
\item $\log{M} \ge \frac{v}{8}\log{\left(1 + \frac{d}{2v}\right)}$. 
\item $\|\omega_j\|_0 = v \ \ \forall \ j \ \in \{1, \dots, M\}$. 
\end{enumerate} 
We choose the appropriate $\eps$ later. First, for a fixed $0 < \eps < 1$, we construct the set $\Theta$ as follows: applying Gilbert-Varshamov Lemma in dimension $p-1$ with sparsity $v = s-1$, we choose $\Omega = \left\{\omega_1, \dots, \omega_M\right\} \in \{0, 1\}^{p-1}$ which satisfies the above conditions (a) - (c). Then, for each $\omega_j \in \Omega$ set $d_j$ as: 
$$
d_j = \frac{\left(1, \frac{\eps}{\sqrt{s-1}}\omega_j\right)}{\sqrt{1 + \eps^2}}
$$
It is immediate that $\|d_j\|_2 = 1$ and $\|d_j\|_0 = s$. Set $\Theta = \{d_j: \omega_j \in \Omega\}$. From condition (c) above we have: 
$$
\left| \Theta \right| : = M \ge \frac{s-1}{8}\log{\left(1 + \frac{p-1}{s-1}\right)} \,.
$$
Further note that for any $d_i \neq d_j \in \Theta$: 
$$
\|d_i - d_j\|^2_2 = \frac{\eps^2}{(s-1)(1+\eps^2)}\|\omega_i - \omega_j\|^2 = \frac{\eps^2}{(s-1)(1+\eps^2)}d_H(\omega_i, \omega_j) = \frac{\eps^2}{2(1 + \eps^2)} \ge \frac{\eps^2}{4} \,.
$$
which proves that $\Theta$ is a $\eps/2$ packing set of $S^{p-1}$. On the other hand, from condition (c) above it is immediate that, for any $\omega_i \neq \omega_j \in \Omega$, we have $d_H(\omega_i, \omega_j) \le 2s$. This implies that for any $d_i, d_j \in \Theta$: 
$$
\|d_i - d_j\|_2^2 = \frac{\eps^2}{(s-1)(1+\eps^2)}d_H(\omega_i, \omega_j)  \le 2\eps^2 \,.
$$
Now for each $d_i \in \Theta$ define the distribution $P_{d_i}$ of $(X, Y)$ as: $X \sim \cN(0, I_p), \xi \sim \cN(0, 1)$, $X$ is independent of $\xi$ and: 
$$
Y \overset{d}{=} \mathds{1}_{X^{\top}d_i > 0} + \xi \,.
$$
Hence for any $d_i \neq d_j \in \Theta$, the Kullback-Liebler divergence between $P_{d_i}$ and $P_{d_j}$ is: 
\begin{align*}
KL(P_{d_i} || P_{d_j}) & = \frac12 \bbE_X\left[\left(\mathds{1}_{X^{\top}d_i > 0} - \mathds{1}_{X^{\top}d_j > 0}\right)^2\right] \\
& = \bbP\left(\s(X^{\top}d_i) \neq \s(X^{\top}d_j)\right) \\
& \le C \|d_i - d_j\|_2 \le \eps \ \sqrt{2C^2} \,.
\end{align*}
for some universal constant $C$. Hence applying Fano's inequality we obtain: 
\begin{align*}
\inf_{\hat d} \sup_{\cP_d} \bbE\left[\left\|\hat d - d\right\|^2\right] & \ge 
\frac{\eps^2}{16} \left(1 - \frac{\frac{n}{M^2}\sum_{i, j: d_i, d_j \in \Theta} KL(P_{d_i} || P_{d_j}) + \log{2}}{\log{(|\Theta|-1)}}\right) \\
& \ge \frac{\eps^2}{16} \left(1 - \frac{n\eps \ \sqrt{2C^2} + \log{2}}{\log{\left(\frac{s-1}{8}\log{\left(1 + \frac{p-1}{s-1}\right)}-1\right)}}\right) \,.
\end{align*}
Taking $\eps = (s\log{(1 + p/s)})/n$ we conclude the proof.

\section{Proof of supplementary lemmas}
\label{app:supp_lemmas}

\subsection{Proof of Lemma \ref{lem:huber_lower_bound}}
As $\xi$ has symmetric distribution around origin, without loss of generality we can assume $\mu > 0$. Hence we have to establish the result for $0 < \mu k$. Note that difference $H_k(\xi + \mu) - H_k(\xi)$ can be decomposed into five terms, depending where $\xi$ lies: 
\begin{equation}
\label{eq:huber_diff_1}
H_k(\xi + \mu) - H_k(\xi) = 
\begin{cases}
\frac{1}{2}\left[(\xi + \mu)^2 - \xi^2\right] \,, & \text{ if } -k \le \xi \le k - \mu \\
\frac{1}{2}(\xi + \mu)^2 - k\left(|\xi| - \frac{k}{2}\right) \,, & \text{ if } -k - \mu \le \xi \le -k \\
k\left(|\xi + \mu| - \frac{k}{2}\right) - \frac{\xi^2}{2} \,, & \text{ if } k-\mu \le \xi \le k \\
K\left(|\xi + \mu| - |\xi|\right) \,, & \text{ if } \xi > k \text{ or } \xi < -k - \mu
\end{cases}
\end{equation}
Now we inspect the regions individually. Note that when $-k - \mu \le \xi \le -k$, we have: 
\begin{align*}
H_k(\xi + \mu) - H_k(\xi) & =  \frac{1}{2}(\xi + \mu)^2 - k\left(|\xi| - \frac{k}{2}\right) \\
& = \frac{1}{2}(\xi + \mu)^2 + k\left(\xi + \frac{k}{2}\right)  \hspace{0.2in} [\because |\xi| = -\xi] \\
& = \frac{\xi^2}{2} + (\mu + k)\xi + \frac{\mu^2 + k^2}{2} \\
& = \frac12 \left(\xi + \mu + k\right)^2 - \mu k 
\end{align*}
When $k - \mu \le \xi \le k$: 
\begin{align*}
H_k(\xi + \mu) - H_k(\xi) & = k\left(|\xi + \mu| - \frac{k}{2}\right) - \frac{\xi^2}{2} \\
& = k\left((\xi + \mu) - \frac{k}{2}\right) - \frac{\xi^2}{2} \\
& = - \frac{\xi^2}{2} + k \xi - \frac{k^2}{2} + k \mu   \\
& = -\frac12 \left(\xi - k\right)^2 + k \mu
\end{align*}
Also, we have: 
$$
k\left(|\xi + \mu| - |\xi|\right) = 
\begin{cases}
k\mu \,, & \text{ if } \xi > k \\
-k\mu \,, & \text{ if } \xi < -k - \mu \,. 
\end{cases}
$$
Hence we can modify equation \eqref{eq:huber_diff_1} as: 
\begin{equation}
\label{eq:huber_diff_2}
H_k(\xi + \mu) - H_k(\xi) = 
\begin{cases}
\mu \xi + \frac{\mu^2}{2} \,, & \text{ if } -k \le \xi \le k - \mu \\
\frac12 \left(\xi + \mu + k\right)^2 - \mu k  \,, & \text{ if } -k - \mu \le \xi \le -k \\
-\frac12 \left(\xi - k\right)^2 + k \mu \,, & \text{ if } k-\mu \le \xi \le k \\
k\mu \,, & \text{ if } \xi > k \\
-k\mu \,, & \text{ if } \xi < -k - \mu \,.
\end{cases}
\end{equation}
Note that the term $-\mu k$ is active on the region $\xi \le -k$ and $\mu k$ is active on the region $\xi \ge k - \mu$. From the symmetry of the distribution of $\xi$, this effect of $-\mu k$ and $\mu k$ on the region $(-\infty, -k)$ and $(k, \infty)$ will cancel each other upon taking expectation and the effect of $\mu k$ on $(k-\mu, k)$ will remain. Hence we have: 
\begin{align*}
\bbE\left[H_k(\xi + \mu) - H_k(\xi) \right] & = \frac{\mu^2}{2}\bbP\left( -k \le \xi \le k - \mu\right) + \mu \bbE\left[\xi \mathds{1}_{ -k \le \xi \le k - \mu}\right] \\
& \qquad + \bbE\left[\frac12 \left(\xi + \mu + k\right)^2 \mathds{1}_{-k - \mu \le \xi \le -k}\right] + \bbE\left[-\frac12 \left(\xi - k\right)^2 \mathds{1}_{k - \mu \le \xi \le k}\right] \\
& = \frac{\mu^2}{2}\bbP\left( -k \le \xi \le k - \mu\right) - \mu \bbE\left[\xi \mathds{1}_{ k - \mu \le \xi \le k}\right] \\
& \qquad + \bbE\left[\frac12 \left(\xi + \mu + k\right)^2 \mathds{1}_{-k - \mu \le \xi \le -k}\right] + \bbE\left[\left(\mu k -\frac12 \left(\xi - k\right)^2\right) \mathds{1}_{k - \mu \le \xi \le k}\right]  \\
& = \frac{\mu^2}{2}\bbP\left( -k \le \xi \le k - \mu\right)  + \bbE\left[\left(\mu k -\mu \xi - \frac12 \left(\xi - k\right)^2\right) \mathds{1}_{k - \mu \le \xi \le k}\right]  \\
& \qquad \qquad \qquad \qquad + \bbE\left[\frac12 \left(\xi + \mu + k\right)^2 \mathds{1}_{-k - \mu \le \xi \le -k}\right] \\
%& = \frac{\mu^2}{2}\bbP\left( -k \le \xi \le k\right)  + \bbE\left[\left(-\frac12 \mu^2  - \mu(\xi - k) - \frac12 \left(\xi - k\right)^2\right) \mathds{1}_{k - \mu \le \xi \le k}\right]  \\
%& \qquad \qquad \qquad \qquad + \bbE\left[\frac12 \left(\xi + \mu + k\right)^2 \mathds{1}_{-k - \mu \le \xi \le -k}\right] \\
%& = \frac{\mu^2}{2}\bbP\left( -k \le \xi \le k\right)  -\frac{1}{2} \bbE\left[\left(\left(\xi -k + \mu\right)^2\right) \mathds{1}_{k - \mu \le \xi \le k}\right]  \\
%& \qquad \qquad \qquad \qquad + \bbE\left[\frac12 \left(\xi + \mu + k\right)^2 \mathds{1}_{-k - \mu \le \xi \le -k}\right] \\
%& \ge  \frac{\mu^2}{2}\bbP\left( -k \le \xi \le k - \mu\right)  + \bbE\left[\left(\mu (k -\xi) - \frac12 \left(\xi - k\right)^2\right) \mathds{1}_{k - \mu \le \xi \le k}\right]  \\
& \ge  \frac{\mu^2}{2}\bbP\left( -k \le \xi \le k - \mu\right)  + \bbE\left[\left(\mu k -\mu \xi - \frac12 \left(\xi - k\right)^2\right) \mathds{1}_{k - \mu \le \xi \le k}\right]  \\
& \ge  \frac{\mu^2}{2}\bbP\left( -k \le \xi \le k - \mu\right)  
\end{align*}
where the last inequality follows from the fact: 
$$
f(\xi) = \mu k -\mu \xi - \frac12 \left(\xi - k\right)^2 \ge 0 \ \ \ \forall \ \ \ \xi \in \left[k-\mu, k\right] \,.
$$
observing the fact that: 
$$
\bbE\left[\tilde H_k(\xi + \mu) - \tilde H_k(\xi) \right]= \frac{k+1}{k}\bbE\left[H_k(\xi + \mu) - H_k(\xi) \right] \ge \bbE\left[H_k(\xi + \mu) - H_k(\xi) \right]
$$
we complete the proof for all $k > 0$. Now for $k = 0$ for $0 < \mu < \delta$, 
\begin{align*}
    \bbE\left[|\xi + \mu| - |\xi|\right] & = -\mu\bbP(\xi \le -\mu) + \mu \bbP(\xi > 0) + \bbE\left[(2\xi + \mu) \mathds{1}_{-\mu \le \xi \le 0}\right] \\\\
    & = \mu \bbP(0 \le \xi \le \mu) + \bbE\left[(-2\xi + \mu) \mathds{1}_{0 \le \xi \le \mu}\right] \\\\
    & = \bbE\left[(-2\xi + 2\mu) \mathds{1}_{0 \le \xi \le \mu}\right]  \\
    & = 2 \int_{0}^{\mu}(\mu - x) f_\xi(x) \ dx \\
    & \ge f_\xi(0) \int_{0}^{\mu}(\mu - x) \ dx \\
    & = \frac{\mu^2}{2} f_{\xi}(0) \,.
\end{align*}
This completes the proof.

\subsection{Proof of Lemma \ref{lem:lower_bound_rw}}
\begin{proof}
Although we assume continuous steps, our proof can be certainly extended to the case when $S_n$ takes value $0$ with positive probability. The proof critically uses Theorem 4 of Chapter 12 of Volume 2 of \cite{feller1957introduction}. To keep the notational similarity with the book, define: 
$$
q_n = \bbP\left(\max_{1 \le i \le n}S_i < 0\right)
$$
and the corresponding generating function $q(s)$ as: 
$$
q(s) = 1 + \sum_{n = 1}^{\infty}s^n q_n \,.
$$
Then from equation (7.22) of Theorem 4, Chapter 12, Vol. 2 of \cite{feller1957introduction} we have: 
\begin{equation}
    \label{eq:relation_q}
    \log{q(s)} = \sum_{n=1}^{\infty}\frac{s^n}{n}\bbP\left(S_n < 0\right) := f(s) \iff q(s) = e^{f(s)}\,.
\end{equation}
Now we need a lower bound on $q_n$. Note that from the property of the generating function we have: 
$$
q_n = n! q^{(n)}(0) \,.
$$
On the other hand from Faa di Bruno's formula:  
\begin{align}
    q^{(n)}(0) = \left.\frac{d^n}{ds^n}e^{f(s)}\right\vert_{s=0} & =  \left.\left[e^{f(s)}\sum \frac{n!}{m_1!1!^{m_2}m_1!2!^{m_2}\cdots m_n!n!^{m_n}} \Pi_{j=1}^n \left(f^{(j)}(s)\right)^{m_j}\right]\right\vert_{s=0} \notag \\
    & = \left. \left[e^{f(s)}\sum \frac{n!}{m_1!m_2!\cdots m_n!} \Pi_{j=1}^n \left(\frac{f^{(j)}(s)}{j!}\right)^{m_j}\right]\right\vert_{s= 0} \notag \\
    \label{eq:fdb} & = \sum \frac{n!}{m_1!m_2!\cdots m_n!} \Pi_{j=1}^n \left(\frac{f^{(j)}(0)}{j!}\right)^{m_j} 
\end{align}
where the sum runs over all the sequences $\{m_j\}_{j=1}^n$ such that: 
$$
\sum_{i=1}^n i m_i = n \,.
$$
Note that $f^{(j)}(0)$ is non-negative for $j$. Hence using only one sequence with $m_1 = \dots = m_{n-1} = 0$ and $m_n = 1$, equation \eqref{eq:fdb} can lower bounded as: 
\begin{align*}
     \left.\frac{d^n}{ds^n}e^{f(s)}\right\vert_{s=0} \ge f^{(n)}(0)\,.
\end{align*}
On the from the expression of $f(s)$ from equation \eqref{eq:relation_q} it is immediate that: 
$$
f^{(n)}(0) = n! \frac{1}{n}\bbP\left(S_n < 0\right) \,.
$$
Combining our findings we have: 
\begin{align*}
    n! q_n = q^{(n)}(0) & = \left.\frac{d^n}{ds^n}e^{f(s)}\right\vert_{s=0} \\
    & \ge f^{(n)}(0) \\
    & = n! \frac{1}{n}\bbP\left(S_n < 0\right) \,.
\end{align*}
This immediately implies $q_n \ge \bbP\left(S_n > 0\right)/n$ which completes our proof. 
\end{proof}

\subsection{Proof of Lemma \ref{lem:lower_random_onesided}}
From the definition of distribution of $\xi$ we have: 
$$
F_{\xi}(t) =
\frac{\frac12 + t^{\gamma}}{1 + t^{\gamma}} \ \ \ \ t \ge 0 \,.
$$
and 
$$
F_{\xi}(-t) = 1 - F_{\xi}(t) =  \frac{1}{2(1 + t^{\gamma})} \,.
$$
Hence it is immediate that for any $t_0 > 0$: 
\begin{equation}
\label{eq:bound_on_F}
\frac{1}{2(1 + t_0^{-\gamma})} \le \sup_{t \ge t_0} t^{\gamma}\bar F_{\xi}(t) = \sup_{t \ge t_0} \frac{t^{\gamma}}{2(1 + t^{\gamma})} \le \frac12
\end{equation}
For any fixed $k \ge 1$:  
\allowdisplaybreaks
\begin{align}
    \bbP\left(M \ge k\right) & = \sum_{j \ge k} \bbP\left(M = j\right) \notag \\
    & = \sum_{j \ge k} \bbP\left(S_i > S_j \ \forall \ 0 \le i \le j-1, S_i > S_j \ \forall \ i \ge j+1\right) \notag \\
    & = \bbP\left(S_1 > 0, S_2>0, \dots\right)  \sum_{j \ge k} \bbP\left(\max_{1 \le i \le j} S_i < 0 \right) \notag \\
    & = p^* \sum_{j \ge k} \bbP\left(\max_{1 \le i \le j} S_j < 0 \right) \hspace{0.2in} [p^* = P\left(S_1 > 0, S_2>0, \dots\right) ] \notag \\
    \label{eq:eq1} & \ge  p^* \sum_{j \ge k} \frac{1}{j}P\left(S_j \le 0\right) 
\end{align}
where the last inequality uses Lemma \ref{lem:lower_bound_rw}. From the symmetry of the distribution of $\xi$ we have: 
\begin{align*}
    \bbP\left(S_j \le 0\right) & = \bbP\left(\sum_{i=1}^j \xi_j \le - j\mu \right) = \bbP\left(\sum_{i=1}^j \xi_j > j\mu \right) \,.
\end{align*}
Set $a_j = j^{1/\gamma}$. Define the event $A_i$ as: 
$$
A_i = \left\{\xi_i > j\mu + (j-1)a_j, \ \xi_l \in [-a_j, j\mu) \ \forall \ 1 \le l \neq i \le j\right\}
$$
Clearly $\{A_i\}'s$ are disjoint events and 
\begin{align}
\bbP\left(\sum_{i=1}^j \xi_j > j\mu \right) & \ge \bbP\left(\cup_{i=1}^j A_i \right) \notag \\
& = \sum_{i=1}^j \bbP(A_i) \notag \\
& = j \bar F\left(j\mu + (j-1)a_j \right) \left(F[-a_j, j\mu)\right)^{j-1} \notag \\
%& = j \bar F(j\mu) \times \left[\left(1 - \bar F(a_j) - \bar F(j\mu)\right)^{j-1} \frac{\bar F\left(j\mu + (j-1)a_j \right)}{\bar F(j\mu)}\right] \notag \\
%& = \frac{1}{2(1 + \mu^{-\alpha})} \times \frac{1}{(j\mu)^{\alpha}} \times j \times \left[\left(1 - \bar F(a_j) - \bar F(j\mu)\right)^{j-1} \frac{\bar F\left(j\mu + (j-1)a_j \right)}{\bar F(j\mu)}\right] \notag \\
%\label{eq:1} & :=  \frac{c_1}{(j\mu)^{\alpha}} \times j \times \left[\left(1 - \bar F(a_j) - \bar F(j\mu)\right)^{j-1}  \frac{\bar F\left(j\mu + (j-1)a_j \right)}{\bar F(j\mu)}\right]  
\label{eq:1}  &  = j \bar F\left(j\mu + (j-1)a_j \right)\left(1 - \bar F(a_j) - \bar F(j\mu)\right)^{j-1} 
\end{align}
Next note that, $j\mu + (j-1)a_j \ge \mu$ for all $j \ge 1$. Therefore from equation \eqref{eq:bound_on_F} we have for all $j \ge 1$: 
$$
\bar F\left(j\mu + (j-1)a_j \right) \ge \frac{\left(j\mu + (j-1)a_j \right)^{-\gamma}}{2(1 + \mu^{-\gamma})} \,,
$$
which further implies: 
\begin{align*}
    j \times \bar F\left(j\mu + (j-1)a_j \right) & \ge \frac{1}{{2(1 + \mu^{-\alpha})}}\frac{j}{\left(j\mu + (j-1)a_j \right)^{\gamma}} \\
    & = \frac{1}{{2(1 + \mu^{-\gamma})}}\frac{j}{\left(j\mu + (j-1)j^{1/\alpha} \right)^{\gamma}} \\
    & =  \frac{1}{{2(1 + \mu^{-\gamma})}}\frac{j}{j^{\gamma + 1}\left(j^{-1/\gamma}\mu + \left(1 - \frac{1}{j}\right) \right)^{\gamma}} \\
    & = \frac{1}{j^{\gamma}} \ \frac{1}{{2(1 + \mu^{-\gamma})}}\ \frac{1}{\left(j^{-1/\gamma}\mu + \left(1 - \frac{1}{j}\right) \right)^{\gamma}} \\
    & \ge \frac{1}{j^{\gamma}} \ \frac{1}{{2(1 + \mu^{-\gamma})}}\ \frac{1}{(\mu + 1)^{\gamma}} := \frac{c_1}{j^{\gamma}} \,.
  \end{align*}
%
%\begin{align*}
%    j \times  \frac{\bar F\left(j\mu + (j-1)a_j \right)}{\bar F(j\mu)} & \ge j \times \frac{2c_1 \left(j\mu + (j-1)a_j \right)^{-\alpha}}{(j\mu)^{-\alpha}} \\
%    & = j \times 2c_1 \times \left(1 + \frac{1}{\mu}\frac{j-1}{j}a_j\right)^{-\alpha} \\
%    & \ge j \times 2c_1 \times \left(1 + \frac{1}{2\mu}a_j\right)^{-\alpha} \\
%    & = j \times 2c_1 \times \left(1 + \frac{1}{2\mu}j^{\frac{1}{\alpha}}\right)^{-\alpha} \\
%    & =  2c_1 \times \left(\frac{1}{2\mu} + j^{-\frac{1}{\alpha}}\right)^{\alpha} \ge \frac{2c_1}{(2\mu)^{\alpha}} := c_2 
%\end{align*}
Next observe that $(j\mu)^{\gamma} \ge j$ for all $j \ge 1$ if $\mu > 1$ or for all $j \ge \mu^{-\gamma/(\gamma - 1)}$ if $\mu \le 1$. Using this in equation \eqref{eq:1} we have for all $j \ge 1 \vee \lceil \mu^{-\gamma/(\gamma - 1)}\rceil$: 
\begin{align}
    \bbP\left(\sum_{i=1}^j \xi_j > j\mu \right) & \ge \frac{c_1}{j^{\gamma}} \left[\left(1 - \bar F(a_j) - \bar F(j\mu)\right)^{j-1} \right] \notag \\
    & \ge  \frac{c_1}{j^{\gamma}} \times \left(1 -\frac{1}{2(1+a_j^{\alpha})} - \frac{1}{2(1 + (j\mu)^\alpha)}\right)^{j-1} \notag \\
    & = \frac{c_1}{j^{\gamma}} \times \left(1 -\frac{1}{2(1+j)} - \frac{1}{2(1 + (j\mu)^\alpha)}\right)^{j-1} \notag \\
     & \ge  \frac{c_1}{j^{\gamma}} \times \left(1 -\frac{1}{(1+j)}\right)^{j-1} \notag \\
 & \ge   \frac{c_1}{j^{\gamma}}  \times \inf_{x \ge 1} \left(1 -\frac{1}{(1+x)}\right)^{x-1} \notag \\
           \label{eq:final_bound_11}  & :=  \frac{c_1c_2}{j^{\gamma}} 
\end{align}
%Denoting by $c_3 = 1 \wedge \inf_{x \ge 1} \left(1 -\frac{1}{(1+x)}\right)^{x-1}$ we obtain from equation \eqref{eq:final_bound_1} and \eqref{eq:final_bound_11} that for $j \ge 1 \vee \mu^{-\alpha/(\alpha - 1)}$: 
%$$
%P\left(\sum_{i=1}^j \xi_j > j\mu \right) \ge  \frac{c_1c_2c_3}{(j\mu)^{\alpha}} \,.
%$$
Using this in equation \eqref{eq:eq1} we obtain: 
\begin{align*}
    P\left(M \ge k\right) & \ge  p^* \sum_{j \ge k} \frac{1}{j}P\left(S_j \le 0\right) \\
    & = p^* \sum_{j \ge k} \frac{1}{j}\left(\sum_{i=1}^j \xi_j > j\mu \right) \\
    & \ge c_1c_2 \times p^* \times \sum_{j \ge k} j^{-(\gamma + 1)} \\
    &  \ge c_1c_2 \times p^* \times \int_{k}^{\infty} x^{-(\gamma + 1)} \ dx \\
    & \ge c_1c_2 \times p^* \times \frac{1}{\gamma k^{\gamma}} \,.
\end{align*}
This completes the proof of lower bound.

\subsection{Proof of Lemma \ref{lem:two_sided_rw_minimizer}}
We have, by symmetry: 
$$
\bbP\left(|M_{ts}| > k\right) = \bbP\left(M_{ts} > k\right) + \bbP\left(M_{ts} < -k\right) = 2\bbP\left(M_{ts} > k\right) \,.
$$
Hence, by virtue of Lemma \ref{lem:lower_random_onesided}, all we need to show is: 
$$
P(M_{ts} = k) \ge p^*\bbP(M_{os} = k) \,.
$$
Towards that end: 
\begin{align*}
    \bbP\left(M_{ts} = k\right) & = \bbP\left(S_K \le S_i \ \forall \ 0 \le i \le k-1, S_k \le S_i \ \forall \ i \ge k+1, S_k \le \inf_{j \ge 1}S_{-j}\right) \\
    & = \bbP\left(S_K \le S_i \ \forall \ 0 \le i \le k-1, S_k \le \inf_{j \ge 1}S_{-j}\right)\bbP\left(S_i \ge 0 \ \forall \ i \ge 1\right) \\
    & = p^* \bbP\left(S_K \le S_i \ \forall \ 0 \le i \le k-1, S_k \le \inf_{j \ge 1}S_{-j}\right) \\
    & \ge p^* \bbP\left(S_K \le S_i \ \forall \ 0 \le i \le k-1, S_k \le \inf_{j \ge 1}S_{-j} \mid \inf_{j \ge 1}S_{-j} > 0\right)\bbP\left(\inf_{j \ge 1}S_{-j} > 0\right) \\
    & = p^* \bbP\left(\inf_{j \ge 1}S_{-j} > 0\right)  \bbP\left(S_K \le S_i \ \forall \ 0 \le i \le k-1\right) \\
    & = \bbP\left(\inf_{j \ge 1}S_{-j} > 0\right)\bbP\left(M_{os} = k\right) \\
    & = p^* \bbP\left(M_{os} = k\right) \,.
\end{align*}
where the last equality follows from the fact: 
$$
\bbP\left(M_{os} = k\right) = p^* \bbP\left(S_K \le S_i \ \forall \ 0 \le i \le k-1\right)
$$
This completes the proof.

\subsection{Proof of Lemma \ref{cor:tsCPP}}
%As the compound Poisson process is symmetric on both positive and negative real axis we have: 
%$$
%\bbP\left(|M_{ts, CPP}| > x\right) = \bbP\left(M_{ts, CPP} > x\right) + \bbP\left(M_{ts, CPP} \le -x\right) = 2 \bbP\left(M_{ts, CPP} > x\right)
%$$
Using same line of arguments as in  Corollary 1:
\begin{align*}
    & \bbP\left(M_{ts, CPP} = x\right) \\
    & = \sum_{k=0}^{\infty}\bbP\left(M_{ts, CPP} = x \mid N_1(x) = k\right)\bbP(N_1(x) = k) \\
    & = \sum_{k=0}^{\infty}\bbP\left(M_{ts} = k \right)\bbP(N_1(x) = k)  \\
    & \ge p^* \sum_{k=0}^{\infty}\bbP\left(M_{os} = k \right)\bbP(N_1(x) = k) \\
    & = \bbP\left(M_{os, CPP} = x\right)
\end{align*}
Hence to establish Corollary \ref{cor:tsCPP}, all we need show: 
$$
\bbP\left(M_{os, CPP} > x\right) \ge \frac{c_0}{2f_X^{\gamma}(d_0)}x^{-\gamma}
$$
for all large $x$, where $M$ is the argmin of one sided compound Poisson process, namely the minimizer of the following: 
$$
X_+(t) = \sum_{i=1}^{N_1(t)}X_i\,, \ \ \ t \in \bbR^+ \,.
$$
Now we have: 
\allowdisplaybreaks
\begin{align*}
    P(M_{os, CPP} > x) & = \sum_{k = 0}^{\infty} \bbP\left(M_{os, CPP} > x \mid N_1(x) = k\right) \bbP\left(N_1(x) = k\right) \\
    & = \sum_{k = 0}^{\infty} \bbP\left(\argmin_{i \ge 0} S_i > k\right) \bbP\left(N_1(x) = k\right) \\
    & \ge \sum_{k = k_0}^{\infty} \bbP\left(\argmin_{i \ge 0} S_i > k\right) \bbP\left(N_1(x) = k\right) \\
     & \ge c_0 \sum_{k = k_0}^{\infty} k^{-\gamma}\bbP\left(N_1(x) = k\right)  \\
     & = c_0 \sum_{k = k_0}^{\infty} \frac{e^{-\Lambda(x)}\Lambda(x)^k}{k!k^{\gamma}} \\
     & \ge c_0 \sum_{k = k_0}^{\infty} \frac{e^{-\Lambda(x)}\Lambda(x)^k}{k!(k+1)(k+2) \dots (k+\gamma)} \\
     & = c_0 \sum_{k = k_0}^{\infty} \frac{e^{-\Lambda(x)}\Lambda(x)^k}{(k+\gamma)!} \\
     & = c_0\Lambda(x)^{-\gamma}\sum_{k = k_0 + \alpha}^{\infty} \frac{e^{-\Lambda(x)}\Lambda(x)^k}{k!} \\
     & =  c_0\Lambda(x)^{-\gamma}\bbP\left(N_1(x) \ge k_0 + \alpha\right) \\
     & \ge \frac{c_0}{2}\Lambda(x)^{-\gamma} = \frac{c_0}{2f_X^{\gamma}(d_0)}x^{-\gamma}
\end{align*}
where the last inequality is valid as long as $\med\left(N_1(X)\right) \ge k_0 + \gamma$. From \cite{adell2005median}, we know as $N_1(x) \sim Poisson(xf_X(d_0))$, we have $\med\left(N_1(x)\right) \ge xf_X(d_0) - \log{2}$. Hence the inequality is valid as long as $x \ge (k_0 + \gamma + \log{2})/f_X(d_0)$. From This completes the proof.

\subsection{Proof of Lemma \ref{lem:finite_sample_tail_bound}}
\begin{proof}
As per our model description, all the parallel change point processes are i.i.d. Therefore $n(\hat d_i - d_{0, i})$ has same distribution across $1 \le i \le m$. Therefore, we henceforth define $F_n$ to be the distribution of $n(\hat d - d_0)$ and drop $i$ from subscript. From the definition of change point estimator, we have: 
\begin{align*}
n(\hat d - d_0) & = \sargmin_{t} \sum_{i=1}^n \left(\xi_i + \frac12\right)\left\{\mathds{1}_{d_0 < X_i \le d_0 + \frac{t}{n}}\right\} \\
& \qquad \qquad \qquad + \sargmin_{t} \sum_{i=1}^n \left(-\xi_i + \frac12\right)\left\{\mathds{1}_{d_0 + \frac{t}{n} < X_i \le d_0}\right\} \\
& = \sargmin_{t} \sum_{i=1}^{N_{n, +}(t)} \left(\xi_i + \frac12\right)\mathds{1}_{t \ge 0} + \sum_{i=1}^{N_{n, -}(t)} \left(-\xi_i + \frac12\right)\mathds{1}_{t < 0} \,.
\end{align*}
Here the count processes $N_{n, +}(t)$ and $N_{n,-}(t)$ are defined as follows: For $t \ge 0$, 
$$
N_{n,+}(t) = \sum_{i=0}^n \mathds{1}_{d_0 \le X_i \le d_0 + \frac{t}{n}}\sim \text{Bin}\left(n, F_X\left(d_0 + \frac{t}{n}\right) - F_X(t)\right)
$$
and for $t < 0$,
$$
N_{n,-}(t) =  \sum_{i=0}^n \mathds{1}_{d_0 + \frac{t}{n} \le X_i \le d_0} \sim \text{Bin}\left(n, F_X(t) - F_X\left(d_0 + \frac{t}{n}\right)\right) \,.
$$
These processes can be though as finite sample approximation of Compound Poisson Process, where we have approximated the Poisson random variables by Binomial random variables. It is immediate that: 
\begin{align*}
    N_{n,+}(t) & \overset{\mathscr{L}}{\implies} \text{Pois}\left(tf(d_0)\right) \,,\\
    N_{n,-}(t) & \overset{\mathscr{L}}{\implies} \text{Pois}\left(-tf(d_0)\right) \,.
\end{align*}
We name the process as \emph{compound Binomial process} and henceforth denote by CBP: 
\begin{equation}
\label{eq:cbp}
CBP(t) = \sum_{i=1}^{N_{n, +}(t)} \left(\xi_i + \frac12\right)\mathds{1}_{t \ge 0} + \sum_{i=1}^{N_{n, -}(t)} \left(-\xi_i + \frac12\right)\mathds{1}_{t < 0}
\end{equation}
Hence we work with the smallest argmin instead of mid-argmin just for some technical simplicity, but all of the following analysis is valid for mid-argmin also. As will be evident later, the thickness of the tail of the distribution $F_n$ (the distribution of $n(\hat d - d_0)$) is closely related to the tail of the minimizer of a random walk with finitely many steps. Therefore, we start by establishing a lower bound on the tail of a n-step random walk $\{S_i\}_{i=0}^n$ with the usual convention $S_0 = 0$ and step distribution $X_i \overset{d}{=} \xi_i + 1/2$. Let $Z_n$ be the minimizer of this random walk. The random variable $Z_n$ is supported on $\{0, 1, \dots, n\}$. Then for any $0 \le k \le n-1$: 
\begin{align}
    \bbP\left(Z_n > k\right) & = \sum_{j = k+1}^n \bbP\left(Z_n = j\right) \notag \\
    & = \sum_{j = k+1}^n\bbP\left(S_i > S_j \ \forall \ 0 \le i \le j-1, S_i > S_j \ \forall \  j+1 \le i \le n\right) \notag \\
    & = \sum_{j = k+1}^n P\left(S_i < 0 \ \forall \ 1 \le i \le j \right) \bbP\left(S_i > 0 \ \forall \  1 \le i \le n - j\right) \notag \\
    & \ge \bbP\left(S_i > 0 \ \forall \ 1 \le i < \infty \right)  \sum_{j = k+1}^n  \bbP\left(\max_{1 \le i \le j} S_i < 0 \right) \notag \\
    & = p^* \sum_{j = k+1}^n \bbP\left(\max_{1 \le i \le j} S_j < 0 \right) \hspace{0.2in} \left[p^* = \bbP\left(\min_{1 \le i < \infty} S_i > 0\right)\right] \notag \\
    \label{eq:eq1_finite-sample} & \ge  p^* \sum_{j = k+1}^n \frac{1}{j}\bbP\left(S_j \le 0\right) 
\end{align}
From equation \eqref{eq:final_bound_11} in the proof of Lemma \ref{lem:lower_random_onesided} we conclude for $j \ge 2^{\gamma/(\gamma - 1)} := k_0$: 
$$
\bbP\left(S_j \le 0\right)  \ge \frac{c_1c_2}{j^{\gamma}} \,.
$$
Using the above bound in equation \eqref{eq:eq1_finite-sample} we conclude: 
\begin{align}
     \bbP\left(Z_n > k\right) & \ge p^* \sum_{j = k+1}^n \frac{1}{j}\bbP\left(S_j \le 0\right) \notag \\
     & \ge  p^* \sum_{j = k+1}^n \frac{1}{j} \frac{c_1c_2}{j^{\gamma}} \notag \\
     & = c_1c_2 p^* \sum_{j = k+1}^n  \frac{1}{j^{\gamma + 1}} \notag \\
     & \ge  c_1c_2 p^*  \int_{k+1}^{n+1} x^{-(\gamma + 1)} \ dx \hspace{0.2in} [\text{Riemann integral lower bound}] \notag \\
     \label{eq:finite_sample_rw_lower_bound} & = \frac{c_1c_2p^*}{\gamma}\left[\frac{1}{(k+1)^{\gamma}} - \frac{1}{(n+1)^{\gamma}}\right]
\end{align}
Now we go back to the random variable of interest $n(\hat d - d_0)$. Let $X_{(i)}$ denotes the $i^{th}$ order statistics of $\{X_i\}_{i \le n}$. If $X_{(i)} < d_0 < X_{(i+1)}$, then from the definition of $n(\hat d - d_0)$, we have a random walk with $i$ steps on the negative axis and a random walk with $n-i$ steps on the positive axis. Therefore the number of steps of random walk on either side of origin is equal to the number of $X_i's$ on the corresponding side of $d_0$. Denote by $R_n$ (and respectively $L_n$), the number of $X_i's$ greater than $d_0$ (respectively less than $d_0$). Hence $R_n \sim \text{Bin}(n, \bar F_X(d_0))$ and $L_n \sim \text{Bin}(n, F_X(d_0))$ with $R_n + L_n = n$. Then we have for any $x > 0$: 
\begin{align*}
    & \bbP\left(n(\hat d - d_0) > x\right) \notag \\
    & = \sum_{r = 0}^n \sum_{k = 0}^{r}\bbP\left(n(\hat d - d_0) > x \mid R_n = r, N_{n, +}(x) = k\right)\bbP(N_{n, +}(x) = k \mid R_n = r) \bbP(R_n = r) 
    \end{align*}
Given $R_n = r$, we have a two sided random walk, with $r$ steps on the positive real line $n-r$ steps on the negative real line. Therefore, the event $n(\hat d - d_0) > x$ given $N_{n, +}(x) = k$ and $R_n = r$ is equivalent to the event that in a two sided random walks with $r$ steps on the right and $n-r$ steps on the left, the argmin is on the right and it happens after $k$ steps. More precisely, if we denote by $S_0 \equiv 0, S_1, \dots, S_r$ to be the random walk on the right side with step distribution $(\xi+ 1/2)$ and $S'_0 \equiv 0, S'_1, \dots, S'_{n-r}$ to be random walk on the left with step distribution $(-\xi + 1/2)$, then the above event corresponds that this two sided random walk is minimized at $S_j$ for some $k+1 \le j \le r$. Therefore we write: 
\begin{align}
    & \bbP\left(n(\hat d - d_0) > x\right) \notag \\
    & = \sum_{r = 0}^n \sum_{k = 0}^{r}\bbP\left(n(\hat d - d_0) > x \mid R_n = r, N_{n, +}(x) = k\right)\bbP(N_{n, +}(x) = k \mid R_n = r) \bbP(R_n = r)  \notag\\
        \label{eq:finite_lower_1}   & = \sum_{r = 0}^n \sum_{k = 0}^{r}\left[\bbP\left(\argmin \text{ of twosided RW } > k \right)  \times \bbP(N_{n, +}(x) = k \mid R_n = r) \bbP(R_n = r)\right]
%    & = \sum_{r = 0}^n \sum_{k = 0}^{r}\left[\bbP\left(\argmin \text{ twosided RW } > k \right)  \times \bbP(N_{n, +}(x) = k \mid R_n = r) \bbP(R_n = r)\right]
\end{align}
Next, we obtain a lower bound on the tail of the minimizer of the two sided random walk. Note that, we have already established a lower bound on the tail of the minimizer of a one sided random walk in equation \eqref{eq:finite_sample_rw_lower_bound}, which we exploit here to get a lower bound on the tail of the minimizer of this two-sided incarnation: 
\begin{align*}
    & \bbP\left(\argmin \text{ twosided RW } > k \right) \\
    & = \sum_{j = k+1}^r  \bbP\left(\argmin \text{ twosided RW } = j \right) \\
    & = \sum_{j = k+1}^r  \bbP\left(S_j < S_0, \dots, S_j < S_{j-1}, S_j < S_{j+1}, \dots, S_j < S_r, S_j < \min_{1 \le i \le n-r} S'_i\right) \\
    & = \sum_{j = k+1}^r  \bbP\left(S_j < S_0, \dots, S_j < S_{j-1}, S_j < \min_{1 \le i \le n-r} S'_i\right)\bbP(S_1 > 0, \dots, S_{r-j} > 0)  \\
    & \ge \sum_{j = k+1}^r  \bbP\left(S_j < S_0, \dots, S_j < S_{j-1}, S_j < \min_{1 \le i \le n-r} S'_i \mid \min_{1 \le i \le n-r} S'_i > 0\right) \times \\
    & \hspace{15em} \bbP\left(\min_{1 \le i \le n-r} S'_i > 0\right) \bbP(S_1 > 0, \dots, S_{r-j} > 0) \\
    & \ge p^* \sum_{j = k+1}^r  \bbP\left(S_j < S_0, \dots, S_j < S_{j-1}, S_j < \min_{1 \le i \le n-r} S'_i \mid \min_{1 \le i \le n-r} S'_i > 0\right) \times \\
    & \hspace{25em} \bbP(S_1 > 0, \dots, S_{r-j} > 0) \\
    &= p^* \sum_{j = k+1}^r  \bbP\left(S_j < S_0, \dots, S_j < S_{j-1}\right) \bbP(S_1 > 0, \dots, S_{r-j} > 0) \\\\
    & = p^* \bbP\left(\argmin \text{ one-sided RW with length }r > k\right) \\\\
    & \ge \frac{c_1c_2(p^*)^2}{\gamma} \left[\frac{1}{(k+1)^{\gamma}} - \frac{1}{(r+1)^{\gamma}}\right] \hspace{0.2in} [\text{From equation }\eqref{eq:finite_sample_rw_lower_bound}]\,.
\end{align*}
For the rest of the calculation we assume $\gamma$ (the number of finite moments of the error distribution $\xi$) is an integer, as all our calculation is valid by replacing $\gamma$ by $\lfloor \gamma \rfloor$. Define a success probability $p_{x, n}$ as: 
$$
p_{x, n} = \bbP\left(X \in \left(d_0 + \frac{x}{n}, d_0\right) \mid X > d_0\right) = \frac{F_X\left(d_0 + \frac{x}{n}, d_0\right) - F_X(d_0)}{1 - F_X(d_0)}\,.
$$
Therefore it is immediate that: 
$$
N_{n, +}(x) \mid R_n = r \sim \text{Bin}\left(r, p_{x, n}\right) \,.
$$
As per our assumption $F_X$ has continuous density $f_X$ with $f_X(d_0) > 0$. Therefore, there exists $\delta_1 > 0$ such that $f_X(t) > f_X(d_0)/2$ for $|t - d_0| \le \delta_1$. Hence, for any $0 \le x \le n \delta_1$, we have:
\begin{equation}
\label{eq:pxn_lb}
p_{x, n} \ge \frac{x}{n} \times \frac{f_X(d_0)}{2(1 - F_X(d_0))} \,.
\end{equation}
On the other, let $f_{\max}$ be the upper bound on $f_X$ on the entire $\bbR$. Then, again from the mean value theorem, we have: 
\begin{equation}
\label{eq:pxn_ub}
p_{x, n} \le \frac{x}{n} \times \frac{f_{\max}}{2(1 - F_X(d_0))} \,.
\end{equation}
Therefore combining equations \eqref{eq:pxn_lb} and \eqref{eq:pxn_ub}, we have: 
\begin{equation}
\label{eq:pxn_lb_ub}
\frac{x}{n} \times \frac{f_X(d_0)}{2(1 - F_X(d_0))} \le p_{x, n} \le \frac{x}{n} \times \frac{f_{\max}}{2(1 - F_X(d_0))} \,.
\end{equation}
We will use the above relations in our rest of the calculation. Going back to equation \eqref{eq:finite_lower_1} we have: 
\begin{align}
    & \bbP\left(n(\hat d - d_0) > x\right) \notag \\ 
    & = \sum_{r = 1}^n \sum_{k = 0}^{r}\left[\bbP\left(\argmin \text{ twosided RW } > k\right) \times \bbP(N_{n, +}(x) = k \mid R_n = r) \bbP(R_n = r)\right] \notag \\ \notag \\
     & = \sum_{r = k_0}^n \sum_{k = k_0}^{r}\left[\bbP\left(\argmin \text{ twosided RW } > k\right) \times \bbP(N_{n, +}(x) = k \mid R_n = r) \bbP(R_n = r)\right] \notag \\ \notag \\
    & \ge  \frac{c_1c_2(p^*)^2}{\gamma} \sum_{r = 1}^n \sum_{k = 0}^{r}\left[\frac{1}{(k+1)^{\gamma}} - \frac{1}{(r+1)^{\gamma}}\right]\bbP(N_{n, +}(x) = k \mid R_n = r) \bbP(R_n = r) \notag \\ \notag \\
    & = \frac{c_1c_2(p^*)^2}{\gamma} \left[\sum_{r = k_0}^n \sum_{k = k_0}^{r}\frac{1}{(k+1)^{\gamma}}\bbP(N_{n, +}(x) = k \mid R_n = r) \bbP(R_n = r)  \right. \notag \\
    & \hspace{15em} \left. -  \sum_{r = k_0}^n \frac{1}{(r+1)^{\gamma}}\bbP(k_0 \le N_{n, +}(x) < r \mid R_n = r) \bbP(R_n = r)\right] \notag \\
    & =\frac{c_1c_2(p^*)^2}{\gamma} \left[\sum_{r = k_0}^n \sum_{k = k_0}^{r}\frac{1}{(k+1)^{\gamma}}\bbP(N_{n, +}(x) = k \mid R_n = r) \bbP(R_n = r)  \right. \notag \\
    & \hspace{15em} \left. -  \sum_{r = 1}^n \frac{1}{(r+1)^{\gamma}}\bbP(R_n = r)\right] \notag \\
    & = \frac{c_1c_2(p^*)^2}{\gamma} \left[\sum_{r = k_0}^n \sum_{k = k_0}^{r}\frac{1}{(k+1)^{\gamma}}\dbinom{r}{k}p_{x, n}^k(1-p_{x, n})^{r-k} \bbP(R_n = r)  \right.\notag  \\
    & \hspace{15em} \left. -  \sum_{r = 1}^n \frac{1}{(r+1)^{\gamma}}\bbP(R_n = r)\right] \notag \\
    \label{eq:finite_lower_2} & \ge \frac{c_1c_2(p^*)^2}{\gamma} \left[\sum_{r = k_0}^n \left\{\sum_{k = k_0}^{r}\frac{1}{(k+1)^{\gamma}}\dbinom{r}{k}p_{x, n}^k(1-p_{x, n})^{r-k}\right\} \bbP(R_n = r)  -  \sum_{r = 1}^n \frac{1}{(r+1)^{\gamma}}\bbP(R_n = r)\right]
\end{align}
where the last inequality is obtained by replacing $1 - p_{x, n}$ by its upper bound 1. The inner sum of the above equation can be analyzed as follows: 
\begin{align*}
    & \sum_{k = k_0}^{r}\frac{1}{(k+1)^{\gamma}}\dbinom{r}{k}p_{x, n}^k(1-p_{x, n})^{r-k} \\
    & = \sum_{k = k_0}^{r}\frac{1}{(k+1)^{\gamma}}\frac{r!}{k!(r-k)!}p_{x, n}^k(1-p_{x, n})^{r-k} \\
    & \ge \sum_{k = k_0}^{r}\frac{1}{(k+1)(k+2)\dots(k+\gamma)}\frac{r!}{k!(r-k)!}p_{x, n}^k(1-p_{x, n})^{r-k} \\
    & \ge \frac{p_{x, n}^{-\gamma}}{(r+1)(r+2)\dots(r+\gamma)}\sum_{k = k_0}^{r}\frac{(r+\gamma)!}{(k+\gamma)!(r-k)!}p_{x, n}^{k+\gamma}(1-p_{x, n})^{r-k}  \\
    & = \frac{p_{x, n}^{-\gamma}}{(r+1)(r+2)\dots(r+\gamma)} \bbP\left(\text{Bin}(r+\gamma, p_{x, n}) \ge k_0 + \gamma\right) \,.
\end{align*}
Putting this back into equation \eqref{eq:finite_lower_2} we obtain:
\begin{align}
     & \bbP\left(n(\hat d - d_0) > x\right) \notag \\ \notag \\
%     & \ge \frac{2^{\gamma}c_1c_2c_3 (p^*)^2}{\gamma} \left[\sum_{r = 1}^n \left\{\sum_{k = 0}^{r}\frac{1}{(k+1)^{\gamma}}\dbinom{r}{k}p_{x, n}^k(1-p_{x, n})^{r-k}\right\} \bbP(R_n = r)  -  \sum_{r = 1}^n \frac{1}{(r+1)^{\gamma}}\bbP(R_n = r)\right] \notag \\ \notag \\
     \label{eq:finite_lower_3} & \ge \frac{c_1c_2(p^*)^2}{\gamma} p_{x, n}^{-\gamma} \left[\sum_{r = k_0}^n\frac{\bbP\left(\text{Bin}(r+\gamma, p_{x, n}) \ge k_0 + \gamma\right)}{(r+1)(r+2)\dots(r+\gamma)}\bbP(R_n = r) - p_{x, n}^{\gamma}\sum_{r = 1}^n \frac{1}{(r+1)^{\gamma}}\bbP(R_n = r) \right]
\end{align}
Now from the properties of the inverse moments of the binomial distribution (see e.g. \cite{cribari2000note}) we have: 
$$
\sum_{r = 1}^n \frac{1}{(r+1)^{\gamma}}\bbP(R_n = r) \le C\left(n\bar F_X(d_0)\right)^{-\gamma} \hspace{0.2in} [\text{Recall }R_n \sim \text{Bin}(n, \bar F_X(d_0))] \,.
$$
On the other hand we have for all $r \ge n \bar F_X(d_0)$: 
\begin{align*}
(r + \gamma)p_{x, n} & \ge (n \bar F_X(d_0) + \gamma) \times  \frac{x}{n} \times \frac{f_X(d_0)}{2(1 - F_X(d_0))} \hspace{0.2in} [\text{Equation }\eqref{eq:pxn_lb}] \\
& \ge \left(\bar F_X(d_0) + \frac{\gamma}{n}\right) \times x \times \frac{f_X(d_0)}{2(1 - F_X(d_0))} \\
& \ge  \bar F_X(d_0) \times x \times \frac{f_X(d_0)}{2(1 - F_X(d_0))} > \gamma + k_0
\end{align*}
for all $x > 2(\gamma + k_0)/f_X(d_0)$. Now we know that the median of $Bin(n, p)$ is $\lfloor np \rfloor$ or $\lceil np \rceil$. For simplicity, we will use the bound here $\bbP(Bin(n, p) \ge np) \ge 1/2$ as it will be valid simply replacing $np$ by $\lfloor np \rfloor$ and this will not alter any of our subsequent analysis. Therefore we have: 
\begin{align*}
    & \sum_{r = k_0}^n\frac{\bbP\left(\text{Bin}(r+\gamma, p_{x, n}) \ge \gamma + k_0\right)}{(r+1)(r+2)\dots(r+\gamma)}\bbP(R_n = r) \\
    & \ge \sum_{r = n\bar F(d_0)}^n\frac{\bbP\left(\text{Bin}(r+\gamma, p_{x, n}) \ge \gamma + k_0\right)}{(r+1)(r+2)\dots(r+\gamma)}\bbP(R_n = r) \\
    & \ge  \sum_{r = n\bar F(d_0)}^n\frac{1}{2(r+1)(r+2)\dots(r+\gamma)}\bbP(R_n = r) \\
    & \ge \frac{1}{2(n+1)(n+2)\dots(n+\gamma)}\sum_{r = n\bar F(d_0)}^n\bbP(R_n = r) \\
    & \ge \frac{1}{4(n+1)(n+2)\dots(n+\gamma)} \,.
\end{align*}
Going back to equation \eqref{eq:finite_lower_3} we have: 
\begin{align*}
    & \bbP\left(n(\hat d - d_0) > x\right) \\
    & \ge \frac{c_1c_2(p^*)^2}{\gamma} p_{x, n}^{-\gamma} \left[\sum_{r = k_0}^n\frac{\bbP\left(\text{Bin}(r+\gamma, p_{x, n}) \ge \gamma\right)}{(r+1)(r+2)\dots(r+\gamma)}\bbP(R_n = r) - p_{x, n}^{\gamma} \sum_{r = 1}^n \frac{1}{(r+1)^{\gamma}}\bbP(R_n = r) \right] \\
    & \ge \frac{c_1c_2(p^*)^2}{\gamma} p_{x, n}^{-\gamma} \left[\frac{1}{4(n+1)(n+2)\dots(n+\gamma)} - p_{x, n}^{\gamma}\frac{C(\bar F(d_0))^{-\gamma}}{n^{\gamma}} \right] \\
    & \ge \frac{c_1c_2(p^*)^2}{\gamma} x^{\gamma} \times \left(\frac{f_{X, \max}}{1 - F_X(d_0)}\right)^{-\gamma} \left[\frac{n^{\gamma}}{4(n+1)(n+2)\dots(n+\gamma)} - (np_{x, n})^{\gamma}\frac{C(\bar F(d_0))^{-\gamma}}{n^{\gamma}} \right] \\
      & \ge \frac{c_1c_2(p^*)^2}{\gamma} x^{\gamma} \times \left(\frac{f_{X, \max}}{1 - F_X(d_0)}\right)^{-\gamma} \left[\frac{1}{4(2^{\gamma})}  - (np_{x, n})^{\gamma}\frac{C(\bar F(d_0))^{-\gamma}}{n^{\gamma}} \right] \\
      & \ge \frac{c_1c_2(p^*)^2}{\gamma} x^{\gamma} \times \left(\frac{f_{X, \max}}{1 - F_X(d_0)}\right)^{-\gamma} \left[\frac{1}{4(2^{\gamma})}  - \left(x \times \frac{f_{X, \max}}{1 - F_X(d_0)}\right)^{\gamma}\frac{C(\bar F(d_0))^{-\gamma}}{n^{\gamma}} \right] \\
      & \ge \frac{c_1c_2(p^*)^2}{\gamma} x^{\gamma} \times \left(\frac{f_{X, \max}}{1 - F_X(d_0)}\right)^{-\gamma} \left[\frac{1}{4(2^{\gamma})}  - \left( \frac{\delta_2 f_{X, \max}}{C(1 - F_X(d_0))^2}\right)^{\gamma}  \right] \hspace{0.2in} [\forall \ x \le \delta_2 n ]\\
      & \ge \frac{c_1c_2(p^*)^2}{\gamma}  \times \left(\frac{f_{X, \max}}{1 - F_X(d_0)}\right)^{-\gamma} \times \frac{1}{2^{\gamma + 3}} \times x^{-\gamma} \\
      & = \frac{c_1c_2(p^*)^2}{\gamma 2^{\gamma + 3}}  \times \left(\frac{f_{X, \max}}{1 - F_X(d_0)}\right)^{-\gamma} \times x^{-\gamma} \,.
%    & \ge \frac{2^{\gamma}c_1c_2c_3 (p^*)^2}{\gamma} \left(\frac{x}{f_-}\right)^{-\gamma} \left[\frac{1}{42^{\gamma}} - \left(\frac{f_+a_n}{n}\right)^{\gamma}C(\bar F(d_0))^{-\gamma} \right] \\
%    & \ge  \frac{f_-^{\gamma}c_1c_2c_3 (p^*)^2}{8\gamma}x^{-\gamma} \,.
\end{align*}
where the last inequality is valid for small enough $\delta_2$, i.e. we choose $\delta_2$ which satisfies: 
$$
\frac{1}{2^{\gamma + 2}} - \left(\frac{\delta_2 f_{X, \max}}{C(1 - F_X(d_0))^2}\right)^{\gamma} \ge \frac{1}{2^{\gamma + 3}} \,.
$$
Therefore we have established that for any $2\gamma/f_X(d_0) \le x \le (\delta_1 \wedge \delta_2)n$:  
\begin{equation}
\label{eq:pos_x_bound}
\bbP\left(n(\hat d - d_0) > x\right) \le \frac{c_1c_2(p^*)^2}{\gamma 2^{\gamma + 3}}  \times \left(\frac{f_{X, \max}}{1 - F_X(d_0)}\right)^{-\gamma} \times x^{-\gamma} \,.
\end{equation}
The calculation for the negative $x$ is similar. As introduced before, $L_n$ denotes the number of $X_i's$ on the left of $d_0$ and $L_n \sim \text{Bin}(n, F(d_0))$. Given $L_n = l$, define $S_0 \equiv 0, S'_1, \dots, S'_l$ to be random walk on the left of origin and $S_0 = 0, S_1, S_2, \dots, S_{n-l}$ on the right of origin. Given $L_n = l, N_{n, -}(x) = k'$, the event $n(\hat d - d_0) < -x$ is equivalent to the event that in a two sided random walk with $l$ steps on the left and $n-l$ steps on the right, the minima occurs on the left and it occurs at one of the steps among $\left\{S'_{k'+1}, S'_{k'+2}, \dots, S'_{l}\right\}$. Using the similar logic as above we obtain for any $2\gamma/f_X(d_0) \le x \le (\delta_1 \wedge \delta_2)n$: 
\begin{equation}
\label{eq:neg_x_bound}
\bbP\left(n(\hat d - d_0) < -x\right) \ge \frac{c_1c_2(p^*)^2}{\gamma 2^{\gamma + 3}}  \times \left(\frac{f_{X, \max}}{1 - F_X(d_0)}\right)^{-\gamma} \times x^{-\gamma} \,.
\end{equation}
Finally, from equation \eqref{eq:pos_x_bound} and \eqref{eq:neg_x_bound} we conclude for any $2\gamma/f_X(d_0) \le x \le (\delta_1 \wedge \delta_2)n$:
$$
\bbP\left(\left|n(\hat d - d_0)\right| > x\right) \ge \frac{c_1c_2(p^*)^2}{\gamma 2^{\gamma + 2}}  \times \left(\frac{f_{X, \max}}{1 - F_X(d_0)}\right)^{-\gamma} \times x^{-\gamma} \,.
$$
%\begin{align}
%    & \bbP\left(n(\hat d - d_0) < -x\right) \notag \\
%    & = \sum_{l = 0}^n \sum_{k' = 0}^{r}\bbP\left(n(\hat d - d_0) < -x \mid L_n = l, N_{n, -}(x) = k'\right)\bbP(N_{n, +}(x) = k' \mid L_n = l) \bbP(L_n = l)  \notag\\
%    & = \sum_{l = 0}^n \sum_{k' = 0}^{r}\left[\bbP\left(\argmin \text{ twosided RW } > k' \mid L_n = l, N_{n, -}(x) = k'\right) \right. \notag \\
%    & \qquad \qquad \qquad \qquad \qquad \qquad \qquad \qquad \left. \times \bbP(N_{n, -}(x) = k \mid L_n = l) \bbP(L_n = l)\right]\notag  \\
%    & = \sum_{l = 0}^n \sum_{k' = 0}^{r}\left[\bbP\left(\argmin \text{ twosided RW } > k' \mid L_n = l\right) \right. \notag \\
%    \label{eq:finite_lower_11} & \qquad \qquad \qquad \qquad \qquad \qquad \qquad \qquad \left. \times \bbP(N_{n, -}(x) = k' \mid L_n = l) \bbP(L_n = l)\right]
%\end{align}
%The rest of the calculation is similar, interchanging the role of $S_i$ and $S_i'$ and hence skipped for brevity. Therefore, we obtain: 
%$$
%\bbP\left(\left| n(\hat d - d_0) \right| > x\right) \ge \frac{f_-^{\gamma}c_1c_2c_3 (p^*)^2}{4\gamma}x^{-\gamma}
%$$
%for all $x \ge \frac{\gamma}{f_- \left(F(d_0) \wedge \bar F(d_0)\right)}$ and $x \le a_n$ for any sequence of $a_n = o(n)$. 
This completes the proof. 
\end{proof}

\subsection{Proof of Lemma \ref{lem:finite_sample_tail_bound_upper}}
We use the same notations as used in the proof of Lemma \ref{lem:finite_sample_tail_bound}.  As $\xi_i$'s are bounded by $b$ with mean $\mu > 0$, by applying Hoeffding's inequality, we have for any $j  \in \bbN$: 
$$
\bbP\left(\bar \xi_j < -\mu\right) \le e^{-\frac{j\mu^2}{4b^2}} := e^{-cj} \,,
$$
with $c = \mu^2/4b^2$. As in Lemma \ref{lem:finite_sample_tail_bound}, we start with establishing an upper bound on the tail on the minimizer of random walk. Let $\{S_j\}_{j =0,\dots, n}$ denotes a $n$-step random walk and let $Z_n$ denotes its minimizer supported on $\{0, 1, \dots, n\}$. We then have: 
\begin{align*}
\bbP\left(Z_n > k\right) & = \sum_{j=k+1}^n \bbP\left(Z_n = j\right) \\
& =  \sum_{j=k+1}^n \bbP\left(S_j < 0, \dots, S_j < S_{j-1}, S_j < S_{j+1}, \dots, S_j < S_n\right) \\
& \le \sum_{j=k+1}^n  \bbP\left(S_j < 0\right) \\
& = \sum_{j=k+1}^n  \bbP\left(\bar \xi_j < - \mu \right) \\
& \le  \sum_{j=k+1}^n e^{-cj} = \frac{e^{-c(k+1)}}{1 - e^{-c}} \,.
\end{align*}
Going back to distribution of $n(\hat d - d_{0})$, as in the proof of Lemma \ref{lem:finite_sample_tail_bound} we have: 
\begin{align}
    & \bbP\left(n(\hat d - d_0) > x\right) \notag \\
    & = \sum_{r = 0}^n \sum_{k = 0}^{r}\bbP\left(n(\hat d - d_0) > x \mid R_n = r, N_{n, +}(x) = k\right)\bbP(N_{n, +}(x) = k \mid R_n = r) \bbP(R_n = r) \notag \\
    & = \sum_{r = 0}^n \sum_{k = 0}^{r}\left[\bbP\left(\argmin \text{ twosided RW } > k \mid R_n = r, N_{n, +}(x) = k\right) \right. \notag \\
    & \qquad \qquad \qquad \qquad \qquad \qquad \qquad \qquad \left. \times \bbP(N_{n, +}(x) = k \mid R_n = r) \bbP(R_n = r)\right]\notag  \\
    & = \sum_{r = 0}^n \sum_{k = 0}^{r}\left[\bbP\left(\argmin \text{ twosided RW } > k \mid R_n = r\right) \right. \notag \\
    \label{eq:finite_upper_1} & \qquad \qquad \qquad \qquad \qquad \qquad \qquad \qquad \left. \times \bbP(N_{n, +}(x) = k \mid R_n = r) \bbP(R_n = r)\right]
    \end{align}
Upper bounding the argmin of two-sided random walk is relatively easier: 
\begin{align*}
& \bbP\left(\argmin \text{ twosided RW } > k \mid R_n = r\right) \\
& = \sum_{j = k+1}^r \bbP\left(\argmin \text{ twosided RW }=j \mid R_n = r\right) \\
& = \sum_{j = k+1}^r \bbP\left(S_j < \min_{0 \le i \le j-1}S_i, S_j < \min_{j+1 \le i \le r}S_i, S_j < \min_{1 \le i \le n-r} S'_i \right) \\
& \le \sum_{j = k+1}^r \bbP\left(S_j < 0\right) \le \frac{e^{-c(k+1)}}{1 - e^{-c}} \,.
\end{align*} 
Using this bound in equation \eqref{eq:finite_upper_1} we obtain for any $0 \le x \le n \delta_1$ (where $\delta_1$ is same as defined in the proof of Lemma \ref{lem:finite_sample_tail_bound}, i.e. we choose $\delta_1 > 0$ such that $f_X(t) \ge f_X(d_0)/2$ for all $|t - d_0| \le \delta_1$): 
\begin{align*}
  & \bbP\left(n(\hat d - d_0) > x\right) \notag \\ 
  & = \sum_{r = 0}^n \sum_{k = 0}^{r}\left[\bbP\left(\argmin \text{ twosided RW } > k \mid R_n = r\right) \right. \notag \\
    & \qquad \qquad \qquad \qquad \qquad \qquad \qquad \qquad \left. \times \bbP(N_{n, +}(x) = k \mid R_n = r) \bbP(R_n = r)\right] \\
    & \le  \sum_{r = 0}^n \sum_{k = 0}^{r} \frac{e^{-c(k+1)}}{1 - e^{-c}} \bbP(N_{n, +}(x) = k \mid R_n = r) \bbP(R_n = r) \\
    & \le \frac{e^{-c}}{1 - e^{-c}} \sum_{r = 0}^n \sum_{k = 0}^{r} e^{-ck} \ \bbP(N_{n, +}(x) = k \mid R_n = r) \bbP(R_n = r) \\
    & = \frac{e^{-c}}{1 - e^{-c}} \sum_{r = 0}^n \left(1 - p_{n, x} + p_{n, x}e^{-c}\right)^r \bbP(R_n = r) \\
    & = \frac{e^{-c}}{1 - e^{-c}} \left(1 - \bar F(d_0) + \bar F(d_0)\left(1 - p_{n, x} + p_{n, x}e^{-c}\right)\right)^n \\
    & =  \frac{e^{-c}}{1 - e^{-c}} \left(1 - \bar F(d_0)p_{n, x}(1 - e^{-c})\right)^n \\
    & \le \frac{e^{-c}}{1 - e^{-c}}   \left(1 - (1 - e^{-c})\frac{xf_X(d_0)}{2n}\right)^n \\
    & \le \frac{e^{-c}}{1 - e^{-c}}  e^{-x\frac{f_X(d_0)}{2}(1 - e^{-c})} \,.
\end{align*}
The calculation for $\bbP(n(\hat d - d_0) < -x)$ for $x > 0$ is similar and hence skipped for brevity. Therefore we obtain for $0 \le |x| \le n\delta_1$: 
$$
\bbP\left(\left|n(\hat d - d_0)\right| > x\right)  \le \frac{2e^{-c}}{1 - e^{-c}}  e^{-x\frac{f_X(d_0)}{2}(1 - e^{-c})} \,.
$$ 
This completes the proof.

\subsection{Proof of Proposition \ref{prop:multiplier_ineq}}
From proposition 5 of \cite{han2019convergence} we have: 
$$
\bbE\left\|\sum_{i=1}^n \xi f(X_i)\right\|_\cF \le \bbE\left[\sum_{k=1}^n\left(\left|\eta_{(k)}\right| - \left|\eta_{(k+1)}\right|\right)\bbE\left\|\sum_{i=1}^k \eps_i f(X_i)\right\|_\cF\right]
$$
where $\left|\eta_{(1)}\right| \ge \left|\eta_{(2)}\right| \ge \dots \left|\eta_{(n)}\right| \ge \left|\eta_{(n+1)}\right| = 0$ are the decreasing order statistics of $\left\{\left|\xi_i - \xi_i'\right|\right\}_{i=1}^n$ where $\left\{\xi_i'\right\}_{1 \le i \le n}$ are i.i.d copy of $\left\{\xi_i\right\}_{1 \le i \le n}$. Hence we have: 
\begin{align*}
\bbE\left\|\sum_{i=1}^n \xi f(X_i)\right\|_\cF & \le \bbE\left[\sum_{k=1}^n\left(\left|\eta_{(k)}\right| - \left|\eta_{(k+1)}\right|\right)\bbE\left\|\sum_{i=1}^k \eps_i f(X_i)\right\|_\cF\right] \\
& \le \bbE\left[\sum_{k=1}^n\left(\left|\eta_{(k)}\right| - \left|\eta_{(k+1)}\right|\right)\left(\varphi_n(k) + b_n\right)\right] \\
& = \bbE\left[\sum_{i=1}^n \int_{\left|\eta_{(k+1)}\right|}^{\left|\eta_{(k)}\right|} \varphi_n(k) \ dt\right] + b_n \bbE\left[\left|\eta_{(1)}\right|\right]  \\
& \le \bbE\left[\int_0^{\infty} \varphi_n\left(\left|\left\{i: \left|\eta_i\right| \ge t\right\}\right|\right) \ dt\right] + b_n \bbE\left[\max_{1 \le i \le n} \left|\xi_i - \xi_i'\right|\right] \\
& \le \int_0^{\infty} \varphi_n\left(\sum_{i=1}^n \bbP\left(|\xi_i - \xi_i'| > t\right)\right) + 2b_n \bbE\left[\max_{1 \le i \le n} \left|\xi_i\right|\right] \hspace{0.2in} [\text{By Jensen's inequality}] \\
& \le \int_0^{\infty} \varphi_n\left(\sum_{i=1}^n \left(\bbP\left(|\xi_i| > t/2\right) + \bbP\left(|\xi'_i| > t/2\right)\right)\right) + 2b_n \bbE\left[\max_{1 \le i \le n} \left|\xi_i\right|\right] \\
& = \int_0^{\infty} \varphi_n\left(2\sum_{i=1}^n \bbP\left(|\xi_i| > t/2\right)\right) + 2b_n \bbE\left[\max_{1 \le i \le n} \left|\xi_i\right|\right] \\
& = 2\int_0^{\infty} \varphi_n\left(2\sum_{i=1}^n \bbP\left(|\xi_i| > t\right)\right) + 2b_n \bbE\left[\max_{1 \le i \le n} \left|\xi_i\right|\right]  \\
& \le 4\int_0^{\infty} \varphi_n\left(\sum_{i=1}^n \bbP\left(|\xi_i| > t\right)\right) + 2b_n \bbE\left[\max_{1 \le i \le n} \left|\xi_i\right|\right]
\end{align*}
where the last inequality follows from the fact that $\varphi_n(0) = 0$ and $\varphi_n$ concave which leads to $\varphi_n(2x) \le 2 \varphi_n(x)$.

\bibliography{AG, AGMV, mybib_2, mybib, D&C&CP}{}
\bibliographystyle{plain}

\end{document}